\def\fmtname{AmS-TeX}

\def\fmtversion{2.2}
\catcode`\@=11
\ifx\amstexloaded@\relax\catcode`\@=\active
  \endinput\else\let\amstexloaded@\relax\fi
\newlinechar=`\^^J
\def\W@{\immediate\write\sixt@@n}
\def\CR@{\W@{^^J\fmtname - Version \fmtversion^^J}}
\CR@ \everyjob{\CR@}
\message{Loading definitions for}
\message{misc utility macros,}
\toksdef\toks@@=2
\long\def\rightappend@#1\to#2{\toks@{\\{#1}}\toks@@
 =\expandafter{#2}\xdef#2{\the\toks@@\the\toks@}\toks@{}\toks@@{}}
\def\alloclist@{}
\newif\ifalloc@
\def\showallocations{{\def\\{\immediate\write\m@ne}\alloclist@}\alloc@true}
\def\alloc@#1#2#3#4#5{\global\advance\count1#1by\@ne
 \ch@ck#1#4#2\allocationnumber=\count1#1
 \global#3#5=\allocationnumber
 \edef\next@{\string#5=\string#2\the\allocationnumber}%
 \expandafter\rightappend@\next@\to\alloclist@}
\newcount\count@@
\newcount\count@@@
\def\FN@{\futurelet\next}
\def\DN@{\def\next@}
\def\DNii@{\def\nextii@}
\def\RIfM@{\relax\ifmmode}
\def\RIfMIfI@{\relax\ifmmode\ifinner}
\def\setboxz@h{\setbox\z@\hbox}
\def\wdz@{\wd\z@}
\def\boxz@{\box\z@}
\def\setbox@ne{\setbox\@ne}
\def\wd@ne{\wd\@ne}
\def\iterate{\body\expandafter\iterate\else\fi}
\def\err@#1{\errmessage{AmS-TeX error: #1}}
\newhelp\defaulthelp@{Sorry, I already gave what help I could...^^J
Maybe you should try asking a human?^^J
An error might have occurred before I noticed any problems.^^J
``If all else fails, read the instructions.''}
\def\Err@{\errhelp\defaulthelp@\err@}
\def\eat@#1{}
\def\in@#1#2{\def\in@@##1#1##2##3\in@@{\ifx\in@##2\in@false\else\in@true\fi}%
 \in@@#2#1\in@\in@@}
\newif\ifin@
\def\space@.{\futurelet\space@\relax}
\space@. %
\newhelp\athelp@
{Only certain combinations beginning with @ make sense to me.^^J
Perhaps you wanted \string\@\space for a printed @?^^J
I've ignored the character or group after @.}
{\catcode`\~=\active 
 \lccode`\~=`\@ \lowercase{\gdef~{\FN@\at@}}}
\def\at@{\let\next@\at@@
 \ifcat\noexpand\next a\else\ifcat\noexpand\next0\else
 \ifcat\noexpand\next\relax\else
   \let\next\at@@@\fi\fi\fi
 \next@}
\def\at@@#1{\expandafter
 \ifx\csname\space @\string#1\endcsname\relax
  \expandafter\at@@@ \else
  \csname\space @\string#1\expandafter\endcsname\fi}
\def\at@@@#1{\errhelp\athelp@ \err@{\Invalid@@ @}}
\def\atdef@#1{\expandafter\def\csname\space @\string#1\endcsname}
\newhelp\defahelp@{If you typed \string\define\space cs instead of
\string\define\string\cs\space^^J
I've substituted an inaccessible control sequence so that your^^J
definition will be completed without mixing me up too badly.^^J
If you typed \string\define{\string\cs} the inaccessible control sequence^^J
was defined to be \string\cs, and the rest of your^^J
definition appears as input.}
\newhelp\defbhelp@{I've ignored your definition, because it might^^J
conflict with other uses that are important to me.}
\def\define{\FN@\define@}
\def\define@{\ifcat\noexpand\next\relax
 \expandafter\define@@\else\errhelp\defahelp@                               
 \err@{\string\define\space must be followed by a control
 sequence}\expandafter\def\expandafter\nextii@\fi}                          
\def\undefined@@@@@@@@@@{}
\def\preloaded@@@@@@@@@@{}
\def\next@@@@@@@@@@{}
\def\define@@#1{\ifx#1\relax\errhelp\defbhelp@                              
 \err@{\string#1\space is already defined}\DN@{\DNii@}\else
 \expandafter\ifx\csname\expandafter\eat@\string                            
 #1@@@@@@@@@@\endcsname\undefined@@@@@@@@@@\errhelp\defbhelp@
 \err@{\string#1\space can't be defined}\DN@{\DNii@}\else
 \expandafter\ifx\csname\expandafter\eat@\string#1\endcsname\relax          
 \global\let#1\undefined\DN@{\def#1}\else\errhelp\defbhelp@
 \err@{\string#1\space is already defined}\DN@{\DNii@}\fi
 \fi\fi\next@}

\def\predefine#1#2{\let#1#2}
\def\undefine#1{\let#1\undefined}
\message{page layout,}
\newdimen\captionwidth@
\captionwidth@\hsize
\advance\captionwidth@-1.5in
\def\pagewidth#1{\hsize#1\relax
 \captionwidth@\hsize\advance\captionwidth@-1.5in}
\def\pageheight#1{\vsize#1\relax}
\def\hcorrection#1{\advance\hoffset#1\relax}
\def\vcorrection#1{\advance\voffset#1\relax}
\message{accents/punctuation,}

\let\graveaccent\`
\let\acuteaccent\'
\let\tildeaccent\~
\let\hataccent\^
\let\underscore\_
\let\B\=
\let\D\.
\let\ic@\/
\def\/{\unskip\ic@}
\def\textfonti{\the\textfont\@ne}
\def\t#1#2{{\edef\next@{\the\font}\textfonti\accent"7F \next@#1#2}}
\def~{\unskip\nobreak\ \ignorespaces}
\def\.{.\spacefactor\@m}
\atdef@;{\leavevmode\null;}
\atdef@:{\leavevmode\null:}
\atdef@?{\leavevmode\null?}
\edef\@{\string @}
\def\relaxnext@{\let\next\relax}
\atdef@-{\relaxnext@\leavevmode
 \DN@{\ifx\next-\DN@-{\FN@\nextii@}\else
  \DN@{\leavevmode\hbox{-}}\fi\next@}%
 \DNii@{\ifx\next-\DN@-{\leavevmode\hbox{---}}\else
  \DN@{\leavevmode\hbox{--}}\fi\next@}%
 \FN@\next@}
\def\srdr@{\kern.16667em}
\def\drsr@{\kern.02778em}
\def\sldl@{\drsr@}
\def\dlsl@{\srdr@}
\atdef@"{\unskip\relaxnext@
 \DN@{\ifx\next\space@\DN@. {\FN@\nextii@}\else
  \DN@.{\FN@\nextii@}\fi\next@.}%
 \DNii@{\ifx\next`\DN@`{\FN@\nextiii@}\else
  \ifx\next\lq\DN@\lq{\FN@\nextiii@}\else
  \DN@####1{\FN@\nextiv@}\fi\fi\next@}%
 \def\nextiii@{\ifx\next`\DN@`{\sldl@``}\else\ifx\next\lq
  \DN@\lq{\sldl@``}\else\DN@{\dlsl@`}\fi\fi\next@}%
 \def\nextiv@{\ifx\next'\DN@'{\srdr@''}\else
  \ifx\next\rq\DN@\rq{\srdr@''}\else\DN@{\drsr@'}\fi\fi\next@}%
 \FN@\next@}

\def\textfontii{\the\textfont\tw@}
\def\lbrace@{\delimiter"4266308 }
\def\rbrace@{\delimiter"5267309 }
\def\{{\RIfM@\lbrace@\else{\textfontii f}\spacefactor\@m\fi}
\def\}{\RIfM@\rbrace@\else
 \let\@sf\empty\ifhmode\edef\@sf{\spacefactor\the\spacefactor}\fi
 {\textfontii g}\@sf\relax\fi}
\let\lbrace\{
\let\rbrace\}
\def\AmSTeX{{\textfontii A\kern-.1667em%
  \lower.5ex\hbox{M}\kern-.125emS}-\TeX\spacefactor1000 }
\message{line and page breaks,}
\def\vmodeerr@#1{\Err@{\string#1\space not allowed between paragraphs}}
\def\mathmodeerr@#1{\Err@{\string#1\space not allowed in math mode}}
\def\linebreak{\RIfM@\mathmodeerr@\linebreak\else
 \ifhmode\unskip\unkern\break\else\vmodeerr@\linebreak\fi\fi}

\newskip\saveskip@
\def\allowlinebreak{\RIfM@\mathmodeerr@\allowlinebreak\else
 \ifhmode\saveskip@\lastskip\unskip
 \allowbreak\ifdim\saveskip@>\z@\hskip\saveskip@\fi
 \else\vmodeerr@\allowlinebreak\fi\fi}
\def\nolinebreak{\RIfM@\mathmodeerr@\nolinebreak\else
 \ifhmode\saveskip@\lastskip\unskip
 \nobreak\ifdim\saveskip@>\z@\hskip\saveskip@\fi
 \else\vmodeerr@\nolinebreak\fi\fi}
\def\newline{\relaxnext@
 \DN@{\RIfM@\expandafter\mathmodeerr@\expandafter\newline\else
  \ifhmode\ifx\next\par\else
  \expandafter\unskip\expandafter\null\expandafter\hfill\expandafter\break\fi
  \else
  \expandafter\vmodeerr@\expandafter\newline\fi\fi}%
 \FN@\next@}
\def\dmatherr@#1{\Err@{\string#1\space not allowed in display math mode}}
\def\nondmatherr@#1{\Err@{\string#1\space not allowed in non-display math
 mode}}
\def\onlydmatherr@#1{\Err@{\string#1\space allowed only in display math mode}}
\def\nonmatherr@#1{\Err@{\string#1\space allowed only in math mode}}
\def\mathbreak{\RIfMIfI@\break\else
 \dmatherr@\mathbreak\fi\else\nonmatherr@\mathbreak\fi}
\def\nomathbreak{\RIfMIfI@\nobreak\else
 \dmatherr@\nomathbreak\fi\else\nonmatherr@\nomathbreak\fi}
\def\allowmathbreak{\RIfMIfI@\allowbreak\else
 \dmatherr@\allowmathbreak\fi\else\nonmatherr@\allowmathbreak\fi}
\def\pagebreak{\RIfM@
 \ifinner\nondmatherr@\pagebreak\else\postdisplaypenalty-\@M\fi
 \else\ifvmode\removelastskip\break\else\vadjust{\break}\fi\fi}
\def\nopagebreak{\RIfM@
 \ifinner\nondmatherr@\nopagebreak\else\postdisplaypenalty\@M\fi
 \else\ifvmode\nobreak\else\vadjust{\nobreak}\fi\fi}
\def\nonvmodeerr@#1{\Err@{\string#1\space not allowed within a paragraph
 or in math}}
\def\vnonvmode@#1#2{\relaxnext@\DNii@{\ifx\next\par\DN@{#1}\else
 \DN@{#2}\fi\next@}%
 \ifvmode\DN@{#1}\else
 \DN@{\FN@\nextii@}\fi\next@}
\def\newpage{\vnonvmode@{\vfill\break}{\nonvmodeerr@\newpage}}
\def\smallpagebreak{\vnonvmode@\smallbreak{\nonvmodeerr@\smallpagebreak}}
\def\medpagebreak{\vnonvmode@\medbreak{\nonvmodeerr@\medpagebreak}}
\def\bigpagebreak{\vnonvmode@\bigbreak{\nonvmodeerr@\bigpagebreak}}
\def\NoBlackBoxes{\global\overfullrule\z@}
\def\BlackBoxes{\global\overfullrule5\p@}
\def\Invalid@#1{\def#1{\Err@{\Invalid@@\string#1}}}
\def\Invalid@@{Invalid use of }
\message{figures,}
\Invalid@\caption
\Invalid@\captionwidth
\newdimen\smallcaptionwidth@
\def\topspace{\mid@false\ins@}
\def\midspace{\mid@true\ins@}
\newif\ifmid@
\def\captionfont@{}
\def\ins@#1{\relaxnext@\allowbreak
 \smallcaptionwidth@\captionwidth@\gdef\thespace@{#1}%
 \DN@{\ifx\next\space@\DN@. {\FN@\nextii@}\else
  \DN@.{\FN@\nextii@}\fi\next@.}%
 \DNii@{\ifx\next\caption\DN@\caption{\FN@\nextiii@}%
  \else\let\next@\nextiv@\fi\next@}%
 \def\nextiv@{\vnonvmode@
  {\ifmid@\expandafter\midinsert\else\expandafter\topinsert\fi
   \vbox to\thespace@{}\endinsert}
  {\ifmid@\nonvmodeerr@\midspace\else\nonvmodeerr@\topspace\fi}}%
 \def\nextiii@{\ifx\next\captionwidth\expandafter\nextv@
  \else\expandafter\nextvi@\fi}%
 \def\nextv@\captionwidth##1##2{\smallcaptionwidth@##1\relax\nextvi@{##2}}%
 \def\nextvi@##1{\def\thecaption@{\captionfont@##1}%
  \DN@{\ifx\next\space@\DN@. {\FN@\nextvii@}\else
   \DN@.{\FN@\nextvii@}\fi\next@.}%
  \FN@\next@}%
 \def\nextvii@{\vnonvmode@
  {\ifmid@\expandafter\midinsert\else
  \expandafter\topinsert\fi\vbox to\thespace@{}\nobreak\smallskip
  \setboxz@h{\noindent\ignorespaces\thecaption@\unskip}%
  \ifdim\wdz@>\smallcaptionwidth@\centerline{\vbox{\hsize\smallcaptionwidth@
   \noindent\ignorespaces\thecaption@\unskip}}%
  \else\centerline{\boxz@}\fi\endinsert}
  {\ifmid@\nonvmodeerr@\midspace
  \else\nonvmodeerr@\topspace\fi}}%
 \FN@\next@}
\message{comments,}
\def\newcodes@{\catcode`\\12\catcode`\{12\catcode`\}12\catcode`\#12%
 \catcode`\%12\relax}
\def\oldcodes@{\catcode`\\0\catcode`\{1\catcode`\}2\catcode`\#6%
 \catcode`\%14\relax}
\def\comment{\newcodes@\endlinechar=10 \comment@}
{\lccode`\0=`\\
\lowercase{\gdef\comment@#1^^J{\comment@@#10endcomment\comment@@@}%
\gdef\comment@@#10endcomment{\FN@\comment@@@}%
\gdef\comment@@@#1\comment@@@{\ifx\next\comment@@@\let\next\comment@
 \else\def\next{\oldcodes@\endlinechar=`\^^M\relax}%
 \fi\next}}}
\def\pr@m@s{\ifx'\next\DN@##1{\prim@s}\else\let\next@\egroup\fi\next@}
\def\prime{{\null\prime@\null}}
\mathchardef\prime@="0230
\let\dsize\displaystyle

\let\ssize\scriptstyle

\message{math spacing,}
\def\,{\RIfM@\mskip\thinmuskip\relax\else\kern.16667em\fi}
\def\!{\RIfM@\mskip-\thinmuskip\relax\else\kern-.16667em\fi}
\let\thinspace\,
\let\negthinspace\!
\def\medspace{\RIfM@\mskip\medmuskip\relax\else\kern.222222em\fi}
\def\negmedspace{\RIfM@\mskip-\medmuskip\relax\else\kern-.222222em\fi}
\def\thickspace{\RIfM@\mskip\thickmuskip\relax\else\kern.27777em\fi}
\let\;\thickspace
\def\negthickspace{\RIfM@\mskip-\thickmuskip\relax\else
 \kern-.27777em\fi}
\atdef@,{\RIfM@\mskip.1\thinmuskip\else\leavevmode\null,\fi}
\atdef@!{\RIfM@\mskip-.1\thinmuskip\else\leavevmode\null!\fi}
\atdef@.{\RIfM@&&\else\leavevmode.\spacefactor3000 \fi}
\def\and{\DOTSB\;\mathchar"3026 \;}

\message{fractions,}
\def\frac#1#2{{#1\over#2}}

\newdimen\ex@
\ex@.2326ex
\Invalid@\thickness
\def\thickfrac{\relaxnext@
 \DN@{\ifx\next\thickness\let\next@\nextii@\else
 \DN@{\nextii@\thickness1}\fi\next@}%
 \DNii@\thickness##1##2##3{{##2\above##1\ex@##3}}%
 \FN@\next@}

\def\thickfracwithdelims#1#2{\relaxnext@\def\ldelim@{#1}\def\rdelim@{#2}%
 \DN@{\ifx\next\thickness\let\next@\nextii@\else
 \DN@{\nextii@\thickness1}\fi\next@}%
 \DNii@\thickness##1##2##3{{##2\abovewithdelims
 \ldelim@\rdelim@##1\ex@##3}}%
 \FN@\next@}

\def\:{\nobreak\hskip.1111em\mathpunct{}\nonscript\mkern-\thinmuskip{:}\hskip
 .3333emplus.0555em\relax}
\def\snug{\unskip\kern-\mathsurround}
\message{smash commands,}
\def\topsmash{\top@true\bot@false\smash@}
\def\botsmash{\top@false\bot@true\smash@}
\newif\iftop@
\newif\ifbot@
\def\smash{\top@true\bot@true\smash@}
\def\smash@{\RIfM@\expandafter\mathpalette\expandafter\mathsm@sh\else
 \expandafter\makesm@sh\fi}
\def\finsm@sh{\iftop@\ht\z@\z@\fi\ifbot@\dp\z@\z@\fi\leavevmode\boxz@}
\message{large operator symbols,}
\def\LimitsOnSums{\global\let\slimits@\displaylimits}
\def\NoLimitsOnSums{\global\let\slimits@\nolimits}
\LimitsOnSums
\mathchardef\coprod@="1360       \def\coprod{\DOTSB\coprod@\slimits@}
\mathchardef\bigvee@="1357       \def\bigvee{\DOTSB\bigvee@\slimits@}
\mathchardef\bigwedge@="1356     \def\bigwedge{\DOTSB\bigwedge@\slimits@}
\mathchardef\biguplus@="1355     \def\biguplus{\DOTSB\biguplus@\slimits@}
\mathchardef\bigcap@="1354       \def\bigcap{\DOTSB\bigcap@\slimits@}
\mathchardef\bigcup@="1353       \def\bigcup{\DOTSB\bigcup@\slimits@}
\mathchardef\prod@="1351         \def\prod{\DOTSB\prod@\slimits@}
\mathchardef\sum@="1350          \def\sum{\DOTSB\sum@\slimits@}
\mathchardef\bigotimes@="134E    \def\bigotimes{\DOTSB\bigotimes@\slimits@}
\mathchardef\bigoplus@="134C     \def\bigoplus{\DOTSB\bigoplus@\slimits@}
\mathchardef\bigodot@="134A      \def\bigodot{\DOTSB\bigodot@\slimits@}
\mathchardef\bigsqcup@="1346     \def\bigsqcup{\DOTSB\bigsqcup@\slimits@}
\message{integrals,}
\def\LimitsOnInts{\global\let\ilimits@\displaylimits}
\def\NoLimitsOnInts{\global\let\ilimits@\nolimits}
\NoLimitsOnInts
\def\int{\DOTSI\intop\ilimits@}
\def\oint{\DOTSI\ointop\ilimits@}
\def\intic@{\mathchoice{\hskip.5em}{\hskip.4em}{\hskip.4em}{\hskip.4em}}
\def\negintic@{\mathchoice
 {\hskip-.5em}{\hskip-.4em}{\hskip-.4em}{\hskip-.4em}}
\def\intkern@{\mathchoice{\!\!\!}{\!\!}{\!\!}{\!\!}}
\def\intdots@{\mathchoice{\plaincdots@}
 {{\cdotp}\mkern1.5mu{\cdotp}\mkern1.5mu{\cdotp}}
 {{\cdotp}\mkern1mu{\cdotp}\mkern1mu{\cdotp}}
 {{\cdotp}\mkern1mu{\cdotp}\mkern1mu{\cdotp}}}
\newcount\intno@
\def\iint{\DOTSI\intno@\tw@\FN@\ints@}
\def\iiint{\DOTSI\intno@\thr@@\FN@\ints@}
\def\iiiint{\DOTSI\intno@4 \FN@\ints@}
\def\idotsint{\DOTSI\intno@\z@\FN@\ints@}
\def\ints@{\findlimits@\ints@@}
\newif\iflimtoken@
\newif\iflimits@
\def\findlimits@{\limtoken@true\ifx\next\limits\limits@true
 \else\ifx\next\nolimits\limits@false\else
 \limtoken@false\ifx\ilimits@\nolimits\limits@false\else
 \ifinner\limits@false\else\limits@true\fi\fi\fi\fi}
\def\multint@{\int\ifnum\intno@=\z@\intdots@                                
 \else\intkern@\fi                                                          
 \ifnum\intno@>\tw@\int\intkern@\fi                                         
 \ifnum\intno@>\thr@@\int\intkern@\fi                                       
 \int}                                                                      
\def\multintlimits@{\intop\ifnum\intno@=\z@\intdots@\else\intkern@\fi
 \ifnum\intno@>\tw@\intop\intkern@\fi
 \ifnum\intno@>\thr@@\intop\intkern@\fi\intop}
\def\ints@@{\iflimtoken@                                                    
 \def\ints@@@{\iflimits@\negintic@\mathop{\intic@\multintlimits@}\limits    
  \else\multint@\nolimits\fi                                                
  \eat@}                                                                    
 \else                                                                      
 \def\ints@@@{\iflimits@\negintic@
  \mathop{\intic@\multintlimits@}\limits\else
  \multint@\nolimits\fi}\fi\ints@@@}
\def\LimitsOnNames{\global\let\nlimits@\displaylimits}
\def\NoLimitsOnNames{\global\let\nlimits@\nolimits@}
\LimitsOnNames
\def\nolimits@{\relaxnext@
 \DN@{\ifx\next\limits\DN@\limits{\nolimits}\else
  \let\next@\nolimits\fi\next@}%
 \FN@\next@}
\message{operator names,}
\def\newmcodes@{\mathcode`\'"27\mathcode`\*"2A\mathcode`\."613A%
 \mathcode`\-"2D\mathcode`\/"2F\mathcode`\:"603A }
\def\operatorname#1{\mathop{\newmcodes@\kern\z@\fam\z@#1}\nolimits@}
\def\operatornamewithlimits#1{\mathop{\newmcodes@\kern\z@\fam\z@#1}\nlimits@}
\def\qopname@#1{\mathop{\fam\z@#1}\nolimits@}
\def\qopnamewl@#1{\mathop{\fam\z@#1}\nlimits@}
\def\arccos{\qopname@{arccos}}
\def\arcsin{\qopname@{arcsin}}
\def\arctan{\qopname@{arctan}}
\def\arg{\qopname@{arg}}
\def\cos{\qopname@{cos}}
\def\cosh{\qopname@{cosh}}
\def\cot{\qopname@{cot}}
\def\coth{\qopname@{coth}}
\def\csc{\qopname@{csc}}
\def\deg{\qopname@{deg}}
\def\det{\qopnamewl@{det}}
\def\dim{\qopname@{dim}}
\def\exp{\qopname@{exp}}
\def\gcd{\qopnamewl@{gcd}}
\def\hom{\qopname@{hom}}
\def\inf{\qopnamewl@{inf}}
\def\injlim{\qopnamewl@{inj\,lim}}
\def\ker{\qopname@{ker}}
\def\lg{\qopname@{lg}}
\def\lim{\qopnamewl@{lim}}
\def\liminf{\qopnamewl@{lim\,inf}}
\def\limsup{\qopnamewl@{lim\,sup}}
\def\ln{\qopname@{ln}}
\def\log{\qopname@{log}}
\def\max{\qopnamewl@{max}}
\def\min{\qopnamewl@{min}}
\def\Pr{\qopnamewl@{Pr}}
\def\projlim{\qopnamewl@{proj\,lim}}
\def\sec{\qopname@{sec}}
\def\sin{\qopname@{sin}}
\def\sinh{\qopname@{sinh}}
\def\sup{\qopnamewl@{sup}}
\def\tan{\qopname@{tan}}
\def\tanh{\qopname@{tanh}}
\def\varinjlim{\mathop{\vtop{\ialign{##\crcr
 \hfil\rm lim\hfil\crcr\noalign{\nointerlineskip}\rightarrowfill\crcr
 \noalign{\nointerlineskip\kern-\ex@}\crcr}}}}
\def\varprojlim{\mathop{\vtop{\ialign{##\crcr
 \hfil\rm lim\hfil\crcr\noalign{\nointerlineskip}\leftarrowfill\crcr
 \noalign{\nointerlineskip\kern-\ex@}\crcr}}}}
\def\varliminf{\mathop{\underline{\vrule height\z@ depth.2exwidth\z@
 \hbox{\rm lim}}}}

\newdimen\buffer@
\buffer@\fontdimen13 \tenex
\newdimen\buffer
\buffer\buffer@

\def\ResetBuffer{\fontdimen13 \tenex\buffer@\global\buffer\buffer@}
\def\shave#1{\mathop{\hbox{$\m@th\fontdimen13 \tenex\z@                     
 \displaystyle{#1}$}}\fontdimen13 \tenex\buffer}

\message{multilevel sub/superscripts,}
\Invalid@\\
\def\Let@{\relax\iffalse{\fi\let\\=\cr\iffalse}\fi}
\Invalid@\vspace
\def\vspace@{\def\vspace##1{\crcr\noalign{\vskip##1\relax}}}
\def\multilimits@{\bgroup\vspace@\Let@
 \baselineskip\fontdimen10 \scriptfont\tw@
 \advance\baselineskip\fontdimen12 \scriptfont\tw@
 \lineskip\thr@@\fontdimen8 \scriptfont\thr@@
 \lineskiplimit\lineskip
 \vbox\bgroup\ialign\bgroup\hfil$\m@th\scriptstyle{##}$\hfil\crcr}
\def\Sb{_\multilimits@}
\def\endSb{\crcr\egroup\egroup\egroup}
\def\Sp{^\multilimits@}

\def\spreadlines#1{\RIfMIfI@\onlydmatherr@\spreadlines\else
 \openup#1\relax\fi\else\onlydmatherr@\spreadlines\fi}
\def\Mathstrut@{\copy\Mathstrutbox@}
\newbox\Mathstrutbox@
\setbox\Mathstrutbox@\null
\setboxz@h{$\m@th($}
\ht\Mathstrutbox@\ht\z@
\dp\Mathstrutbox@\dp\z@
\message{matrices,}
\newdimen\spreadmlines@
\def\spreadmatrixlines#1{\RIfMIfI@
 \onlydmatherr@\spreadmatrixlines\else
 \spreadmlines@#1\relax\fi\else\onlydmatherr@\spreadmatrixlines\fi}
\def\matrix{\null\,\vcenter\bgroup\Let@\vspace@
 \normalbaselines\openup\spreadmlines@\ialign
 \bgroup\hfil$\m@th##$\hfil&&\quad\hfil$\m@th##$\hfil\crcr
 \Mathstrut@\crcr\noalign{\kern-\baselineskip}}
\def\endmatrix{\crcr\Mathstrut@\crcr\noalign{\kern-\baselineskip}\egroup
 \egroup\,}
\def\format{\crcr\egroup\iffalse{\fi\ifnum`}=0 \fi\format@}
\newtoks\hashtoks@
\hashtoks@{#}
\def\format@#1\\{\def\preamble@{#1}%
 \def\l{$\m@th\the\hashtoks@$\hfil}%
 \def\c{\hfil$\m@th\the\hashtoks@$\hfil}%
 \def\r{\hfil$\m@th\the\hashtoks@$}%
 \edef\preamble@@{\preamble@}\ifnum`{=0 \fi\iffalse}\fi
 \ialign\bgroup\span\preamble@@\crcr}
\def\smallmatrix{\null\,\vcenter\bgroup\vspace@\Let@
 \baselineskip9\ex@\lineskip\ex@
 \ialign\bgroup\hfil$\m@th\scriptstyle{##}$\hfil&&\thickspace\hfil
 $\m@th\scriptstyle{##}$\hfil\crcr}
\def\endsmallmatrix{\crcr\egroup\egroup\,}
\def\pmatrix{\left(\matrix}
\def\endpmatrix{\endmatrix\right)}

\newmuskip\dotsspace@
\dotsspace@1.5mu
\def\strip@#1 {#1}
\def\spacehdots#1\for#2{\multispan{#2}\xleaders
 \hbox{$\m@th\mkern\strip@#1 \dotsspace@.\mkern\strip@#1 \dotsspace@$}\hfill}
\def\hdotsfor#1{\spacehdots\@ne\for{#1}}
\def\multispan@#1{\omit\mscount#1\unskip\loop\ifnum\mscount>\@ne\sp@n\repeat}
\def\spaceinnerhdots#1\for#2\after#3{\multispan@{\strip@#2 }#3\xleaders
 \hbox{$\m@th\mkern\strip@#1 \dotsspace@.\mkern\strip@#1 \dotsspace@$}\hfill}
\def\innerhdotsfor#1\after#2{\spaceinnerhdots\@ne\for#1\after{#2}}
\def\cases{\bgroup\spreadmlines@\jot\left\{\,\matrix\format\l&\quad\l\\}
\def\endcases{\endmatrix\right.\egroup}
\message{multiline displays,}
\newif\ifinany@
\newif\ifinalign@
\newif\ifingather@
\def\strut@{\copy\strutbox@}
\newbox\strutbox@
\setbox\strutbox@\hbox{\vrule height8\p@ depth3\p@ width\z@}
\def\topaligned{\null\,\vtop\aligned@}
\def\botaligned{\null\,\vbox\aligned@}
\def\aligned{\null\,\vcenter\aligned@}
\def\aligned@{\bgroup\vspace@\Let@
 \ifinany@\else\openup\jot\fi\ialign
 \bgroup\hfil\strut@$\m@th\displaystyle{##}$&
 $\m@th\displaystyle{{}##}$\hfil\crcr}
\def\endaligned{\crcr\egroup\egroup}

\def\alignedat#1{\null\,\vcenter\bgroup\doat@{#1}\vspace@\Let@
 \ifinany@\else\openup\jot\fi\ialign\bgroup\span\preamble@@\crcr}
\newcount\atcount@
\def\doat@#1{\toks@{\hfil\strut@$\m@th
 \displaystyle{\the\hashtoks@}$&$\m@th\displaystyle
 {{}\the\hashtoks@}$\hfil}
 \atcount@#1\relax\advance\atcount@\m@ne                                    
 \loop\ifnum\atcount@>\z@\toks@=\expandafter{\the\toks@&\hfil$\m@th
 \displaystyle{\the\hashtoks@}$&$\m@th
 \displaystyle{{}\the\hashtoks@}$\hfil}\advance
  \atcount@\m@ne\repeat                                                     
 \xdef\preamble@{\the\toks@}\xdef\preamble@@{\preamble@}}

\def\gathered{\null\,\vcenter\bgroup\vspace@\Let@
 \ifinany@\else\openup\jot\fi\ialign
 \bgroup\hfil\strut@$\m@th\displaystyle{##}$\hfil\crcr}
\def\endgathered{\crcr\egroup\egroup}
\newif\iftagsleft@
\def\TagsOnLeft{\global\tagsleft@true}
\def\TagsOnRight{\global\tagsleft@false}
\TagsOnLeft
\newif\ifmathtags@
\def\TagsAsMath{\global\mathtags@true}
\def\TagsAsText{\global\mathtags@false}
\TagsAsText
\def\tagform@#1{\hbox{\rm(\ignorespaces#1\unskip)}}
\def\thetag{\leavevmode\tagform@}
\def\tag#1$${\iftagsleft@\leqno\else\eqno\fi                                
 \maketag@#1\maketag@                                                       
 $$}                                                                        
\def\maketag@{\FN@\maketag@@}
\def\maketag@@{\ifx\next"\expandafter\maketag@@@\else\expandafter\maketag@@@@
 \fi}
\def\maketag@@@"#1"#2\maketag@{\hbox{\rm#1}}                                
\def\maketag@@@@#1\maketag@{\ifmathtags@\tagform@{$\m@th#1$}\else
 \tagform@{#1}\fi}
\interdisplaylinepenalty\@M
\def\allowdisplaybreaks{\RIfMIfI@
 \onlydmatherr@\allowdisplaybreaks\else
 \interdisplaylinepenalty\z@\fi\else\onlydmatherr@\allowdisplaybreaks\fi}
\Invalid@\allowdisplaybreak
\Invalid@\displaybreak
\Invalid@\intertext
\def\allowdisplaybreak@{\def\allowdisplaybreak{\crcr\noalign{\allowbreak}}}
\def\displaybreak@{\def\displaybreak{\crcr\noalign{\break}}}
\def\intertext@{\def\intertext##1{\crcr\noalign{%
 \penalty\postdisplaypenalty \vskip\belowdisplayskip
 \vbox{\normalbaselines\noindent##1}%
 \penalty\predisplaypenalty \vskip\abovedisplayskip}}}
\newskip\centering@
\centering@\z@ plus\@m\p@
\def\align{\relax\ifingather@\DN@{\csname align (in
  \string\gather)\endcsname}\else
 \ifmmode\ifinner\DN@{\onlydmatherr@\align}\else
  \let\next@\align@\fi
 \else\DN@{\onlydmatherr@\align}\fi\fi\next@}
\newhelp\andhelp@
{An extra & here is so disastrous that you should probably exit^^J
and fix things up.}
\newif\iftag@
\newcount\and@
\def\align@{\inalign@true\inany@true
 \vspace@\allowdisplaybreak@\displaybreak@\intertext@
 \def\tag{\global\tag@true\ifnum\and@=\z@\DN@{&&}\else
          \DN@{&}\fi\next@}%
 \iftagsleft@\DN@{\csname align \endcsname}\else
  \DN@{\csname align \space\endcsname}\fi\next@}
\def\Tag@{\iftag@\else\errhelp\andhelp@\err@{Extra & on this line}\fi}
\newdimen\lwidth@
\newdimen\rwidth@
\newdimen\maxlwidth@
\newdimen\maxrwidth@
\newdimen\totwidth@
\def\measure@#1\endalign{\lwidth@\z@\rwidth@\z@\maxlwidth@\z@\maxrwidth@\z@
 \global\and@\z@                                                            
 \setbox@ne\vbox                                                            
  {\everycr{\noalign{\global\tag@false\global\and@\z@}}\Let@                
  \halign{\setboxz@h{$\m@th\displaystyle{\@lign##}$}
   \global\lwidth@\wdz@                                                     
   \ifdim\lwidth@>\maxlwidth@\global\maxlwidth@\lwidth@\fi                  
   \global\advance\and@\@ne                                                 
   &\setboxz@h{$\m@th\displaystyle{{}\@lign##}$}\global\rwidth@\wdz@        
   \ifdim\rwidth@>\maxrwidth@\global\maxrwidth@\rwidth@\fi                  
   \global\advance\and@\@ne                                                
   &\Tag@
   \eat@{##}\crcr#1\crcr}}
 \totwidth@\maxlwidth@\advance\totwidth@\maxrwidth@}                       
\def\displ@y@{\global\dt@ptrue\openup\jot
 \everycr{\noalign{\global\tag@false\global\and@\z@\ifdt@p\global\dt@pfalse
 \vskip-\lineskiplimit\vskip\normallineskiplimit\else
 \penalty\interdisplaylinepenalty\fi}}}
\def\black@#1{\noalign{\ifdim#1>\displaywidth
 \dimen@\prevdepth\nointerlineskip                                          
 \vskip-\ht\strutbox@\vskip-\dp\strutbox@                                   
 \vbox{\noindent\hbox to#1{\strut@\hfill}}
 \prevdepth\dimen@                                                          
 \fi}}
\expandafter\def\csname align \space\endcsname#1\endalign
 {\measure@#1\endalign\global\and@\z@                                       
 \ifingather@\everycr{\noalign{\global\and@\z@}}\else\displ@y@\fi           
 \Let@\tabskip\centering@                                                   
 \halign to\displaywidth
  {\hfil\strut@\setboxz@h{$\m@th\displaystyle{\@lign##}$}
  \global\lwidth@\wdz@\boxz@\global\advance\and@\@ne                        
  \tabskip\z@skip                                                           
  &\setboxz@h{$\m@th\displaystyle{{}\@lign##}$}
  \global\rwidth@\wdz@\boxz@\hfill\global\advance\and@\@ne                  
  \tabskip\centering@                                                       
  &\setboxz@h{\@lign\strut@\maketag@##\maketag@}
  \dimen@\displaywidth\advance\dimen@-\totwidth@
  \divide\dimen@\tw@\advance\dimen@\maxrwidth@\advance\dimen@-\rwidth@     
  \ifdim\dimen@<\tw@\wdz@\llap{\vtop{\normalbaselines\null\boxz@}}
  \else\llap{\boxz@}\fi                                                    
  \tabskip\z@skip                                                          
  \crcr#1\crcr                                                             
  \black@\totwidth@}}                                                      
\newdimen\lineht@
\expandafter\def\csname align \endcsname#1\endalign{\measure@#1\endalign
 \global\and@\z@
 \ifdim\totwidth@>\displaywidth\let\displaywidth@\totwidth@\else
  \let\displaywidth@\displaywidth\fi                                        
 \ifingather@\everycr{\noalign{\global\and@\z@}}\else\displ@y@\fi
 \Let@\tabskip\centering@\halign to\displaywidth
  {\hfil\strut@\setboxz@h{$\m@th\displaystyle{\@lign##}$}%
  \global\lwidth@\wdz@\global\lineht@\ht\z@                                 
  \boxz@\global\advance\and@\@ne
  \tabskip\z@skip&\setboxz@h{$\m@th\displaystyle{{}\@lign##}$}%
  \global\rwidth@\wdz@\ifdim\ht\z@>\lineht@\global\lineht@\ht\z@\fi         
  \boxz@\hfil\global\advance\and@\@ne
  \tabskip\centering@&\kern-\displaywidth@                                  
  \setboxz@h{\@lign\strut@\maketag@##\maketag@}%
  \dimen@\displaywidth\advance\dimen@-\totwidth@
  \divide\dimen@\tw@\advance\dimen@\maxlwidth@\advance\dimen@-\lwidth@
  \ifdim\dimen@<\tw@\wdz@
   \rlap{\vbox{\normalbaselines\boxz@\vbox to\lineht@{}}}\else
   \rlap{\boxz@}\fi
  \tabskip\displaywidth@\crcr#1\crcr\black@\totwidth@}}
\expandafter\def\csname align (in \string\gather)\endcsname
  #1\endalign{\vcenter{\align@#1\endalign}}
\Invalid@\endalign
\newif\ifxat@
\def\alignat{\RIfMIfI@\DN@{\onlydmatherr@\alignat}\else
 \DN@{\csname alignat \endcsname}\fi\else
 \DN@{\onlydmatherr@\alignat}\fi\next@}
\newif\ifmeasuring@
\newbox\savealignat@
\expandafter\def\csname alignat \endcsname#1#2\endalignat                   
 {\inany@true\xat@false
 \def\tag{\global\tag@true\count@#1\relax\multiply\count@\tw@
  \xdef\tag@{}\loop\ifnum\count@>\and@\xdef\tag@{&\tag@}\advance\count@\m@ne
  \repeat\tag@}%
 \vspace@\allowdisplaybreak@\displaybreak@\intertext@
 \displ@y@\measuring@true                                                   
 \setbox\savealignat@\hbox{$\m@th\displaystyle\Let@
  \attag@{#1}
  \vbox{\halign{\span\preamble@@\crcr#2\crcr}}$}%
 \measuring@false                                                           
 \Let@\attag@{#1}
 \tabskip\centering@\halign to\displaywidth
  {\span\preamble@@\crcr#2\crcr                                             
  \black@{\wd\savealignat@}}}                                               
\Invalid@\endalignat
\def\xalignat{\RIfMIfI@
 \DN@{\onlydmatherr@\xalignat}\else
 \DN@{\csname xalignat \endcsname}\fi\else
 \DN@{\onlydmatherr@\xalignat}\fi\next@}
\expandafter\def\csname xalignat \endcsname#1#2\endxalignat
 {\inany@true\xat@true
 \def\tag{\global\tag@true\def\tag@{}\count@#1\relax\multiply\count@\tw@
  \loop\ifnum\count@>\and@\xdef\tag@{&\tag@}\advance\count@\m@ne\repeat\tag@}%
 \vspace@\allowdisplaybreak@\displaybreak@\intertext@
 \displ@y@\measuring@true\setbox\savealignat@\hbox{$\m@th\displaystyle\Let@
 \attag@{#1}\vbox{\halign{\span\preamble@@\crcr#2\crcr}}$}%
 \measuring@false\Let@
 \attag@{#1}\tabskip\centering@\halign to\displaywidth
 {\span\preamble@@\crcr#2\crcr\black@{\wd\savealignat@}}}
\def\attag@#1{\let\Maketag@\maketag@\let\TAG@\Tag@                          
 \let\Tag@=0\let\maketag@=0
 \ifmeasuring@\def\llap@##1{\setboxz@h{##1}\hbox to\tw@\wdz@{}}%
  \def\rlap@##1{\setboxz@h{##1}\hbox to\tw@\wdz@{}}\else
  \let\llap@\llap\let\rlap@\rlap\fi                                         
 \toks@{\hfil\strut@$\m@th\displaystyle{\@lign\the\hashtoks@}$\tabskip\z@skip
  \global\advance\and@\@ne&$\m@th\displaystyle{{}\@lign\the\hashtoks@}$\hfil
  \ifxat@\tabskip\centering@\fi\global\advance\and@\@ne}
 \iftagsleft@
  \toks@@{\tabskip\centering@&\Tag@\kern-\displaywidth
   \rlap@{\@lign\maketag@\the\hashtoks@\maketag@}%
   \global\advance\and@\@ne\tabskip\displaywidth}\else
  \toks@@{\tabskip\centering@&\Tag@\llap@{\@lign\maketag@
   \the\hashtoks@\maketag@}\global\advance\and@\@ne\tabskip\z@skip}\fi      
 \atcount@#1\relax\advance\atcount@\m@ne
 \loop\ifnum\atcount@>\z@
 \toks@=\expandafter{\the\toks@&\hfil$\m@th\displaystyle{\@lign
  \the\hashtoks@}$\global\advance\and@\@ne
  \tabskip\z@skip&$\m@th\displaystyle{{}\@lign\the\hashtoks@}$\hfil\ifxat@
  \tabskip\centering@\fi\global\advance\and@\@ne}\advance\atcount@\m@ne
 \repeat                                                                    
 \xdef\preamble@{\the\toks@\the\toks@@}
 \xdef\preamble@@{\preamble@}
 \let\maketag@\Maketag@\let\Tag@\TAG@}                                      
\Invalid@\endxalignat
\def\xxalignat{\RIfMIfI@
 \DN@{\onlydmatherr@\xxalignat}\else\DN@{\csname xxalignat
  \endcsname}\fi\else
 \DN@{\onlydmatherr@\xxalignat}\fi\next@}
\expandafter\def\csname xxalignat \endcsname#1#2\endxxalignat{\inany@true
 \vspace@\allowdisplaybreak@\displaybreak@\intertext@
 \displ@y\setbox\savealignat@\hbox{$\m@th\displaystyle\Let@
 \xxattag@{#1}\vbox{\halign{\span\preamble@@\crcr#2\crcr}}$}%
 \Let@\xxattag@{#1}\tabskip\z@skip\halign to\displaywidth
 {\span\preamble@@\crcr#2\crcr\black@{\wd\savealignat@}}}
\def\xxattag@#1{\toks@{\tabskip\z@skip\hfil\strut@
 $\m@th\displaystyle{\the\hashtoks@}$&%
 $\m@th\displaystyle{{}\the\hashtoks@}$\hfil\tabskip\centering@&}%
 \atcount@#1\relax\advance\atcount@\m@ne\loop\ifnum\atcount@>\z@
 \toks@=\expandafter{\the\toks@&\hfil$\m@th\displaystyle{\the\hashtoks@}$%
  \tabskip\z@skip&$\m@th\displaystyle{{}\the\hashtoks@}$\hfil
  \tabskip\centering@}\advance\atcount@\m@ne\repeat
 \xdef\preamble@{\the\toks@\tabskip\z@skip}\xdef\preamble@@{\preamble@}}
\Invalid@\endxxalignat
\newdimen\gwidth@
\newdimen\gmaxwidth@
\def\gmeasure@#1\endgather{\gwidth@\z@\gmaxwidth@\z@\setbox@ne\vbox{\Let@
 \halign{\setboxz@h{$\m@th\displaystyle{##}$}\global\gwidth@\wdz@
 \ifdim\gwidth@>\gmaxwidth@\global\gmaxwidth@\gwidth@\fi
 &\eat@{##}\crcr#1\crcr}}}
\def\gather{\RIfMIfI@\DN@{\onlydmatherr@\gather}\else
 \ingather@true\inany@true\def\tag{&}%
 \vspace@\allowdisplaybreak@\displaybreak@\intertext@
 \displ@y\Let@
 \iftagsleft@\DN@{\csname gather \endcsname}\else
  \DN@{\csname gather \space\endcsname}\fi\fi
 \else\DN@{\onlydmatherr@\gather}\fi\next@}
\expandafter\def\csname gather \space\endcsname#1\endgather
 {\gmeasure@#1\endgather\tabskip\centering@
 \halign to\displaywidth{\hfil\strut@\setboxz@h{$\m@th\displaystyle{##}$}%
 \global\gwidth@\wdz@\boxz@\hfil&
 \setboxz@h{\strut@{\maketag@##\maketag@}}%
 \dimen@\displaywidth\advance\dimen@-\gwidth@
 \ifdim\dimen@>\tw@\wdz@\llap{\boxz@}\else
 \llap{\vtop{\normalbaselines\null\boxz@}}\fi
 \tabskip\z@skip\crcr#1\crcr\black@\gmaxwidth@}}
\newdimen\glineht@
\expandafter\def\csname gather \endcsname#1\endgather{\gmeasure@#1\endgather
 \ifdim\gmaxwidth@>\displaywidth\let\gdisplaywidth@\gmaxwidth@\else
 \let\gdisplaywidth@\displaywidth\fi\tabskip\centering@\halign to\displaywidth
 {\hfil\strut@\setboxz@h{$\m@th\displaystyle{##}$}%
 \global\gwidth@\wdz@\global\glineht@\ht\z@\boxz@\hfil&\kern-\gdisplaywidth@
 \setboxz@h{\strut@{\maketag@##\maketag@}}%
 \dimen@\displaywidth\advance\dimen@-\gwidth@
 \ifdim\dimen@>\tw@\wdz@\rlap{\boxz@}\else
 \rlap{\vbox{\normalbaselines\boxz@\vbox to\glineht@{}}}\fi
 \tabskip\gdisplaywidth@\crcr#1\crcr\black@\gmaxwidth@}}
\newif\ifctagsplit@
\def\CenteredTagsOnSplits{\global\ctagsplit@true}
\def\TopOrBottomTagsOnSplits{\global\ctagsplit@false}
\TopOrBottomTagsOnSplits
\def\split{\relax\ifinany@\let\next@\insplit@\else
 \ifmmode\ifinner\def\next@{\onlydmatherr@\split}\else
 \let\next@\outsplit@\fi\else
 \def\next@{\onlydmatherr@\split}\fi\fi\next@}
\def\insplit@{\global\setbox\z@\vbox\bgroup\vspace@\Let@\ialign\bgroup
 \hfil\strut@$\m@th\displaystyle{##}$&$\m@th\displaystyle{{}##}$\hfill\crcr}
\def\endsplit{\crcr\egroup\egroup\iftagsleft@\expandafter\lendsplit@\else
 \expandafter\rendsplit@\fi}
\def\rendsplit@{\global\setbox9 \vbox
 {\unvcopy\z@\global\setbox8 \lastbox\unskip}
 \setbox@ne\hbox{\unhcopy8 \unskip\global\setbox\tw@\lastbox
 \unskip\global\setbox\thr@@\lastbox}
 \global\setbox7 \hbox{\unhbox\tw@\unskip}
 \ifinalign@\ifctagsplit@                                                   
  \gdef\split@{\hbox to\wd\thr@@{}&
   \vcenter{\vbox{\moveleft\wd\thr@@\boxz@}}}
 \else\gdef\split@{&\vbox{\moveleft\wd\thr@@\box9}\crcr
  \box\thr@@&\box7}\fi                                                      
 \else                                                                      
  \ifctagsplit@\gdef\split@{\vcenter{\boxz@}}\else
  \gdef\split@{\box9\crcr\hbox{\box\thr@@\box7}}\fi
 \fi
 \split@}                                                                   
\def\lendsplit@{\global\setbox9\vtop{\unvcopy\z@}
 \setbox@ne\vbox{\unvcopy\z@\global\setbox8\lastbox}
 \setbox@ne\hbox{\unhcopy8\unskip\setbox\tw@\lastbox
  \unskip\global\setbox\thr@@\lastbox}
 \ifinalign@\ifctagsplit@                                                   
  \gdef\split@{\hbox to\wd\thr@@{}&
  \vcenter{\vbox{\moveleft\wd\thr@@\box9}}}
  \else                                                                     
  \gdef\split@{\hbox to\wd\thr@@{}&\vbox{\moveleft\wd\thr@@\box9}}\fi
 \else
  \ifctagsplit@\gdef\split@{\vcenter{\box9}}\else
  \gdef\split@{\box9}\fi
 \fi\split@}
\def\outsplit@#1$${\align\insplit@#1\endalign$$}
\newdimen\multlinegap@
\multlinegap@1em
\newdimen\multlinetaggap@
\multlinetaggap@1em
\def\MultlineGap#1{\global\multlinegap@#1\relax}
\def\multlinegap#1{\RIfMIfI@\onlydmatherr@\multlinegap\else
 \multlinegap@#1\relax\fi\else\onlydmatherr@\multlinegap\fi}
\def\nomultlinegap{\multlinegap{\z@}}
\def\multline{\RIfMIfI@
 \DN@{\onlydmatherr@\multline}\else
 \DN@{\multline@}\fi\else
 \DN@{\onlydmatherr@\multline}\fi\next@}
\newif\iftagin@
\def\tagin@#1{\tagin@false\in@\tag{#1}\ifin@\tagin@true\fi}
\def\multline@#1$${\inany@true\vspace@\allowdisplaybreak@\displaybreak@
 \tagin@{#1}\iftagsleft@\DN@{\multline@l#1$$}\else
 \DN@{\multline@r#1$$}\fi\next@}
\newdimen\mwidth@
\def\rmmeasure@#1\endmultline{%
 \def\shoveleft##1{##1}\def\shoveright##1{##1}
 \setbox@ne\vbox{\Let@\halign{\setboxz@h
  {$\m@th\@lign\displaystyle{}##$}\global\mwidth@\wdz@
  \crcr#1\crcr}}}
\newdimen\mlineht@
\newif\ifzerocr@
\newif\ifonecr@
\def\lmmeasure@#1\endmultline{\global\zerocr@true\global\onecr@false
 \everycr{\noalign{\ifonecr@\global\onecr@false\fi
  \ifzerocr@\global\zerocr@false\global\onecr@true\fi}}
  \def\shoveleft##1{##1}\def\shoveright##1{##1}%
 \setbox@ne\vbox{\Let@\halign{\setboxz@h
  {$\m@th\@lign\displaystyle{}##$}\ifonecr@\global\mwidth@\wdz@
  \global\mlineht@\ht\z@\fi\crcr#1\crcr}}}
\newbox\mtagbox@
\newdimen\ltwidth@
\newdimen\rtwidth@
\def\multline@l#1$${\iftagin@\DN@{\lmultline@@#1$$}\else
 \DN@{\setbox\mtagbox@\null\ltwidth@\z@\rtwidth@\z@
  \lmultline@@@#1$$}\fi\next@}
\def\lmultline@@#1\endmultline\tag#2$${%
 \setbox\mtagbox@\hbox{\maketag@#2\maketag@}
 \lmmeasure@#1\endmultline\dimen@\mwidth@\advance\dimen@\wd\mtagbox@
 \advance\dimen@\multlinetaggap@                                            
 \ifdim\dimen@>\displaywidth\ltwidth@\z@\else\ltwidth@\wd\mtagbox@\fi       
 \lmultline@@@#1\endmultline$$}
\def\lmultline@@@{\displ@y
 \def\shoveright##1{##1\hfilneg\hskip\multlinegap@}%
 \def\shoveleft##1{\setboxz@h{$\m@th\displaystyle{}##1$}%
  \setbox@ne\hbox{$\m@th\displaystyle##1$}%
  \hfilneg
  \iftagin@
   \ifdim\ltwidth@>\z@\hskip\ltwidth@\hskip\multlinetaggap@\fi
  \else\hskip\multlinegap@\fi\hskip.5\wd@ne\hskip-.5\wdz@##1}
  \halign\bgroup\Let@\hbox to\displaywidth
   {\strut@$\m@th\displaystyle\hfil{}##\hfil$}\crcr
   \hfilneg                                                                 
   \iftagin@                                                                
    \ifdim\ltwidth@>\z@                                                     
     \box\mtagbox@\hskip\multlinetaggap@                                    
    \else
     \rlap{\vbox{\normalbaselines\hbox{\strut@\box\mtagbox@}%
     \vbox to\mlineht@{}}}\fi                                               
   \else\hskip\multlinegap@\fi}                                             
\def\multline@r#1$${\iftagin@\DN@{\rmultline@@#1$$}\else
 \DN@{\setbox\mtagbox@\null\ltwidth@\z@\rtwidth@\z@
  \rmultline@@@#1$$}\fi\next@}
\def\rmultline@@#1\endmultline\tag#2$${\ltwidth@\z@
 \setbox\mtagbox@\hbox{\maketag@#2\maketag@}%
 \rmmeasure@#1\endmultline\dimen@\mwidth@\advance\dimen@\wd\mtagbox@
 \advance\dimen@\multlinetaggap@
 \ifdim\dimen@>\displaywidth\rtwidth@\z@\else\rtwidth@\wd\mtagbox@\fi
 \rmultline@@@#1\endmultline$$}
\def\rmultline@@@{\displ@y
 \def\shoveright##1{##1\hfilneg\iftagin@\ifdim\rtwidth@>\z@
  \hskip\rtwidth@\hskip\multlinetaggap@\fi\else\hskip\multlinegap@\fi}%
 \def\shoveleft##1{\setboxz@h{$\m@th\displaystyle{}##1$}%
  \setbox@ne\hbox{$\m@th\displaystyle##1$}%
  \hfilneg\hskip\multlinegap@\hskip.5\wd@ne\hskip-.5\wdz@##1}%
 \halign\bgroup\Let@\hbox to\displaywidth
  {\strut@$\m@th\displaystyle\hfil{}##\hfil$}\crcr
 \hfilneg\hskip\multlinegap@}
\def\endmultline{\iftagsleft@\expandafter\lendmultline@\else
 \expandafter\rendmultline@\fi}
\def\lendmultline@{\hfilneg\hskip\multlinegap@\crcr\egroup}
\def\rendmultline@{\iftagin@                                                
 \ifdim\rtwidth@>\z@                                                        
  \hskip\multlinetaggap@\box\mtagbox@                                       
 \else\llap{\vtop{\normalbaselines\null\hbox{\strut@\box\mtagbox@}}}\fi     
 \else\hskip\multlinegap@\fi                                                
 \hfilneg\crcr\egroup}
\def\bmod{\mskip-\medmuskip\mkern5mu\mathbin{\fam\z@ mod}\penalty900
 \mkern5mu\mskip-\medmuskip}
\def\pmod#1{\allowbreak\ifinner\mkern8mu\else\mkern18mu\fi
 ({\fam\z@ mod}\,\,#1)}
\def\pod#1{\allowbreak\ifinner\mkern8mu\else\mkern18mu\fi(#1)}
\def\mod#1{\allowbreak\ifinner\mkern12mu\else\mkern18mu\fi{\fam\z@ mod}\,\,#1}
\message{continued fractions,}
\newcount\cfraccount@
\def\cfrac{\bgroup\bgroup\advance\cfraccount@\@ne\strut
 \iffalse{\fi\def\\{\over\displaystyle}\iffalse}\fi}
\def\lcfrac{\bgroup\bgroup\advance\cfraccount@\@ne\strut
 \iffalse{\fi\def\\{\hfill\over\displaystyle}\iffalse}\fi}
\def\rcfrac{\bgroup\bgroup\advance\cfraccount@\@ne\strut\hfill
 \iffalse{\fi\def\\{\over\displaystyle}\iffalse}\fi}
\def\gloop@#1\repeat{\gdef\body{#1}\iterate}
\def\endcfrac{\gloop@\ifnum\cfraccount@>\z@\global\advance\cfraccount@\m@ne
 \egroup\hskip-\nulldelimiterspace\egroup\repeat}
\message{compound symbols,}
\def\binrel@#1{\setboxz@h{\thinmuskip0mu
  \medmuskip\m@ne mu\thickmuskip\@ne mu$#1\m@th$}%
 \setbox@ne\hbox{\thinmuskip0mu\medmuskip\m@ne mu\thickmuskip
  \@ne mu${}#1{}\m@th$}%
 \setbox\tw@\hbox{\hskip\wd@ne\hskip-\wdz@}}
\def\overset#1\to#2{\binrel@{#2}\ifdim\wd\tw@<\z@
 \mathbin{\mathop{\kern\z@#2}\limits^{#1}}\else\ifdim\wd\tw@>\z@
 \mathrel{\mathop{\kern\z@#2}\limits^{#1}}\else
 {\mathop{\kern\z@#2}\limits^{#1}}{}\fi\fi}
\def\underset#1\to#2{\binrel@{#2}\ifdim\wd\tw@<\z@
 \mathbin{\mathop{\kern\z@#2}\limits_{#1}}\else\ifdim\wd\tw@>\z@
 \mathrel{\mathop{\kern\z@#2}\limits_{#1}}\else
 {\mathop{\kern\z@#2}\limits_{#1}}{}\fi\fi}
\def\oversetbrace#1\to#2{\overbrace{#2}^{#1}}
\def\undersetbrace#1\to#2{\underbrace{#2}_{#1}}
\def\sideset#1\and#2\to#3{%
 \setbox@ne\hbox{$\dsize{\vphantom{#3}}#1{#3}\m@th$}%
 \setbox\tw@\hbox{$\dsize{#3}#2\m@th$}%
 \hskip\wd@ne\hskip-\wd\tw@\mathop{\hskip\wd\tw@\hskip-\wd@ne
  {\vphantom{#3}}#1{#3}#2}}
\def\rightarrowfill@#1{\setboxz@h{$#1-\m@th$}\ht\z@\z@
  $#1\m@th\copy\z@\mkern-6mu\cleaders
  \hbox{$#1\mkern-2mu\box\z@\mkern-2mu$}\hfill
  \mkern-6mu\mathord\rightarrow$}
\def\leftarrowfill@#1{\setboxz@h{$#1-\m@th$}\ht\z@\z@
  $#1\m@th\mathord\leftarrow\mkern-6mu\cleaders
  \hbox{$#1\mkern-2mu\copy\z@\mkern-2mu$}\hfill
  \mkern-6mu\box\z@$}
\def\leftrightarrowfill@#1{\setboxz@h{$#1-\m@th$}\ht\z@\z@
  $#1\m@th\mathord\leftarrow\mkern-6mu\cleaders
  \hbox{$#1\mkern-2mu\box\z@\mkern-2mu$}\hfill
  \mkern-6mu\mathord\rightarrow$}
\def\overrightarrow{\mathpalette\overrightarrow@}
\def\overrightarrow@#1#2{\vbox{\ialign{##\crcr\rightarrowfill@#1\crcr
 \noalign{\kern-\ex@\nointerlineskip}$\m@th\hfil#1#2\hfil$\crcr}}}

\def\overleftarrow{\mathpalette\overleftarrow@}
\def\overleftarrow@#1#2{\vbox{\ialign{##\crcr\leftarrowfill@#1\crcr
 \noalign{\kern-\ex@\nointerlineskip}$\m@th\hfil#1#2\hfil$\crcr}}}
\def\overleftrightarrow{\mathpalette\overleftrightarrow@}
\def\overleftrightarrow@#1#2{\vbox{\ialign{##\crcr\leftrightarrowfill@#1\crcr
 \noalign{\kern-\ex@\nointerlineskip}$\m@th\hfil#1#2\hfil$\crcr}}}
\def\underrightarrow{\mathpalette\underrightarrow@}
\def\underrightarrow@#1#2{\vtop{\ialign{##\crcr$\m@th\hfil#1#2\hfil$\crcr
 \noalign{\nointerlineskip}\rightarrowfill@#1\crcr}}}

\def\underleftarrow{\mathpalette\underleftarrow@}
\def\underleftarrow@#1#2{\vtop{\ialign{##\crcr$\m@th\hfil#1#2\hfil$\crcr
 \noalign{\nointerlineskip}\leftarrowfill@#1\crcr}}}
\def\underleftrightarrow{\mathpalette\underleftrightarrow@}
\def\underleftrightarrow@#1#2{\vtop{\ialign{##\crcr$\m@th\hfil#1#2\hfil$\crcr
 \noalign{\nointerlineskip}\leftrightarrowfill@#1\crcr}}}
\message{various kinds of dots,}
\let\DOTSI\relax
\let\DOTSB\relax

\newif\ifmath@
{\uccode`7=`\\ \uccode`8=`m \uccode`9=`a \uccode`0=`t \uccode`!=`h
 \uppercase{\gdef\math@#1#2#3#4#5#6\math@{\global\math@false\ifx 7#1\ifx 8#2%
 \ifx 9#3\ifx 0#4\ifx !#5\xdef\meaning@{#6}\global\math@true\fi\fi\fi\fi\fi}}}
\newif\ifmathch@
{\uccode`7=`c \uccode`8=`h \uccode`9=`\"
 \uppercase{\gdef\mathch@#1#2#3#4#5#6\mathch@{\global\mathch@false
  \ifx 7#1\ifx 8#2\ifx 9#5\global\mathch@true\xdef\meaning@{9#6}\fi\fi\fi}}}
\newcount\classnum@
\def\getmathch@#1.#2\getmathch@{\classnum@#1 \divide\classnum@4096
 \ifcase\number\classnum@\or\or\gdef\thedots@{\dotsb@}\or
 \gdef\thedots@{\dotsb@}\fi}
\newif\ifmathbin@
{\uccode`4=`b \uccode`5=`i \uccode`6=`n
 \uppercase{\gdef\mathbin@#1#2#3{\relaxnext@
  \DNii@##1\mathbin@{\ifx\space@\next\global\mathbin@true\fi}%
 \global\mathbin@false\DN@##1\mathbin@{}%
 \ifx 4#1\ifx 5#2\ifx 6#3\DN@{\FN@\nextii@}\fi\fi\fi\next@}}}
\newif\ifmathrel@
{\uccode`4=`r \uccode`5=`e \uccode`6=`l
 \uppercase{\gdef\mathrel@#1#2#3{\relaxnext@
  \DNii@##1\mathrel@{\ifx\space@\next\global\mathrel@true\fi}%
 \global\mathrel@false\DN@##1\mathrel@{}%
 \ifx 4#1\ifx 5#2\ifx 6#3\DN@{\FN@\nextii@}\fi\fi\fi\next@}}}
\newif\ifmacro@
{\uccode`5=`m \uccode`6=`a \uccode`7=`c
 \uppercase{\gdef\macro@#1#2#3#4\macro@{\global\macro@false
  \ifx 5#1\ifx 6#2\ifx 7#3\global\macro@true
  \xdef\meaning@{\macro@@#4\macro@@}\fi\fi\fi}}}
\def\macro@@#1->#2\macro@@{#2}
\newif\ifDOTS@
\newcount\DOTSCASE@
{\uccode`6=`\\ \uccode`7=`D \uccode`8=`O \uccode`9=`T \uccode`0=`S
 \uppercase{\gdef\DOTS@#1#2#3#4#5{\global\DOTS@false\DN@##1\DOTS@{}%
  \ifx 6#1\ifx 7#2\ifx 8#3\ifx 9#4\ifx 0#5\let\next@\DOTS@@\fi\fi\fi\fi\fi
  \next@}}}
{\uccode`3=`B \uccode`4=`I \uccode`5=`X
 \uppercase{\gdef\DOTS@@#1{\relaxnext@
  \DNii@##1\DOTS@{\ifx\space@\next\global\DOTS@true\fi}%
  \DN@{\FN@\nextii@}%
  \ifx 3#1\global\DOTSCASE@\z@\else
  \ifx 4#1\global\DOTSCASE@\@ne\else
  \ifx 5#1\global\DOTSCASE@\tw@\else\DN@##1\DOTS@{}%
  \fi\fi\fi\next@}}}
\newif\ifnot@
{\uccode`5=`\\ \uccode`6=`n \uccode`7=`o \uccode`8=`t
 \uppercase{\gdef\not@#1#2#3#4{\relaxnext@
  \DNii@##1\not@{\ifx\space@\next\global\not@true\fi}%
 \global\not@false\DN@##1\not@{}%
 \ifx 5#1\ifx 6#2\ifx 7#3\ifx 8#4\DN@{\FN@\nextii@}\fi\fi\fi
 \fi\next@}}}
\newif\ifkeybin@
\def\keybin@{\keybin@true
 \ifx\next+\else\ifx\next=\else\ifx\next<\else\ifx\next>\else\ifx\next-\else
 \ifx\next*\else\ifx\next:\else\keybin@false\fi\fi\fi\fi\fi\fi\fi}
\def\dots{\RIfM@\expandafter\mdots@\else\expandafter\tdots@\fi}
\def\tdots@{\unskip\relaxnext@
 \DN@{$\m@th\mathinner{\ldotp\ldotp\ldotp}\,
   \ifx\next,\,$\else\ifx\next.\,$\else\ifx\next;\,$\else\ifx\next:\,$\else
   \ifx\next?\,$\else\ifx\next!\,$\else$ \fi\fi\fi\fi\fi\fi}%
 \ \FN@\next@}
\def\mdots@{\FN@\mdots@@}
\def\mdots@@{\gdef\thedots@{\dotso@}
 \ifx\next\boldkey\gdef\thedots@\boldkey{\boldkeydots@}\else                
 \ifx\next\boldsymbol\gdef\thedots@\boldsymbol{\boldsymboldots@}\else       
 \ifx,\next\gdef\thedots@{\dotsc}
 \else\ifx\not\next\gdef\thedots@{\dotsb@}
 \else\keybin@
 \ifkeybin@\gdef\thedots@{\dotsb@}
 \else\xdef\meaning@{\meaning\next..........}\xdef\meaning@@{\meaning@}
  \expandafter\math@\meaning@\math@
  \ifmath@
   \expandafter\mathch@\meaning@\mathch@
   \ifmathch@\expandafter\getmathch@\meaning@\getmathch@\fi                 
  \else\expandafter\macro@\meaning@@\macro@                                 
  \ifmacro@                                                                
   \expandafter\not@\meaning@\not@\ifnot@\gdef\thedots@{\dotsb@}
  \else\expandafter\DOTS@\meaning@\DOTS@
  \ifDOTS@
   \ifcase\number\DOTSCASE@\gdef\thedots@{\dotsb@}%
    \or\gdef\thedots@{\dotsi}\else\fi                                      
  \else\expandafter\math@\meaning@\math@                                   
  \ifmath@\expandafter\mathbin@\meaning@\mathbin@
  \ifmathbin@\gdef\thedots@{\dotsb@}
  \else\expandafter\mathrel@\meaning@\mathrel@
  \ifmathrel@\gdef\thedots@{\dotsb@}
  \fi\fi\fi\fi\fi\fi\fi\fi\fi\fi\fi\fi
 \thedots@}
\def\plainldots@{\mathinner{\ldotp\ldotp\ldotp}}
\def\plaincdots@{\mathinner{\cdotp\cdotp\cdotp}}
\def\dotsi{\!\plaincdots@}
\let\dotsb@\plaincdots@
\newif\ifextra@
\newif\ifrightdelim@
\def\rightdelim@{\global\rightdelim@true                                    
 \ifx\next)\else                                                            
 \ifx\next]\else
 \ifx\next\rbrack\else
 \ifx\next\}\else
 \ifx\next\rbrace\else
 \ifx\next\rangle\else
 \ifx\next\rceil\else
 \ifx\next\rfloor\else
 \ifx\next\rgroup\else
 \ifx\next\rmoustache\else
 \ifx\next\right\else
 \ifx\next\bigr\else
 \ifx\next\biggr\else
 \ifx\next\Bigr\else                                                        
 \ifx\next\Biggr\else\global\rightdelim@false
 \fi\fi\fi\fi\fi\fi\fi\fi\fi\fi\fi\fi\fi\fi\fi}
\def\extra@{%
 \global\extra@false\rightdelim@\ifrightdelim@\global\extra@true            
 \else\ifx\next$\global\extra@true                                          
 \else\xdef\meaning@{\meaning\next..........}
 \expandafter\macro@\meaning@\macro@\ifmacro@                               
 \expandafter\DOTS@\meaning@\DOTS@
 \ifDOTS@
 \ifnum\DOTSCASE@=\tw@\global\extra@true                                    
 \fi\fi\fi\fi\fi}
\newif\ifbold@
\def\dotso@{\relaxnext@
 \ifbold@
  \let\next\delayed@
  \DNii@{\extra@\plainldots@\ifextra@\,\fi}%
 \else
  \DNii@{\DN@{\extra@\plainldots@\ifextra@\,\fi}\FN@\next@}%
 \fi
 \nextii@}
\def\extrap@#1{%
 \ifx\next,\DN@{#1\,}\else
 \ifx\next;\DN@{#1\,}\else
 \ifx\next.\DN@{#1\,}\else\extra@
 \ifextra@\DN@{#1\,}\else
 \let\next@#1\fi\fi\fi\fi\next@}
\def\ldots{\DN@{\extrap@\plainldots@}%
 \FN@\next@}
\def\cdots{\DN@{\extrap@\plaincdots@}%
 \FN@\next@}

\def\dotsc{\relaxnext@
 \DN@{\ifx\next;\plainldots@\,\else
  \ifx\next.\plainldots@\,\else\extra@\plainldots@
  \ifextra@\,\fi\fi\fi}%
 \FN@\next@}
\def\cdot{\mathchar"2201 }

\def\Longleftrightarrow{\DOTSB\Leftarrow\joinrel\Rightarrow}

\def\iff{\DOTSB\;\Longleftrightarrow\;}
\message{special superscripts,}
\def\dddot#1{{\mathop{#1}\limits^{\vbox to-1.4\ex@{\kern-\tw@\ex@
 \hbox{\rm...}\vss}}}}
\def\ddddot#1{{\mathop{#1}\limits^{\vbox to-1.4\ex@{\kern-\tw@\ex@
 \hbox{\rm....}\vss}}}}
\def\sphat{^{\mathchoice{}{}%
 {\,\,\botsmash{\hbox{\lower4\ex@\hbox{$\m@th\widehat{\null}$}}}}%
 {\,\botsmash{\hbox{\lower3\ex@\hbox{$\m@th\hat{\null}$}}}}}}

\def\spacute{^{\!\botsmash{\hbox{\lower\@ne ex\hbox{\'{}}}}}}
\def\spgrave{^{\mathchoice{}{}{}{\!}%
 \botsmash{\hbox{\lower\@ne ex\hbox{\`{}}}}}}
\def\spdot{^{\hbox{\raise\ex@\hbox{\rm.}}}}
\def\spddot{^{\hbox{\raise\ex@\hbox{\rm..}}}}
\def\spdddot{^{\hbox{\raise\ex@\hbox{\rm...}}}}
\def\spddddot{^{\hbox{\raise\ex@\hbox{\rm....}}}}
\def\spbreve{^{\!\botsmash{\hbox{\lower4\ex@\hbox{\u{}}}}}}

\message{\string\text,}
\def\textonlyfont@#1#2{\def#1{\RIfM@
 \Err@{Use \string#1\space only in text}\else#2\fi}}
\textonlyfont@\rm\tenrm
\textonlyfont@\it\tenit
\textonlyfont@\sl\tensl
\textonlyfont@\bf\tenbf
\def\oldnos#1{\RIfM@{\mathcode`\,="013B \fam\@ne#1}\else
 \leavevmode\hbox{$\m@th\mathcode`\,="013B \fam\@ne#1$}\fi}
\def\text{\RIfM@\expandafter\text@\else\expandafter\text@@\fi}
\def\text@@#1{\leavevmode\hbox{#1}}
\def\mathhexbox@#1#2#3{\text{$\m@th\mathchar"#1#2#3$}}
\def\dag{{\mathhexbox@279}}
\def\ddag{{\mathhexbox@27A}}
\def\S{{\mathhexbox@278}}
\def\P{{\mathhexbox@27B}}
\newif\iffirstchoice@
\firstchoice@true
\def\text@#1{\mathchoice
 {\hbox{\everymath{\displaystyle}\def\textfonti{\the\textfont\@ne}%
  \def\textfontii{\the\textfont\tw@}\textdef@@ T#1}}
 {\hbox{\firstchoice@false
  \everymath{\textstyle}\def\textfonti{\the\textfont\@ne}%
  \def\textfontii{\the\textfont\tw@}\textdef@@ T#1}}
 {\hbox{\firstchoice@false
  \everymath{\scriptstyle}\def\textfonti{\the\scriptfont\@ne}%
  \def\textfontii{\the\scriptfont\tw@}\textdef@@ S\rm#1}}
 {\hbox{\firstchoice@false
  \everymath{\scriptscriptstyle}\def\textfonti
  {\the\scriptscriptfont\@ne}%
  \def\textfontii{\the\scriptscriptfont\tw@}\textdef@@ s\rm#1}}}
\def\textdef@@#1{\textdef@#1\rm\textdef@#1\bf\textdef@#1\sl\textdef@#1\it}
\def\rmfam{0}
\def\textdef@#1#2{%
 \DN@{\csname\expandafter\eat@\string#2fam\endcsname}%
 \if S#1\edef#2{\the\scriptfont\next@\relax}%
 \else\if s#1\edef#2{\the\scriptscriptfont\next@\relax}%
 \else\edef#2{\the\textfont\next@\relax}\fi\fi}
\scriptfont\itfam\tenit \scriptscriptfont\itfam\tenit
\scriptfont\slfam\tensl \scriptscriptfont\slfam\tensl
\newif\iftopfolded@
\newif\ifbotfolded@
\def\topfoldedtext{\topfolded@true\botfolded@false\foldedtext@}
\def\botfoldedtext{\botfolded@true\topfolded@false\foldedtext@}
\def\foldedtext{\topfolded@false\botfolded@false\foldedtext@}
\Invalid@\foldedwidth
\def\foldedtext@{\relaxnext@
 \DN@{\ifx\next\foldedwidth\let\next@\nextii@\else
  \DN@{\nextii@\foldedwidth{.3\hsize}}\fi\next@}%
 \DNii@\foldedwidth##1##2{\setbox\z@\vbox
  {\normalbaselines\hsize##1\relax
  \tolerance1600 \noindent\ignorespaces##2}\ifbotfolded@\boxz@\else
  \iftopfolded@\vtop{\unvbox\z@}\else\vcenter{\boxz@}\fi\fi}%
 \FN@\next@}
\message{math font commands,}
\def\bold{\RIfM@\expandafter\bold@\else
 \expandafter\nonmatherr@\expandafter\bold\fi}
\def\bold@#1{{\bold@@{#1}}}
\def\bold@@#1{\fam\bffam\relax#1}
\def\slanted{\RIfM@\expandafter\slanted@\else
 \expandafter\nonmatherr@\expandafter\slanted\fi}
\def\slanted@#1{{\slanted@@{#1}}}
\def\slanted@@#1{\fam\slfam\relax#1}
\def\roman{\RIfM@\expandafter\roman@\else
 \expandafter\nonmatherr@\expandafter\roman\fi}
\def\roman@#1{{\roman@@{#1}}}
\def\roman@@#1{\fam\rmfam\relax#1}
\def\italic{\RIfM@\expandafter\italic@\else
 \expandafter\nonmatherr@\expandafter\italic\fi}
\def\italic@#1{{\italic@@{#1}}}
\def\italic@@#1{\fam\itfam\relax#1}
\def\Cal{\RIfM@\expandafter\Cal@\else
 \expandafter\nonmatherr@\expandafter\Cal\fi}
\def\Cal@#1{{\Cal@@{#1}}}
\def\Cal@@#1{\noaccents@\fam\tw@#1}
\mathchardef\Gamma="0000
\mathchardef\Delta="0001
\mathchardef\Theta="0002
\mathchardef\Lambda="0003
\mathchardef\Xi="0004
\mathchardef\Pi="0005
\mathchardef\Sigma="0006
\mathchardef\Upsilon="0007
\mathchardef\Phi="0008
\mathchardef\Psi="0009
\mathchardef\Omega="000A
\mathchardef\varGamma="0100
\mathchardef\varDelta="0101
\mathchardef\varTheta="0102
\mathchardef\varLambda="0103
\mathchardef\varXi="0104
\mathchardef\varPi="0105
\mathchardef\varSigma="0106
\mathchardef\varUpsilon="0107
\mathchardef\varPhi="0108
\mathchardef\varPsi="0109
\mathchardef\varOmega="010A
\let\alloc@@\alloc@
\def\hexnumber@#1{\ifcase#1 0\or 1\or 2\or 3\or 4\or 5\or 6\or 7\or 8\or
 9\or A\or B\or C\or D\or E\or F\fi}
\def\loadmsam{%
 \font@\tenmsa=msam10
 \font@\sevenmsa=msam7
 \font@\fivemsa=msam5
 \alloc@@8\fam\chardef\sixt@@n\msafam
 \textfont\msafam=\tenmsa
 \scriptfont\msafam=\sevenmsa
 \scriptscriptfont\msafam=\fivemsa
 \edef\next{\hexnumber@\msafam}%
 \mathchardef\dabar@"0\next39
 \edef\dashrightarrow{\mathrel{\dabar@\dabar@\mathchar"0\next4B}}%
 \edef\dashleftarrow{\mathrel{\mathchar"0\next4C\dabar@\dabar@}}%
 \let\dasharrow\dashrightarrow
 \edef\ulcorner{\delimiter"4\next70\next70 }%
 \edef\urcorner{\delimiter"5\next71\next71 }%
 \edef\llcorner{\delimiter"4\next78\next78 }%
 \edef\lrcorner{\delimiter"5\next79\next79 }%
 \edef\yen{{\noexpand\mathhexbox@\next55}}%
 \edef\checkmark{{\noexpand\mathhexbox@\next58}}%
 \edef\circledR{{\noexpand\mathhexbox@\next72}}%
 \edef\maltese{{\noexpand\mathhexbox@\next7A}}%
 \global\let\loadmsam\empty}%
\def\loadmsbm{%
 \font@\tenmsb=msbm10 \font@\sevenmsb=msbm7 \font@\fivemsb=msbm5
 \alloc@@8\fam\chardef\sixt@@n\msbfam
 \textfont\msbfam=\tenmsb
 \scriptfont\msbfam=\sevenmsb \scriptscriptfont\msbfam=\fivemsb
 \global\let\loadmsbm\empty
 }
\def\widehat#1{\ifx\undefined\msbfam \DN@{362}%
  \else \setboxz@h{$\m@th#1$}%
    \edef\next@{\ifdim\wdz@>\tw@ em%
        \hexnumber@\msbfam 5B%
      \else 362\fi}\fi
  \mathaccent"0\next@{#1}}
\def\widetilde#1{\ifx\undefined\msbfam \DN@{365}%
  \else \setboxz@h{$\m@th#1$}%
    \edef\next@{\ifdim\wdz@>\tw@ em%
        \hexnumber@\msbfam 5D%
      \else 365\fi}\fi
  \mathaccent"0\next@{#1}}
\message{\string\newsymbol,}
\def\newsymbol#1#2#3#4#5{\define#1{}%
  \count@#2\relax \advance\count@\m@ne 
 \ifcase\count@
   \ifx\undefined\msafam\loadmsam\fi \let\next@\msafam
 \or \ifx\undefined\msbfam\loadmsbm\fi \let\next@\msbfam
 \else  \Err@{\Invalid@@\string\newsymbol}\let\next@\tw@\fi
 \mathchardef#1="#3\hexnumber@\next@#4#5\space}
\def\UseAMSsymbols{\loadmsam\loadmsbm}
\def\Bbb{\RIfM@\expandafter\Bbb@\else
 \expandafter\nonmatherr@\expandafter\Bbb\fi}
\def\Bbb@#1{{\Bbb@@{#1}}}
\def\Bbb@@#1{\noaccents@\fam\msbfam\relax#1}
\message{bold Greek and bold symbols,}
\def\loadbold{%
 \font@\tencmmib=cmmib10 \font@\sevencmmib=cmmib7 \font@\fivecmmib=cmmib5
 \skewchar\tencmmib'177 \skewchar\sevencmmib'177 \skewchar\fivecmmib'177
 \alloc@@8\fam\chardef\sixt@@n\cmmibfam
 \textfont\cmmibfam\tencmmib
 \scriptfont\cmmibfam\sevencmmib \scriptscriptfont\cmmibfam\fivecmmib
 \font@\tencmbsy=cmbsy10 \font@\sevencmbsy=cmbsy7 \font@\fivecmbsy=cmbsy5
 \skewchar\tencmbsy'60 \skewchar\sevencmbsy'60 \skewchar\fivecmbsy'60
 \alloc@@8\fam\chardef\sixt@@n\cmbsyfam
 \textfont\cmbsyfam\tencmbsy
 \scriptfont\cmbsyfam\sevencmbsy \scriptscriptfont\cmbsyfam\fivecmbsy
 \let\loadbold\empty
}
\def\boldnotloaded#1{\Err@{\ifcase#1\or First\else Second\fi
       bold symbol font not loaded}}
\def\mathchari@#1#2#3{\ifx\undefined\cmmibfam
    \boldnotloaded@\@ne
  \else\mathchar"#1\hexnumber@\cmmibfam#2#3\space \fi}
\def\mathcharii@#1#2#3{\ifx\undefined\cmbsyfam
    \boldnotloaded\tw@
  \else \mathchar"#1\hexnumber@\cmbsyfam#2#3\space\fi}
\edef\bffam@{\hexnumber@\bffam}
\def\boldkey#1{\ifcat\noexpand#1A%
  \ifx\undefined\cmmibfam \boldnotloaded\@ne
  \else {\fam\cmmibfam#1}\fi
 \else
 \ifx#1!\mathchar"5\bffam@21 \else
 \ifx#1(\mathchar"4\bffam@28 \else\ifx#1)\mathchar"5\bffam@29 \else
 \ifx#1+\mathchar"2\bffam@2B \else\ifx#1:\mathchar"3\bffam@3A \else
 \ifx#1;\mathchar"6\bffam@3B \else\ifx#1=\mathchar"3\bffam@3D \else
 \ifx#1?\mathchar"5\bffam@3F \else\ifx#1[\mathchar"4\bffam@5B \else
 \ifx#1]\mathchar"5\bffam@5D \else
 \ifx#1,\mathchari@63B \else
 \ifx#1-\mathcharii@200 \else
 \ifx#1.\mathchari@03A \else
 \ifx#1/\mathchari@03D \else
 \ifx#1<\mathchari@33C \else
 \ifx#1>\mathchari@33E \else
 \ifx#1*\mathcharii@203 \else
 \ifx#1|\mathcharii@06A \else
 \ifx#10\bold0\else\ifx#11\bold1\else\ifx#12\bold2\else\ifx#13\bold3\else
 \ifx#14\bold4\else\ifx#15\bold5\else\ifx#16\bold6\else\ifx#17\bold7\else
 \ifx#18\bold8\else\ifx#19\bold9\else
  \Err@{\string\boldkey\space can't be used with #1}%
 \fi\fi\fi\fi\fi\fi\fi\fi\fi\fi\fi\fi\fi\fi\fi
 \fi\fi\fi\fi\fi\fi\fi\fi\fi\fi\fi\fi\fi\fi}
\def\boldsymbol#1{%
 \DN@{\Err@{You can't use \string\boldsymbol\space with \string#1}#1}%
 \ifcat\noexpand#1A%
   \let\next@\relax
   \ifx\undefined\cmmibfam \boldnotloaded\@ne
   \else {\fam\cmmibfam#1}\fi
 \else
  \xdef\meaning@{\meaning#1.........}%
  \expandafter\math@\meaning@\math@
  \ifmath@
   \expandafter\mathch@\meaning@\mathch@
   \ifmathch@
    \expandafter\boldsymbol@@\meaning@\boldsymbol@@
   \fi
  \else
   \expandafter\macro@\meaning@\macro@
   \expandafter\delim@\meaning@\delim@
   \ifdelim@
    \expandafter\delim@@\meaning@\delim@@
   \else
    \boldsymbol@{#1}%
   \fi
  \fi
 \fi
 \next@}
\def\mathhexboxii@#1#2{\ifx\undefined\cmbsyfam
    \boldnotloaded\tw@
  \else \mathhexbox@{\hexnumber@\cmbsyfam}{#1}{#2}\fi}
\def\boldsymbol@#1{\let\next@\relax\let\next#1%
 \ifx\next\cdot\mathcharii@201 \else
 \ifx\next\prime{{\null\mathcharii@030 \null}}\else
 \ifx\next\lbrack\mathchar"4\bffam@5B \else
 \ifx\next\rbrack\mathchar"5\bffam@5D \else
 \ifx\next\{\mathcharii@466 \else
 \ifx\next\lbrace\mathcharii@466 \else
 \ifx\next\}\mathcharii@567 \else
 \ifx\next\rbrace\mathcharii@567 \else
 \ifx\next\surd{{\mathcharii@170}}\else
 \ifx\next\S{{\mathhexboxii@78}}\else
 \ifx\next\P{{\mathhexboxii@7B}}\else
 \ifx\next\dag{{\mathhexboxii@79}}\else
 \ifx\next\ddag{{\mathhexboxii@7A}}\else
 \DN@{\Err@{You can't use \string\boldsymbol\space with \string#1}#1}%
 \fi\fi\fi\fi\fi\fi\fi\fi\fi\fi\fi\fi\fi}
\def\boldsymbol@@#1.#2\boldsymbol@@{\classnum@#1 \count@@@\classnum@        
 \divide\classnum@4096 \count@\classnum@                                    
 \multiply\count@4096 \advance\count@@@-\count@ \count@@\count@@@           
 \divide\count@@@\@cclvi \count@\count@@                                    
 \multiply\count@@@\@cclvi \advance\count@@-\count@@@                       
 \divide\count@@@\@cclvi                                                    
 \multiply\classnum@4096 \advance\classnum@\count@@                         
 \ifnum\count@@@=\z@                                                        
  \count@"\bffam@ \multiply\count@\@cclvi
  \advance\classnum@\count@
  \DN@{\mathchar\number\classnum@}%
 \else
  \ifnum\count@@@=\@ne                                                      
   \ifx\undefined\cmmibfam \DN@{\boldnotloaded\@ne}%
   \else \count@\cmmibfam \multiply\count@\@cclvi
     \advance\classnum@\count@
     \DN@{\mathchar\number\classnum@}\fi
  \else
   \ifnum\count@@@=\tw@                                                    
     \ifx\undefined\cmbsyfam
       \DN@{\boldnotloaded\tw@}%
     \else
       \count@\cmbsyfam \multiply\count@\@cclvi
       \advance\classnum@\count@
       \DN@{\mathchar\number\classnum@}%
     \fi
  \fi
 \fi
\fi}
\newif\ifdelim@
\newcount\delimcount@
{\uccode`6=`\\ \uccode`7=`d \uccode`8=`e \uccode`9=`l
 \uppercase{\gdef\delim@#1#2#3#4#5\delim@
  {\delim@false\ifx 6#1\ifx 7#2\ifx 8#3\ifx 9#4\delim@true
   \xdef\meaning@{#5}\fi\fi\fi\fi}}}
\def\delim@@#1"#2#3#4#5#6\delim@@{\if#32%
\let\next@\relax
 \ifx\undefined\cmbsyfam \boldnotloaded\@ne
 \else \mathcharii@#2#4#5\space \fi\fi}
\def\vert{\delimiter"026A30C }
\def\Vert{\delimiter"026B30D }
\let\|\Vert
\def\backslash{\delimiter"026E30F }
\def\boldkeydots@#1{\bold@true\let\next=#1\let\delayed@=#1\mdots@@
 \boldkey#1\bold@false}  
\def\boldsymboldots@#1{\bold@true\let\next#1\let\delayed@#1\mdots@@
 \boldsymbol#1\bold@false}
\message{Euler fonts,}
\def\loadeufm{\loadmathfont{eufm}}

\def\frak{\mathfont@\frak}

\def\loadmathfont#1{%
   \expandafter\font@\csname ten#1\endcsname=#110
   \expandafter\font@\csname seven#1\endcsname=#17
   \expandafter\font@\csname five#1\endcsname=#15
   \edef\next{\noexpand\alloc@@8\fam\chardef\sixt@@n
     \expandafter\noexpand\csname#1fam\endcsname}%
   \next
   \textfont\csname#1fam\endcsname \csname ten#1\endcsname
   \scriptfont\csname#1fam\endcsname \csname seven#1\endcsname
   \scriptscriptfont\csname#1fam\endcsname \csname five#1\endcsname
   \expandafter\def\csname #1\expandafter\endcsname\expandafter{%
      \expandafter\mathfont@\csname#1\endcsname}%
 \expandafter\gdef\csname load#1\endcsname{}%
}
\def\mathfont@#1{\RIfM@\expandafter\mathfont@@\expandafter#1\else
  \expandafter\nonmatherr@\expandafter#1\fi}
\def\mathfont@@#1#2{{\mathfont@@@#1{#2}}}
\def\mathfont@@@#1#2{\noaccents@
   \fam\csname\expandafter\eat@\string#1fam\endcsname
   \relax#2}
\message{math accents,}
\def\accentclass@{7}
\def\noaccents@{\def\accentclass@{0}}
\def\makeacc@#1#2{\def#1{\mathaccent"\accentclass@#2 }}
\makeacc@\hat{05E}
\makeacc@\check{014}
\makeacc@\tilde{07E}
\makeacc@\acute{013}
\makeacc@\grave{012}
\makeacc@\dot{05F}
\makeacc@\ddot{07F}
\makeacc@\breve{015}
\makeacc@\bar{016}

\newcount\skewcharcount@
\newcount\familycount@
\def\theskewchar@{\familycount@\@ne
 \global\skewcharcount@\the\skewchar\textfont\@ne                           
 \ifnum\fam>\m@ne\ifnum\fam<16
  \global\familycount@\the\fam\relax
  \global\skewcharcount@\the\skewchar\textfont\the\fam\relax\fi\fi          
 \ifnum\skewcharcount@>\m@ne
  \ifnum\skewcharcount@<128
  \multiply\familycount@256
  \global\advance\skewcharcount@\familycount@
  \global\advance\skewcharcount@28672
  \mathchar\skewcharcount@\else
  \global\skewcharcount@\m@ne\fi\else
 \global\skewcharcount@\m@ne\fi}                                            
\newcount\pointcount@
\def\getpoints@#1.#2\getpoints@{\pointcount@#1 }
\newdimen\accentdimen@
\newcount\accentmu@
\def\dimentomu@{\multiply\accentdimen@ 100
 \expandafter\getpoints@\the\accentdimen@\getpoints@
 \multiply\pointcount@18
 \divide\pointcount@\@m
 \global\accentmu@\pointcount@}
\def\Makeacc@#1#2{\def#1{\RIfM@\DN@{\mathaccent@
 {"\accentclass@#2 }}\else\DN@{\nonmatherr@{#1}}\fi\next@}}
\def\unbracefonts@{\let\Cal@\Cal@@\let\roman@\roman@@\let\bold@\bold@@
 \let\slanted@\slanted@@}
\def\mathaccent@#1#2{\ifnum\fam=\m@ne\xdef\thefam@{1}\else
 \xdef\thefam@{\the\fam}\fi                                                 
 \accentdimen@\z@                                                           
 \setboxz@h{\unbracefonts@$\m@th\fam\thefam@\relax#2$}
 \ifdim\accentdimen@=\z@\DN@{\mathaccent#1{#2}}
  \setbox@ne\hbox{\unbracefonts@$\m@th\fam\thefam@\relax#2\theskewchar@$}
  \setbox\tw@\hbox{$\m@th\ifnum\skewcharcount@=\m@ne\else
   \mathchar\skewcharcount@\fi$}
  \global\accentdimen@\wd@ne\global\advance\accentdimen@-\wdz@
  \global\advance\accentdimen@-\wd\tw@                                     
  \global\multiply\accentdimen@\tw@
  \dimentomu@\global\advance\accentmu@\@ne                                 
 \else\DN@{{\mathaccent#1{#2\mkern\accentmu@ mu}%
    \mkern-\accentmu@ mu}{}}\fi                                             
 \next@}\Makeacc@\Hat{05E}
\Makeacc@\Check{014}
\Makeacc@\Tilde{07E}
\Makeacc@\Acute{013}
\Makeacc@\Grave{012}
\Makeacc@\Dot{05F}
\Makeacc@\Ddot{07F}
\Makeacc@\Breve{015}
\Makeacc@\Bar{016}
\def\Vec{\RIfM@\DN@{\mathaccent@{"017E }}\else
 \DN@{\nonmatherr@\Vec}\fi\next@}
\def\accentedsymbol#1#2{\csname newbox\expandafter\endcsname
  \csname\expandafter\eat@\string#1@box\endcsname
 \expandafter\setbox\csname\expandafter\eat@
  \string#1@box\endcsname\hbox{$\m@th#2$}\define
  #1{\copy\csname\expandafter\eat@\string#1@box\endcsname{}}}
\message{roots,}
\def\sqrt#1{\radical"270370 {#1}}
\let\underline@\underline
\let\overline@\overline
\def\underline#1{\underline@{#1}}
\def\overline#1{\overline@{#1}}
\Invalid@\leftroot
\Invalid@\uproot
\newcount\uproot@
\newcount\leftroot@
\def\root{\relaxnext@
  \DN@{\ifx\next\uproot\let\next@\nextii@\else
   \ifx\next\leftroot\let\next@\nextiii@\else
   \let\next@\plainroot@\fi\fi\next@}%
  \DNii@\uproot##1{\uproot@##1\relax\FN@\nextiv@}%
  \def\nextiv@{\ifx\next\space@\DN@. {\FN@\nextv@}\else
   \DN@.{\FN@\nextv@}\fi\next@.}%
  \def\nextv@{\ifx\next\leftroot\let\next@\nextvi@\else
   \let\next@\plainroot@\fi\next@}%
  \def\nextvi@\leftroot##1{\leftroot@##1\relax\plainroot@}%
   \def\nextiii@\leftroot##1{\leftroot@##1\relax\FN@\nextvii@}%
  \def\nextvii@{\ifx\next\space@
   \DN@. {\FN@\nextviii@}\else
   \DN@.{\FN@\nextviii@}\fi\next@.}%
  \def\nextviii@{\ifx\next\uproot\let\next@\nextix@\else
   \let\next@\plainroot@\fi\next@}%
  \def\nextix@\uproot##1{\uproot@##1\relax\plainroot@}%
  \bgroup\uproot@\z@\leftroot@\z@\FN@\next@}
\def\plainroot@#1\of#2{\setbox\rootbox\hbox{$\m@th\scriptscriptstyle{#1}$}%
 \mathchoice{\r@@t\displaystyle{#2}}{\r@@t\textstyle{#2}}
 {\r@@t\scriptstyle{#2}}{\r@@t\scriptscriptstyle{#2}}\egroup}
\def\r@@t#1#2{\setboxz@h{$\m@th#1\sqrt{#2}$}%
 \dimen@\ht\z@\advance\dimen@-\dp\z@
 \setbox@ne\hbox{$\m@th#1\mskip\uproot@ mu$}\advance\dimen@ 1.667\wd@ne
 \mkern-\leftroot@ mu\mkern5mu\raise.6\dimen@\copy\rootbox
 \mkern-10mu\mkern\leftroot@ mu\boxz@}
\def\boxed#1{\setboxz@h{$\m@th\displaystyle{#1}$}\dimen@.4\ex@
 \advance\dimen@3\ex@\advance\dimen@\dp\z@
 \hbox{\lower\dimen@\hbox{%
 \vbox{\hrule height.4\ex@
 \hbox{\vrule width.4\ex@\hskip3\ex@\vbox{\vskip3\ex@\boxz@\vskip3\ex@}%
 \hskip3\ex@\vrule width.4\ex@}\hrule height.4\ex@}%
 }}}
\message{commutative diagrams,}
\let\ampersand@\relax
\newdimen\minaw@
\minaw@11.11128\ex@
\newdimen\minCDaw@
\minCDaw@2.5pc
\def\minCDarrowwidth#1{\RIfMIfI@\onlydmatherr@\minCDarrowwidth
 \else\minCDaw@#1\relax\fi\else\onlydmatherr@\minCDarrowwidth\fi}
\newif\ifCD@
\def\CD{\bgroup\vspace@\relax\let\ampersand@&\iffalse}\fi
 \CD@true\vcenter\bgroup\Let@\tabskip\z@skip\baselineskip20\ex@
 \lineskip3\ex@\lineskiplimit3\ex@\halign\bgroup
 &\hfill$\m@th##$\hfill\crcr}
\def\endCD{\crcr\egroup\egroup\egroup}
\newdimen\bigaw@
\atdef@>#1>#2>{\ampersand@                                                  
 \setboxz@h{$\m@th\ssize\;{#1}\;\;$}
 \setbox@ne\hbox{$\m@th\ssize\;{#2}\;\;$}
 \setbox\tw@\hbox{$\m@th#2$}
 \ifCD@\global\bigaw@\minCDaw@\else\global\bigaw@\minaw@\fi                 
 \ifdim\wdz@>\bigaw@\global\bigaw@\wdz@\fi
 \ifdim\wd@ne>\bigaw@\global\bigaw@\wd@ne\fi                                
 \ifCD@\enskip\fi                                                           
 \ifdim\wd\tw@>\z@
  \mathrel{\mathop{\hbox to\bigaw@{\rightarrowfill@\displaystyle}}%
    \limits^{#1}_{#2}}
 \else\mathrel{\mathop{\hbox to\bigaw@{\rightarrowfill@\displaystyle}}%
    \limits^{#1}}\fi                                                        
 \ifCD@\enskip\fi                                                          
 \ampersand@}                                                              
\atdef@<#1<#2<{\ampersand@\setboxz@h{$\m@th\ssize\;\;{#1}\;$}%
 \setbox@ne\hbox{$\m@th\ssize\;\;{#2}\;$}\setbox\tw@\hbox{$\m@th#2$}%
 \ifCD@\global\bigaw@\minCDaw@\else\global\bigaw@\minaw@\fi
 \ifdim\wdz@>\bigaw@\global\bigaw@\wdz@\fi
 \ifdim\wd@ne>\bigaw@\global\bigaw@\wd@ne\fi
 \ifCD@\enskip\fi
 \ifdim\wd\tw@>\z@
  \mathrel{\mathop{\hbox to\bigaw@{\leftarrowfill@\displaystyle}}%
       \limits^{#1}_{#2}}\else
  \mathrel{\mathop{\hbox to\bigaw@{\leftarrowfill@\displaystyle}}%
       \limits^{#1}}\fi
 \ifCD@\enskip\fi\ampersand@}
\begingroup
 \catcode`\~=\active \lccode`\~=`\@
 \lowercase{%
  \global\atdef@)#1)#2){~>#1>#2>}
  \global\atdef@(#1(#2({~<#1<#2<}}
\endgroup
\atdef@ A#1A#2A{\llap{$\m@th\vcenter{\hbox
 {$\ssize#1$}}$}\Big\uparrow\rlap{$\m@th\vcenter{\hbox{$\ssize#2$}}$}&&}
\atdef@ V#1V#2V{\llap{$\m@th\vcenter{\hbox
 {$\ssize#1$}}$}\Big\downarrow\rlap{$\m@th\vcenter{\hbox{$\ssize#2$}}$}&&}
\atdef@={&\enskip\mathrel
 {\vbox{\hrule width\minCDaw@\vskip3\ex@\hrule width
 \minCDaw@}}\enskip&}
\atdef@|{\Big\Vert&&}
\atdef@\vert{\Big\Vert&&}
\def\pretend#1\haswidth#2{\setboxz@h{$\m@th\scriptstyle{#2}$}\hbox
 to\wdz@{\hfill$\m@th\scriptstyle{#1}$\hfill}}
\message{poor man's bold,}
\def\pmb{\RIfM@\expandafter\mathpalette\expandafter\pmb@\else
 \expandafter\pmb@@\fi}
\def\pmb@@#1{\leavevmode\setboxz@h{#1}%
   \dimen@-\wdz@
   \kern-.5\ex@\copy\z@
   \kern\dimen@\kern.25\ex@\raise.4\ex@\copy\z@
   \kern\dimen@\kern.25\ex@\box\z@
}
\def\binrel@@#1{\ifdim\wd2<\z@\mathbin{#1}\else\ifdim\wd\tw@>\z@
 \mathrel{#1}\else{#1}\fi\fi}
\newdimen\pmbraise@
\def\pmb@#1#2{\setbox\thr@@\hbox{$\m@th#1{#2}$}%
 \setbox4\hbox{$\m@th#1\mkern.5mu$}\pmbraise@\wd4\relax
 \binrel@{#2}%
 \dimen@-\wd\thr@@
   \binrel@@{%
   \mkern-.8mu\copy\thr@@
   \kern\dimen@\mkern.4mu\raise\pmbraise@\copy\thr@@
   \kern\dimen@\mkern.4mu\box\thr@@
}}
\def\documentstyle#1{\W@{}\input #1.sty\relax}
\message{syntax check,}
\font\dummyft@=dummy
\fontdimen1 \dummyft@=\z@
\fontdimen2 \dummyft@=\z@
\fontdimen3 \dummyft@=\z@
\fontdimen4 \dummyft@=\z@
\fontdimen5 \dummyft@=\z@
\fontdimen6 \dummyft@=\z@
\fontdimen7 \dummyft@=\z@
\fontdimen8 \dummyft@=\z@
\fontdimen9 \dummyft@=\z@
\fontdimen10 \dummyft@=\z@
\fontdimen11 \dummyft@=\z@
\fontdimen12 \dummyft@=\z@
\fontdimen13 \dummyft@=\z@
\fontdimen14 \dummyft@=\z@
\fontdimen15 \dummyft@=\z@
\fontdimen16 \dummyft@=\z@
\fontdimen17 \dummyft@=\z@
\fontdimen18 \dummyft@=\z@
\fontdimen19 \dummyft@=\z@
\fontdimen20 \dummyft@=\z@
\fontdimen21 \dummyft@=\z@
\fontdimen22 \dummyft@=\z@
\def\fontlist@{\\{\tenrm}\\{\sevenrm}\\{\fiverm}\\{\teni}\\{\seveni}%
 \\{\fivei}\\{\tensy}\\{\sevensy}\\{\fivesy}\\{\tenex}\\{\tenbf}\\{\sevenbf}%
 \\{\fivebf}\\{\tensl}\\{\tenit}}
\def\font@#1=#2 {\rightappend@#1\to\fontlist@\font#1=#2 }
\def\dodummy@{{\def\\##1{\global\let##1\dummyft@}\fontlist@}}
\def\nopages@{\output{\setbox\z@\box\@cclv \deadcycles\z@}%
 \alloc@5\toks\toksdef\@cclvi\output}
\let\galleys\nopages@
\newif\ifsyntax@
\newcount\countxviii@
\def\syntax{\syntax@true\dodummy@\countxviii@\count18
 \loop\ifnum\countxviii@>\m@ne\textfont\countxviii@=\dummyft@
 \scriptfont\countxviii@=\dummyft@\scriptscriptfont\countxviii@=\dummyft@
 \advance\countxviii@\m@ne\repeat                                           
 \dummyft@\tracinglostchars\z@\nopages@\frenchspacing\hbadness\@M}
\def\first@#1#2\end{#1}
\def\printoptions{\W@{Do you want S(yntax check),
  G(alleys) or P(ages)?}%
 \message{Type S, G or P, followed by <return>: }%
 \begingroup 
 \endlinechar\m@ne 
 \read\m@ne to\ans@
 \edef\ans@{\uppercase{\def\noexpand\ans@{%
   \expandafter\first@\ans@ P\end}}}%
 \expandafter\endgroup\ans@
 \if\ans@ P
 \else \if\ans@ S\syntax
 \else \if\ans@ G\galleys
 \else\message{? Unknown option: \ans@; using the `pages' option.}%
 \fi\fi\fi}
\def\alloc@#1#2#3#4#5{\global\advance\count1#1by\@ne
 \ch@ck#1#4#2\allocationnumber=\count1#1
 \global#3#5=\allocationnumber
 \ifalloc@\wlog{\string#5=\string#2\the\allocationnumber}\fi}
\def\document{\def\alloclist@{}\def\fontlist@{}}
\let\enddocument\bye

\let\proclaim\undefined
\let\footnote\undefined
\let\=\undefined
\let\>\undefined

\catcode`\@=\active
\message{... finished}
\def\filename{amsppt.sty}
\def\fileversion{2.2}
\def\filedate{2001/08/07}
\expandafter\ifx\csname amsppt.sty\endcsname\endinput
  \expandafter\def\csname amsppt.sty\endcsname{2.2 (2001/08/07)}\fi
\xdef\fileversiontest{\fileversion\space(\filedate)}
\expandafter\ifx\csname\filename\endcsname\fileversiontest
  \message{[already loaded]}\endinput\fi
\expandafter\ifx\csname\filename\endcsname\relax 
  \else\errmessage{Discrepancy in `\filename' file versions:
     version \csname\filename\endcsname\space already loaded, trying
     now to load version \fileversiontest}\fi
\expandafter\xdef\csname\filename\endcsname{%
  \catcode`\noexpand\@=\the\catcode`\@
  \expandafter\gdef\csname\filename\endcsname{%
     \fileversion\space(\filedate)}}

\catcode`\@=11
\def\savecat#1{%
  \expandafter\xdef\csname\string#1savedcat\endcsname{\the\catcode`#1}}
\def\restorecat#1{\catcode`#1=\csname\string#1savedcat\endcsname}
\message{version \fileversion\space(\filedate):}
\expandafter\ifx\csname styname\endcsname\relax
  
\fi
\message{Loading utility definitions,}
\def\identity@#1{#1}
\def\nofrills@@#1{%
 \DN@{#1}%
 \ifx\next\nofrills \let\frills@\eat@
   \expandafter\expandafter\expandafter\next@\expandafter\eat@
  \else \let\frills@\identity@\expandafter\next@\fi}
\def\nofrillscheck#1{\def\nofrills@{\nofrills@@{#1}}%
  \futurelet\next\nofrills@}
\Invalid@\usualspace
\def\addto#1#2{\csname \expandafter\eat@\string#1@\endcsname
  \expandafter{\the\csname \expandafter\eat@\string#1@\endcsname#2}}
\newdimen\bigsize@
\def\big@#1#2{{\hbox{$\left#2\vcenter to#1\bigsize@{}%
  \right.\nulldelimiterspace\z@\m@th$}}}
\def\big{\big@\@ne}
\def\Big{\big@{1.5}}
\def\bigg{\big@\tw@}
\def\Bigg{\big@{2.5}}
\def\raggedcenter@{\leftskip\z@ plus.4\hsize \rightskip\leftskip
  \parfillskip\z@ \parindent\z@ \spaceskip.3333em \xspaceskip.5em
  \pretolerance9999\tolerance9999 \exhyphenpenalty\@M
  \hyphenpenalty\@M \let\\\linebreak}
\def\uppercasetext@#1{%
   {\spaceskip1.3\fontdimen2\the\font plus1.3\fontdimen3\the\font
    \def\ss{SS}\let\i=I\let\j=J\let\ae\AE\let\oe\OE
    \let\o\O\let\aa\AA\let\l\L
    \skipmath@#1$\skipmath@$}}
\def\skipmath@#1$#2${\uppercase{#1}%
  \ifx\skipmath@#2\else$#2$\expandafter\skipmath@\fi}
\def\add@missing#1{\expandafter\ifx\envir@end#1%
  \Err@{You seem to have a missing or misspelled
  \expandafter\string\envir@end ...}%
  \envir@end
\fi}
\newtoks\revert@
\def\envir@stack#1{\toks@\expandafter{\envir@end}%
  \edef\next@{\def\noexpand\envir@end{\the\toks@}%
    \revert@{\the\revert@}}%
  \revert@\expandafter{\next@}%
  \def\envir@end{#1}}
\begingroup
\catcode`\ =11
\gdef\revert@envir#1{\expandafter\ifx\envir@end#1%
\the\revert@%
\else\ifx\envir@end\enddocument \Err@{Extra \string#1}%
\else\expandafter\add@missing\envir@end\revert@envir#1%
\fi\fi}
\xdef\enddocument {\string\enddocument}%
\global\let\envir@end\enddocument 
\endgroup\relax
\def\first@#1#2\end{#1}
\def\true@{TT}
\def\false@{TF}
\def\empty@{}
\begingroup  \catcode`\-=3
\long\gdef\notempty#1{%
  \expandafter\ifx\first@#1-\end-\empty@ \false@\else \true@\fi}
\endgroup
\message{more fonts,}
\font@\tensmc=cmcsc10 \relax
\let\sevenex=\tenex 
\font@\sevenit=cmti7 \relax
\font@\eightrm=cmr8 \relax 
\font@\sixrm=cmr6 \relax 
\font@\eighti=cmmi8 \relax     \skewchar\eighti='177 
\font@\sixi=cmmi6 \relax       \skewchar\sixi='177   
\font@\eightsy=cmsy8 \relax    \skewchar\eightsy='60 
\font@\sixsy=cmsy6 \relax      \skewchar\sixsy='60   
\let\eightex=\tenex 
\font@\eightbf=cmbx8 \relax 
\font@\sixbf=cmbx6 \relax   
\font@\eightit=cmti8 \relax 
\font@\eightsl=cmsl8 \relax 
\font@\eighttt=cmtt8 \relax 
\let\eightsmc=\nullfont 
\newtoks\tenpoint@
\def\tenpoint{\normalbaselineskip12\p@
 \abovedisplayskip12\p@ plus3\p@ minus9\p@
 \belowdisplayskip\abovedisplayskip
 \abovedisplayshortskip\z@ plus3\p@
 \belowdisplayshortskip7\p@ plus3\p@ minus4\p@
 \textonlyfont@\rm\tenrm \textonlyfont@\it\tenit
 \textonlyfont@\sl\tensl \textonlyfont@\bf\tenbf
 \textonlyfont@\smc\tensmc \textonlyfont@\tt\tentt
 \ifsyntax@ \def\big##1{{\hbox{$\left##1\right.$}}}%
  \let\Big\big \let\bigg\big \let\Bigg\big
 \else
   \textfont\z@\tenrm  \scriptfont\z@\sevenrm
       \scriptscriptfont\z@\fiverm
   \textfont\@ne\teni  \scriptfont\@ne\seveni
       \scriptscriptfont\@ne\fivei
   \textfont\tw@\tensy \scriptfont\tw@\sevensy
       \scriptscriptfont\tw@\fivesy
   \textfont\thr@@\tenex \scriptfont\thr@@\sevenex
        \scriptscriptfont\thr@@\sevenex
   \textfont\itfam\tenit \scriptfont\itfam\sevenit
        \scriptscriptfont\itfam\sevenit
   \textfont\bffam\tenbf \scriptfont\bffam\sevenbf
        \scriptscriptfont\bffam\fivebf
   \setbox\strutbox\hbox{\vrule height8.5\p@ depth3.5\p@ width\z@}%
   \setbox\strutbox@\hbox{\lower.5\normallineskiplimit\vbox{%
        \kern-\normallineskiplimit\copy\strutbox}}%
   \setbox\z@\vbox{\hbox{$($}\kern\z@}\bigsize@1.2\ht\z@
  \fi
  \normalbaselines\rm\dotsspace@1.5mu\ex@.2326ex\jot3\ex@
  \the\tenpoint@}
\newtoks\eightpoint@
\def\eightpoint{\normalbaselineskip10\p@
 \abovedisplayskip10\p@ plus2.4\p@ minus7.2\p@
 \belowdisplayskip\abovedisplayskip
 \abovedisplayshortskip\z@ plus2.4\p@
 \belowdisplayshortskip5.6\p@ plus2.4\p@ minus3.2\p@
 \textonlyfont@\rm\eightrm \textonlyfont@\it\eightit
 \textonlyfont@\sl\eightsl \textonlyfont@\bf\eightbf
 \textonlyfont@\smc\eightsmc \textonlyfont@\tt\eighttt
 \ifsyntax@\def\big##1{{\hbox{$\left##1\right.$}}}%
  \let\Big\big \let\bigg\big \let\Bigg\big
 \else
  \textfont\z@\eightrm \scriptfont\z@\sixrm
       \scriptscriptfont\z@\fiverm
  \textfont\@ne\eighti \scriptfont\@ne\sixi
       \scriptscriptfont\@ne\fivei
  \textfont\tw@\eightsy \scriptfont\tw@\sixsy
       \scriptscriptfont\tw@\fivesy
  \textfont\thr@@\eightex \scriptfont\thr@@\sevenex
   \scriptscriptfont\thr@@\sevenex
  \textfont\itfam\eightit \scriptfont\itfam\sevenit
   \scriptscriptfont\itfam\sevenit
  \textfont\bffam\eightbf \scriptfont\bffam\sixbf
   \scriptscriptfont\bffam\fivebf
 \setbox\strutbox\hbox{\vrule height7\p@ depth3\p@ width\z@}%
 \setbox\strutbox@\hbox{\raise.5\normallineskiplimit\vbox{%
   \kern-\normallineskiplimit\copy\strutbox}}%
 \setbox\z@\vbox{\hbox{$($}\kern\z@}\bigsize@1.2\ht\z@
 \fi
 \normalbaselines\eightrm\dotsspace@1.5mu\ex@.2326ex\jot3\ex@
 \the\eightpoint@}
\def\linespacing#1{%
  \addto\tenpoint{\normalbaselineskip=#1\normalbaselineskip
    \normalbaselines
    \setbox\strutbox=\hbox{\vrule height.7\normalbaselineskip
      depth.3\normalbaselineskip width\z@}%
    \setbox\strutbox@\hbox{\raise.5\normallineskiplimit
      \vbox{\kern-\normallineskiplimit\copy\strutbox}}%
  }%
  \addto\eightpoint{\normalbaselineskip=#1\normalbaselineskip
    \normalbaselines
    \setbox\strutbox=\hbox{\vrule height.7\normalbaselineskip
      depth.3\normalbaselineskip width\z@}%
    \setbox\strutbox@\hbox{\raise.5\normallineskiplimit
      \vbox{\kern-\normallineskiplimit\copy\strutbox}}%
  }%
}
\def\extrafont@#1#2#3{\font#1=#2#3\relax}
\newif\ifPSAMSFonts
\def\PSAMSFonts{%
  \def\extrafont@##1##2##3{%
    \font##1=##2%
      \ifnum##3=9 10 at9pt%
      \else\ifnum##3=8 10 at8pt%
      \else\ifnum##3=6 7 at6pt%
              \else ##3\fi\fi\fi\relax}%
  \font@\eightsmc=cmcsc10 at 8pt
  \font@\eightex=cmex10 at 8pt
  \font@\sevenex=cmex10 at 7pt
  \PSAMSFontstrue
}
\def\loadextrasizes@#1#2#3#4#5#6#7{%
 \ifx\undefined#1%
 \else \extrafont@{#4}{#2}{8}\extrafont@{#6}{#2}{6}%
   \ifsyntax@
   \else
     \addto\tenpoint{\textfont#1#3\scriptfont#1#5%
       \scriptscriptfont#1#7}%
    \addto\eightpoint{\textfont#1#4\scriptfont#1#6%
       \scriptscriptfont#1#7}%
   \fi
 \fi
}
\newtoks\loadextrafonts@@
\def\loadextrafonts@{%
  \loadextrasizes@\msafam{msam}%
    \tenmsa\eightmsa\sevenmsa\sixmsa\fivemsa
  \loadextrasizes@\msbfam{msbm}%
    \tenmsb\eightmsb\sevenmsb\sixmsb\fivemsb
  \loadextrasizes@\eufmfam{eufm}%
    \teneufm\eighteufm\seveneufm\sixeufm\fiveeufm
  \loadextrasizes@\eufbfam{eufb}%
    \teneufb\eighteufb\seveneufb\sixeufb\fiveeufb
  \loadextrasizes@\eusmfam{eusm}%
    \teneusm\eighteusm\seveneusm\sixeusm\fiveeusm
  \loadextrasizes@\eusbfam{eusb}%
    \teneusb\eighteusb\seveneusb\sixeusb\fiveeusb
  \loadextrasizes@\eurmfam{eurm}%
    \teneurm\eighteurm\seveneurm\sixeurm\fiveeurm
  \loadextrasizes@\eurbfam{eurb}%
    \teneurb\eighteurb\seveneurb\sixeurb\fiveeurb
  \loadextrasizes@\cmmibfam{cmmib}%
    \tencmmib\eightcmmib\sevencmmib\sixcmmib\fivecmmib
  \loadextrasizes@\cmbsyfam{cmbsy}%
    \tencmbsy\eightcmbsy\sevencmbsy\sixcmbsy\fivecmbsy
  \let\loadextrafonts@\empty@
  \ifPSAMSFonts
  \else
    \font@\eightsmc=cmcsc8 \relax
    \font@\eightex=cmex8 \relax
    \font@\sevenex=cmex7 \relax
  \fi
  \the\loadextrafonts@@
}
\message{page dimension settings,}
\parindent1pc
\newdimen\normalparindent \normalparindent\parindent
\normallineskiplimit\p@
\newdimen\indenti \indenti=2pc
\let\magnification=\mag
\topskip10pt \relax
\message{top matter,}
\def\topmatter{\loadextrafonts@ \let\topmatter\relax}
\def\chapterno@{\uppercase\expandafter{\romannumeral\chaptercount@}}
\newcount\chaptercount@
\def\chapter{\let\savedef@\chapter
  \def\chapter##1{\let\chapter\savedef@
  \leavevmode\hskip-\leftskip
   \rlap{\vbox to\z@{\vss\centerline{\eightpoint
   \frills@{CHAPTER\space\afterassignment\chapterno@
       \global\chaptercount@=}%
   ##1\unskip}\baselineskip2pc\null}}\hskip\leftskip}%
 \nofrillscheck\chapter}
\newbox\titlebox@
\def\title{\let\savedef@\title
  \def\title##1\endtitle{\let\title\savedef@
    \global\setbox\titlebox@\vtop{\tenpoint\bf
      \raggedcenter@
      \baselineskip1.3\baselineskip
      \frills@\uppercasetext@{##1}\endgraf}%
    \ifmonograph@
      \edef\next{\the\leftheadtoks}\ifx\next\empty@ \leftheadtext{##1}\fi
    \fi
    \edef\next{\the\rightheadtoks}\ifx\next\empty@ \rightheadtext{##1}\fi
  }%
  \nofrillscheck\title}
\newbox\authorbox@
\def\author#1\endauthor{\global\setbox\authorbox@
  \vbox{\tenpoint\smc\raggedcenter@ #1\endgraf}\relaxnext@
  \edef\next{\the\leftheadtoks}%
  \ifx\next\empty@\leftheadtext{#1}\fi}
\newbox\affilbox@
\def\affil#1\endaffil{\global\setbox\affilbox@
  \vbox{\tenpoint\raggedcenter@#1\endgraf}}
\newcount\addresscount@
\addresscount@\z@
\def\addressfont@{\eightpoint}
\def\address#1\endaddress{\global\advance\addresscount@\@ne
  \expandafter\gdef\csname address\number\addresscount@\endcsname
  {\nobreak\vskip12\p@ minus6\p@\indent\addressfont@\smc#1\par}}
\def\curraddr{\let\savedef@\curraddr
  \def\curraddr##1\endcurraddr{\let\curraddr\savedef@
    \if\notempty{##1}%
      \toks@\expandafter\expandafter\expandafter{%
        \csname address\number\addresscount@\endcsname}%
      \toks@@{##1}%
      \expandafter\xdef\csname address\number\addresscount@\endcsname
        {\the\toks@\endgraf\noexpand\nobreak
          \indent\noexpand\addressfont@{\noexpand\rm
          \frills@{{\noexpand\it Current address\noexpand\/}:\space}%
          \def\noexpand\usualspace{\space}\the\toks@@\unskip}}%
    \fi}%
  \nofrillscheck\curraddr}
\def\email{\let\savedef@\email
  \def\email##1\endemail{\let\email\savedef@
    \if\notempty{##1}%
      \toks@{\def\usualspace{{\it\enspace}}\endgraf\indent\addressfont@}%
      \toks@@{{\tt##1}\par}%
      \expandafter\xdef\csname email\number\addresscount@\endcsname
      {\the\toks@\frills@{{\noexpand\it E-mail address\noexpand\/}:%
        \noexpand\enspace}\the\toks@@}%
    \fi}%
  \nofrillscheck\email}

\def\urladdr{\let\savedef@\urladdr
  \def\urladdr##1\endurladdr{\let\urladdr\savedef@
    \if\notempty{##1}%
      \toks@{\def\usualspace{{\it\enspace}}\endgraf\indent\eightpoint}%
      \toks@@{\tt##1\par}%
      \expandafter\xdef\csname urladdr\number\addresscount@\endcsname
      {\the\toks@\frills@{{\noexpand\it URL\noexpand\/}:%
        \noexpand\enspace}\the\toks@@}%
    \fi}%
  \nofrillscheck\urladdr}
\def\thedate@{}
\def\date#1\enddate{\gdef\thedate@{\tenpoint#1\unskip}}
\def\thethanks@{}
\def\thanks#1\endthanks{%
  \if\notempty{#1}%
    \ifx\thethanks@\empty@ \gdef\thethanks@{\eightpoint#1}%
    \else
      \expandafter\gdef\expandafter\thethanks@\expandafter{%
       \thethanks@\endgraf#1}%
    \fi
  \fi}
\def\thekeywords@{}
\def\keywords{\let\savedef@\keywords
  \def\keywords##1\endkeywords{\let\keywords\savedef@
    \if\notempty{##1}%
      \toks@{\def\usualspace{{\it\enspace}}\eightpoint}%
      \toks@@{##1\unskip.}%
      \edef\thekeywords@{\the\toks@\frills@{{\noexpand\it
        Key words and phrases.\noexpand\enspace}}\the\toks@@}%
    \fi}%
  \nofrillscheck\keywords}
\def\xci@{1991}
\def\mm@{2000}
\def\subjclassyear#1{%
   \def\subjyear@{#1}%
   \ifx\subjyear@\mm@
   \else \ifx\subjyear@\xci@
   \else \message{AmS-TeX warning: Unknown edition (#1) of
       Mathematics Subject Classification; using 1991 edition}%
     \def\subjyear@{1991}%
   \fi\fi}
\subjclassyear{1991}
\def\thesubjclass@{}
\def\subjclass{\let\savedef@\subjclass
 \def\subjclass##1\endsubjclass{\let\subjclass\savedef@
   \toks@{\def\usualspace{{\rm\enspace}}\eightpoint}%
   \toks@@{##1\unskip.}%
   \edef\thesubjclass@{\the\toks@
     \frills@{{\noexpand\rm\noexpand\subjyear@\noexpand\space
       {\noexpand\it Mathematics Subject Classification}.\noexpand\enspace}}%
     \the\toks@@}}%
  \nofrillscheck\subjclass}
\newbox\abstractbox@
\def\abstract{\let\savedef@\abstract
  \def\abstract{\let\abstract\savedef@
    \setbox\abstractbox@\vbox\bgroup\noindent$$\vbox\bgroup
      \def\envir@end{\endabstract}\advance\hsize-2\indenti
      \def\usualspace{\enspace}\eightpoint \noindent
      \frills@{{\smc Abstract.\enspace}}}%
  \nofrillscheck\abstract}
\def\endabstract{\par\unskip\egroup$$\egroup}
\def\widestnumber{\begingroup \let\head\relax\let\subhead\relax
  \let\subsubhead\relax \expandafter\endgroup\setwidest@}
\def\setwidest@#1#2{%
   \ifx#1\head\setbox\tocheadbox@\hbox{#2.\enspace}%
   \else\ifx#1\subhead\setbox\tocsubheadbox@\hbox{#2.\enspace}%
   \else\ifx#1\subsubhead\setbox\tocsubheadbox@\hbox{#2.\enspace}%
   \else\ifx#1\key
       \if C\refstyle@ \else\refstyle A\fi
       \setboxz@h{\refsfont@\keyformat{#2}}%
       \refindentwd\wd\z@
   \else\ifx#1\no\refstyle C%
       \setboxz@h{\refsfont@\keyformat{#2}}%
       \refindentwd\wd\z@
   \else\ifx#1\page\setbox\z@\hbox{\quad\bf#2}%
       \pagenumwd\wd\z@
   \else\ifx#1\item
       \setboxz@h{(#2)}\rosteritemwd\wdz@
   \else\message{\string\widestnumber\space not defined for this
      option (\string#1)}%
\fi\fi\fi\fi\fi\fi\fi}
\newif\ifmonograph@
\def\Monograph{\monograph@true \let\headmark\rightheadtext
  \let\varindent@\indent \def\headfont@{\bf}\def\proclaimheadfont@{\smc}%
  \def\remarkheadfont@{\smc}}
\let\varindent@\noindent
\newbox\tocheadbox@    \newbox\tocsubheadbox@
\newbox\tocbox@
\newdimen\pagenumwd
\def\toc{\toc@{Contents}}
\def\newtocdefs{%
   \def \title##1\endtitle
       {\penaltyandskip@\z@\smallskipamount
        \hangindent\wd\tocheadbox@\noindent{\bf##1}}%
   \def \chapter##1{%
        Chapter \uppercase\expandafter{%
              \romannumeral##1.\unskip}\enspace}%
   \def \specialhead##1\endspecialhead
       {\par\hangindent\wd\tocheadbox@ \noindent##1\par}%
   \def \head##1 ##2\endhead
       {\par\hangindent\wd\tocheadbox@ \noindent
        \if\notempty{##1}\hbox to\wd\tocheadbox@{\hfil##1\enspace}\fi
        ##2\par}%
   \def \subhead##1 ##2\endsubhead
       {\par\vskip-\parskip {\normalbaselines
        \advance\leftskip\wd\tocheadbox@
        \hangindent\wd\tocsubheadbox@ \noindent
        \if\notempty{##1}%
              \hbox to\wd\tocsubheadbox@{##1\unskip\hfil}\fi
         ##2\par}}%
   \def \subsubhead##1 ##2\endsubsubhead
       {\par\vskip-\parskip {\normalbaselines
        \advance\leftskip\wd\tocheadbox@
        \hangindent\wd\tocsubheadbox@ \noindent
        \if\notempty{##1}%
              \hbox to\wd\tocsubheadbox@{##1\unskip\hfil}\fi
        ##2\par}}}
\def\toc@#1{\relaxnext@
 \DN@{\ifx\next\nofrills\DN@\nofrills{\nextii@}%
      \else\DN@{\nextii@{{#1}}}\fi
      \next@}%
 \DNii@##1{%
\ifmonograph@\bgroup\else\setbox\tocbox@\vbox\bgroup
   \centerline{\headfont@\ignorespaces##1\unskip}\nobreak
   \vskip\belowheadskip \fi
   \def\page####1%
       {\unskip\penalty\z@\null\hfil
        \rlap{\hbox to\pagenumwd{\quad\hfil####1}}%
              \hfilneg\penalty\@M}%
   \setbox\tocheadbox@\hbox{0.\enspace}%
   \setbox\tocsubheadbox@\hbox{0.0.\enspace}%
   \leftskip\indenti \rightskip\leftskip
   \setboxz@h{\bf\quad000}\pagenumwd\wd\z@
   \advance\rightskip\pagenumwd
   \newtocdefs
 }%
 \FN@\next@}
\def\endtoc{\par\egroup}
\let\pretitle\relax
\let\preauthor\relax
\let\preaffil\relax
\let\predate\relax
\let\preabstract\relax
\let\prepaper\relax
\def\dedicatory #1\enddedicatory{\def\preabstract{{\medskip
  \eightpoint\it \raggedcenter@#1\endgraf}}}
\def\thetranslator@{}
\def\translator{%
  \let\savedef@\translator
  \def\translator##1\endtranslator{\let\translator\savedef@
    \edef\thetranslator@{\noexpand\nobreak\noexpand\medskip
      \noexpand\line{\noexpand\eightpoint\hfil
      \frills@{Translated by \uppercase}{##1}\qquad\qquad}%
       \noexpand\nobreak}}%
  \nofrillscheck\translator}
\outer\def\endtopmatter{\add@missing\endabstract
  \edef\next{\the\leftheadtoks}%
  \ifx\next\empty@
    \expandafter\leftheadtext\expandafter{\the\rightheadtoks}%
  \fi
  \ifmonograph@\else
    \ifx\thesubjclass@\empty@\else \makefootnote@{}{\thesubjclass@}\fi
    \ifx\thekeywords@\empty@\else \makefootnote@{}{\thekeywords@}\fi
    \ifx\thethanks@\empty@\else \makefootnote@{}{\thethanks@}\fi
  \fi
  \inslogo@
  \pretitle
  \begingroup 
  \ifmonograph@ \topskip7pc \else \topskip4pc \fi
  \box\titlebox@
  \endgroup
  \preauthor
  \ifvoid\authorbox@\else \vskip2.5pcplus1pc\unvbox\authorbox@\fi
  \preaffil
  \ifvoid\affilbox@\else \vskip1pcplus.5pc\unvbox\affilbox@\fi
  \predate
  \ifx\thedate@\empty@\else
    \vskip1pcplus.5pc\line{\hfil\thedate@\hfil}\fi
  \preabstract
  \ifvoid\abstractbox@\else
    \vskip1.5pcplus.5pc\unvbox\abstractbox@ \fi
  \ifvoid\tocbox@\else\vskip1.5pcplus.5pc\unvbox\tocbox@\fi
  \prepaper
  \vskip2pcplus1pc\relax
}
\newif\ifdocument@  \document@false
\def\document{\document@true
  \loadextrafonts@
  \let\fontlist@\relax\let\alloclist@\relax
  \tenpoint}
\message{section heads,}
\newskip\aboveheadskip       \aboveheadskip\bigskipamount
\newdimen\belowheadskip      \belowheadskip6\p@
\def\headfont@{\smc}
\def\penaltyandskip@#1#2{\par\skip@#2\relax
  \ifdim\lastskip<\skip@\relax\removelastskip
      \ifnum#1=\z@\else\penalty@#1\relax\fi\vskip\skip@
  \else\ifnum#1=\z@\else\penalty@#1\relax\fi\fi}
\def\nobreak{\penalty\@M
  \ifvmode\gdef\penalty@{\global\let\penalty@\penalty\count@@@}%
  \everypar{\global\let\penalty@\penalty\everypar{}}\fi}
\let\penalty@\penalty
\def\heading#1\endheading{\head#1\endhead}
\def\subheading{\DN@{\ifx\next\nofrills
    \expandafter\subheading@
  \else \expandafter\subheading@\expandafter\empty@
  \fi}%
  \FN@\next@
}
\def\subheading@#1#2{\subhead#1#2\endsubhead}
\newskip\abovespecialheadskip
\abovespecialheadskip=\aboveheadskip
\def\specialheadfont@{\bf}
\outer\def\specialhead{%
  \add@missing\endroster \add@missing\enddefinition
  \add@missing\enddemo \add@missing\endexample
  \add@missing\endproclaim
  \penaltyandskip@{-200}\abovespecialheadskip
  \begingroup\interlinepenalty\@M\rightskip\z@ plus\hsize
  \let\\\linebreak
  \specialheadfont@\noindent}
\def\endspecialhead{\par\endgroup\nobreak\vskip\belowheadskip}
\outer\def\head#1\endhead{%
  \add@missing\endroster \add@missing\enddefinition
  \add@missing\enddemo \add@missing\endexample
  \add@missing\endproclaim
  \penaltyandskip@{-200}\aboveheadskip
  {\headfont@\raggedcenter@\interlinepenalty\@M
  #1\endgraf}\headmark{#1}%
  \nobreak
  \vskip\belowheadskip}
\let\headmark\eat@
\def\restoredef@#1{\relax\let#1\savedef@\let\savedef@\relax}
\newskip\subheadskip       \subheadskip\medskipamount
\def\subheadfont@{\bf}
\outer\def\subhead{%
  \add@missing\endroster \add@missing\enddefinition
  \add@missing\enddemo \add@missing\endexample
  \add@missing\endproclaim
  \let\savedef@\subhead \let\subhead\relax
  \def\subhead##1\endsubhead{\restoredef@\subhead
    \penaltyandskip@{-100}\subheadskip
    {\def\usualspace{\/{\subheadfont@\enspace}}%
     \varindent@\subheadfont@\ignorespaces##1\unskip\frills@{.\enspace}}%
    \ignorespaces}%
  \nofrillscheck\subhead}
\newskip\subsubheadskip       \subsubheadskip\medskipamount
\def\subsubheadfont@{\it}
\outer\def\subsubhead{%
  \add@missing\endroster \add@missing\enddefinition
  \add@missing\enddemo
  \add@missing\endexample \add@missing\endproclaim
  \let\savedef@\subsubhead \let\subsubhead\relax
  \def\subsubhead##1\endsubsubhead{\restoredef@\subsubhead
    \penaltyandskip@{-50}\subsubheadskip
    {\def\usualspace{\/{\subsubheadfont@\enspace}}%
     \subsubheadfont@##1\unskip\frills@{.\enspace}}\ignorespaces}%
  \nofrillscheck\subsubhead}
\message{theorems/proofs/definitions/remarks,}
\def\proclaimheadfont@{\bf}
\def\proclaimfont{\it}
\newskip\preproclaimskip  \preproclaimskip=\medskipamount
\newskip\postproclaimskip \postproclaimskip=\medskipamount
\outer\def\proclaim{%
  \let\savedef@\proclaim \let\proclaim\relax
  \add@missing\endroster \add@missing\enddefinition
  \add@missing\endproclaim \envir@stack\endproclaim
  \def\proclaim##1{\restoredef@\proclaim
    \penaltyandskip@{-100}\preproclaimskip
    {\def\usualspace{\/{\proclaimheadfont@\enspace}}%
     \varindent@\proclaimheadfont@\ignorespaces##1\unskip
     \frills@{.\enspace}}%
    \proclaimfont\ignorespaces}%
  \nofrillscheck\proclaim}
\def\endproclaim{\revert@envir\endproclaim \par\rm
  \vskip\postproclaimskip}
\def\remarkheadfont@{\it}
\def\remarkfont{\rm}
\newskip\remarkskip  \remarkskip=\medskipamount
\def\remark{\let\savedef@\remark \let\remark\relax
  \add@missing\endroster \add@missing\endproclaim
  \envir@stack\endremark
  \def\remark##1{\restoredef@\remark
    \penaltyandskip@\z@\remarkskip
    {\def\usualspace{\/{\remarkheadfont@\enspace}}%
     \varindent@\remarkheadfont@\ignorespaces##1\unskip
     \frills@{.\enspace}}%
    \remarkfont\ignorespaces}%
  \nofrillscheck\remark}
\def\endremark{\par\revert@envir\endremark}
\def\qed{\ifhmode\unskip\nobreak\fi\quad
  \ifmmode\square\else$\m@th\square$\fi}
\newskip\postdemoskip  \postdemoskip=\medskipamount
\newif\if@qedhere
\def\demo{%
  \@qedherefalse
  \DN@{\ifx\next\nofrills
    \DN@####1####2{\remark####1{####2}\envir@stack\enddemo
      \ignorespaces}%
  \else
    \DN@####1{\remark{####1}\envir@stack\enddemo\ignorespaces}%
  \fi
  \next@}%
\FN@\next@}
\def\enddemo{\par\revert@envir\enddemo \endremark\vskip\postdemoskip}
\def\definitionfont{\rm}
\newskip\predefinitionskip  \predefinitionskip=\medskipamount
\newskip\postdefinitionskip  \postdefinitionskip=\medskipamount
\def\definition{\let\savedef@\definition \let\definition\relax
  \add@missing\endproclaim \add@missing\endroster
  \add@missing\enddefinition \envir@stack\enddefinition
  \def\definition##1{\restoredef@\definition
    \penaltyandskip@{-100}\predefinitionskip
    {\def\usualspace{\/{\proclaimheadfont@\enspace}}%
     \varindent@\proclaimheadfont@\ignorespaces##1\unskip
     \frills@{.\proclaimheadfont@\enspace}}%
    \definitionfont\ignorespaces}%
  \nofrillscheck\definition}
\def\enddefinition{\revert@envir\enddefinition
  \par\vskip\postdefinitionskip}
\def\example{%
  \DN@{\ifx\next\nofrills
    \DN@####1####2{\definition####1{####2}\envir@stack\endexample
      \ignorespaces}%
  \else
    \DN@####1{\definition{####1}\envir@stack\endexample\ignorespaces}%
  \fi
  \next@}%
\FN@\next@}
\def\endexample{\revert@envir\endexample \enddefinition }
\message{rosters,}
\newdimen\rosteritemwd
\rosteritemwd16pt 
\newcount\rostercount@
\newif\iffirstitem@
\let\plainitem@\item
\newtoks\everypartoks@
\def\par@{\everypartoks@\expandafter{\the\everypar}\everypar{}}
\def\leftskip@{}
\def\roster{%
  \envir@stack\endroster
  \edef\leftskip@{\leftskip\the\leftskip}%
  \relaxnext@
  \rostercount@\z@
  \def\item{\FN@\rosteritem@}
  \DN@{\ifx\next\runinitem\let\next@\nextii@
    \else\let\next@\nextiii@
    \fi\next@}%
  \DNii@\runinitem
    {\unskip
     \DN@{\ifx\next[\let\next@\nextii@
       \else\ifx\next"\let\next@\nextiii@\else\let\next@\nextiv@\fi
       \fi\next@}%
     \DNii@[####1]{\rostercount@####1\relax
       \enspace\therosteritem{\number\rostercount@}~\ignorespaces}%
     \def\nextiii@"####1"{\enspace{\rm####1}~\ignorespaces}%
     \def\nextiv@{\enspace\therosteritem1\rostercount@\@ne~}%
     \par@\firstitem@false
     \FN@\next@}
  \def\nextiii@{\par\par@
    \penalty\@m\smallskip\vskip-\parskip
    \firstitem@true}%
  \FN@\next@}
\def\rosteritem@{\iffirstitem@\firstitem@false
  \else\par\vskip-\parskip\fi
 \leftskip\rosteritemwd \advance\leftskip\normalparindent
 \advance\leftskip.5em \noindent
 \DNii@[##1]{\rostercount@##1\relax\itembox@}%
 \def\nextiii@"##1"{\def\therosteritem@{\rm##1}\itembox@}%
 \def\nextiv@{\advance\rostercount@\@ne\itembox@}%
 \def\therosteritem@{\therosteritem{\number\rostercount@}}%
 \ifx\next[\let\next@\nextii@\else\ifx\next"\let\next@\nextiii@\else
  \let\next@\nextiv@\fi\fi\next@}
\def\itembox@{\llap{\hbox to\rosteritemwd{\hss
  \kern\z@ 
  \therosteritem@}\enspace}\ignorespaces}
\def\therosteritem#1{\rom{(\ignorespaces#1\unskip)}}
\newif\ifnextRunin@
\def\endroster{\relaxnext@
 \revert@envir\endroster 
 \par\leftskip@ 
 \global\rosteritemwd16\p@ 
 \penalty-50 \vskip-\parskip\smallskip 
 \DN@{\ifx\next\Runinitem\let\next@\relax
  \else\nextRunin@false\let\item\plainitem@ 
   \ifx\next\par 
    \DN@\par{\everypar\expandafter{\the\everypartoks@}}%
   \else 
    \DN@{\noindent\everypar\expandafter{\the\everypartoks@}}%
  \fi\fi\next@}%
 \FN@\next@}
\newcount\rosterhangafter@
\def\Runinitem#1\roster\runinitem{\relaxnext@
  \envir@stack\endroster
 \rostercount@\z@
 \def\item{\FN@\rosteritem@}%
 \def\runinitem@{#1}%
 \DN@{\ifx\next[\let\next\nextii@\else\ifx\next"\let\next\nextiii@
  \else\let\next\nextiv@\fi\fi\next}%
 \DNii@[##1]{\rostercount@##1\relax
  \def\item@{\therosteritem{\number\rostercount@}}\nextv@}%
 \def\nextiii@"##1"{\def\item@{{\rm##1}}\nextv@}%
 \def\nextiv@{\advance\rostercount@\@ne
  \def\item@{\therosteritem{\number\rostercount@}}\nextv@}%
 \def\nextv@{\setbox\z@\vbox
  {\ifnextRunin@\noindent\fi
  \runinitem@\unskip\enspace\item@~\par
  \global\rosterhangafter@\prevgraf}%
  \firstitem@false
  \ifnextRunin@\else\par\fi
  \hangafter\rosterhangafter@\hangindent3\normalparindent
  \ifnextRunin@\noindent\fi
  \runinitem@\unskip\enspace
  \item@~\ifnextRunin@\else\par@\fi
  \nextRunin@true\ignorespaces}
 \FN@\next@}
\message{footnotes,}
\def\footmarkform@#1{$\m@th^{#1}$}
\let\thefootnotemark\footmarkform@
\def\makefootnote@#1#2{\insert\footins
 {\interlinepenalty\interfootnotelinepenalty
 \eightpoint\splittopskip\ht\strutbox\splitmaxdepth\dp\strutbox
 \floatingpenalty\@MM\leftskip\z@skip\rightskip\z@skip
 \spaceskip\z@skip\xspaceskip\z@skip
 \leavevmode{#1}\footstrut\ignorespaces#2\unskip\lower\dp\strutbox
 \vbox to\dp\strutbox{}}}
\newcount\footmarkcount@
\footmarkcount@\z@
\def\footnotemark{\let\@sf\empty@\relaxnext@
 \ifhmode\edef\@sf{\spacefactor\the\spacefactor}\/\fi
 \DN@{\ifx[\next\let\next@\nextii@\else
  \ifx"\next\let\next@\nextiii@\else
  \let\next@\nextiv@\fi\fi\next@}%
 \DNii@[##1]{\footmarkform@{##1}\@sf}%
 \def\nextiii@"##1"{{##1}\@sf}%
 \def\nextiv@{\iffirstchoice@\global\advance\footmarkcount@\@ne\fi
  \footmarkform@{\number\footmarkcount@}\@sf}%
 \FN@\next@}
\def\footnotetext{\relaxnext@
 \DN@{\ifx[\next\let\next@\nextii@\else
  \ifx"\next\let\next@\nextiii@\else
  \let\next@\nextiv@\fi\fi\next@}%
 \DNii@[##1]##2{\makefootnote@{\footmarkform@{##1}}{##2}}%
 \def\nextiii@"##1"##2{\makefootnote@{##1}{##2}}%
 \def\nextiv@##1{\makefootnote@{\footmarkform@%
  {\number\footmarkcount@}}{##1}}%
 \FN@\next@}
\def\footnote{\let\@sf\empty@\relaxnext@
 \ifhmode\edef\@sf{\spacefactor\the\spacefactor}\/\fi
 \DN@{\ifx[\next\let\next@\nextii@\else
  \ifx"\next\let\next@\nextiii@\else
  \let\next@\nextiv@\fi\fi\next@}%
 \DNii@[##1]##2{\footnotemark[##1]\footnotetext[##1]{##2}}%
 \def\nextiii@"##1"##2{\footnotemark"##1"\footnotetext"##1"{##2}}%
 \def\nextiv@##1{\footnotemark\footnotetext{##1}}%
 \FN@\next@}
\def\adjustfootnotemark#1{\advance\footmarkcount@#1\relax}
\def\footnoterule{\kern-4\p@
  \hrule width5pc\kern 3.6\p@}
\message{figures and captions,}
\def\captionfont@{\smc}
\def\topcaption#1#2\endcaption{%
  {\dimen@\hsize \advance\dimen@-\captionwidth@
   \rm\raggedcenter@ \advance\leftskip.5\dimen@ \rightskip\leftskip
  {\captionfont@#1}%
  \if\notempty{#2}\if\notempty{#1}.\enspace\fi\ignorespaces#2\fi
  \endgraf}\nobreak\bigskip}
\def\botcaption#1#2\endcaption{%
  \nobreak\bigskip
  \setboxz@h{\captionfont@#1\if\notempty{#2}\if\notempty{#1}.\enspace\fi
    \rm\ignorespaces#2\fi}%
  {\dimen@\hsize \advance\dimen@-\captionwidth@
   \leftskip.5\dimen@ \rightskip\leftskip
   \noindent \ifdim\wdz@>\captionwidth@
   \else\hfil\fi
  {\captionfont@#1}%
  \if\notempty{#2}\if\notempty{#1}.\enspace\fi\rm\ignorespaces#2\fi\endgraf}}
\def\@ins{\par\begingroup\def\vspace##1{\vskip##1\relax}%
  \def\captionwidth##1{\captionwidth@##1\relax}%
  \setbox\z@\vbox\bgroup} 
\message{miscellaneous,}
\def\block{\RIfMIfI@\nondmatherr@\block\fi
       \else\ifvmode\noindent$$\predisplaysize\hsize
         \else$$\fi
  \def\endblock{\par\egroup$$}\fi
  \vbox\bgroup\advance\hsize-2\indenti\noindent}
\def\endblock{\par\egroup}
\def\cite#1{\rom{[{\citefont@\m@th#1}]}}
\def\citefont@{\rm}
\def\rom#1{\leavevmode
  \edef\prevskip@{\ifdim\lastskip=\z@ \else\hskip\the\lastskip\relax\fi}%
  \unskip
  \edef\prevpenalty@{\ifnum\lastpenalty=\z@ \else
    \penalty\the\lastpenalty\relax\fi}%
  \unpenalty \/\prevpenalty@ \prevskip@ {\rm #1}}
\message{references,}
\def\refsfont@{\eightpoint}
\def\refsheadfont@{\headfont@}
\newdimen\refindentwd
\setboxz@h{\refsfont@ 00.\enspace}
\refindentwd\wdz@
\def\Refsname{References}
\outer\def\Refs{\add@missing\endroster \add@missing\endproclaim
 \let\savedef@\Refs \let\Refs\relax 
 \def\Refs##1{\restoredef@\Refs
   \if\notempty{##1}\penaltyandskip@{-200}\aboveheadskip
     \begingroup \raggedcenter@\refsheadfont@
       \ignorespaces##1\endgraf\endgroup
     \penaltyandskip@\@M\belowheadskip
   \fi
   \begingroup\def\envir@end{\endRefs}\refsfont@\sfcode`\.\@m
   }%
 \nofrillscheck{\csname Refs\expandafter\endcsname
  \frills@{{\Refsname}}}}
\def\endRefs{\par 
  \endgroup}
\newif\ifbook@ \newif\ifprocpaper@
\def\nofrills{%
  \expandafter\ifx\envir@end\endref
    \let\do\relax
    \xdef\nofrills@list{\nofrills@list\do\curbox}%
  \else\errmessage{\Invalid@@ \string\nofrills}%
  \fi}%
\def\defaultreftexts{\gdef\edtext{ed.}\gdef\pagestext{pp.}%
  \gdef\voltext{vol.}\gdef\issuetext{no.}}
\defaultreftexts
\def\ref{\par
  \begingroup \def\envir@end{\endref}%
  \noindent\hangindent\refindentwd
  \def\par{\add@missing\endref}%
  \let\orig@footnote\footnote
  \def\footnote{\message{AmS-TeX warning: \string\footnote\space within
       a reference will disappear;^^J
       use \string\footnotemark\space \string\footnotetext\space instead}%
    \orig@footnote}%
  \global\let\nofrills@list\empty@
  \refbreaks
  \procpaper@false \book@false \moreref@false
  \def\curbox{\z@}\setbox\z@\vbox\bgroup
}
\let\keyhook@\empty@
\def\endref{%
  \setbox\tw@\box\thr@@
  \makerefbox?\thr@@{\endgraf\egroup}%
  \endref@
  \endgraf
  \endgroup
  \keyhook@
  \global\let\keyhook@\empty@ 
}
\def\key{\gdef\key{\makerefbox\key\keybox@\empty@}\key} \newbox\keybox@
\def\no{\gdef\no{\makerefbox\no\keybox@\empty@}%
  \gdef\keyhook@{\refstyle C}\no}
\def\by{\makerefbox\by\bybox@\empty@} \newbox\bybox@
\def\bysame{\by\hbox to3em{\hrulefill}\thinspace\kern\z@}
\def\paper{\makerefbox\paper\paperbox@\it} \newbox\paperbox@
\def\paperinfo{\makerefbox\paperinfo\paperinfobox@\empty@}%
  \newbox\paperinfobox@
\def\jour{\makerefbox\jour\jourbox@
  {\aftergroup\book@false \aftergroup\procpaper@false}} \newbox\jourbox@
\def\issue{\makerefbox\issue\issuebox@\empty@} \newbox\issuebox@
\def\yr{\makerefbox\yr\yrbox@\empty@} \newbox\yrbox@
\def\pages{\makerefbox\pages\pagesbox@\empty@} \newbox\pagesbox@
\def\page{\gdef\pagestext{p.}\makerefbox\page\pagesbox@\empty@}
\def\ed{\makerefbox\ed\edbox@\empty@} \newbox\edbox@
\def\eds{\gdef\edtext{eds.}\makerefbox\eds\edbox@\empty@}
\def\book{\makerefbox\book\bookbox@
  {\it\aftergroup\book@true \aftergroup\procpaper@false}}
  \newbox\bookbox@
\def\bookinfo{\makerefbox\bookinfo\bookinfobox@\empty@}%
  \newbox\bookinfobox@
\def\publ{\makerefbox\publ\publbox@\empty@} \newbox\publbox@
\def\publaddr{\makerefbox\publaddr\publaddrbox@\empty@}%
  \newbox\publaddrbox@
\def\inbook{\makerefbox\inbook\bookbox@
  {\aftergroup\procpaper@true \aftergroup\book@false}}
\def\procinfo{\makerefbox\procinfo\procinfobox@\empty@}%
  \newbox\procinfobox@
\def\finalinfo{\makerefbox\finalinfo\finalinfobox@\empty@}%
  \newbox\finalinfobox@
\def\miscnote{\makerefbox\miscnote\miscnotebox@\empty@}%
  \newbox\miscnotebox@

\def\lang{\makerefbox\lang\langbox@\empty@} \newbox\langbox@
\newbox\morerefbox@
\def\vol{\makerefbox\vol\volbox@{\ifbook@ \else
  \ifprocpaper@\else\bf\fi\fi}}
\newbox\volbox@
\define\MR#1{\makerefbox\MR\MRbox@\empty@
  \def\next@##1:##2:##3\next@{\ifx @##2\empty@##1\else{\bf##1:}##2\fi}%
  MR \next@#1:@:\next@}
\newbox\MRbox@
\define\AMSPPS#1{\makerefbox\AMSPPS\MRbox@\empty@ AMS\-PPS \##1}
\define\CMP#1{\makerefbox\CMP\MRbox@\empty@ CMP #1}
\newbox\holdoverbox
\def\makerefbox#1#2#3{\endgraf
  \setbox\z@\lastbox
  \global\setbox\@ne\hbox{\unhbox\holdoverbox
    \ifvoid\z@\else\unhbox\z@\unskip\unskip\unpenalty\fi}%
  \egroup
  \setbox\curbox\box\ifdim\wd\@ne>\z@ \@ne \else\voidb@x\fi
  \ifvoid#2\else\Err@{Redundant \string#1; duplicate use, or
     mutually exclusive information already given}\fi
  \def\curbox{#2}\setbox\curbox\vbox\bgroup \hsize\maxdimen \noindent
  #3}
\def\refbreaks{%
  \def\refconcat##1{\setbox\z@\lastbox \setbox\holdoverbox\hbox{%
       \unhbox\holdoverbox \unhbox\z@\unskip\unskip\unpenalty##1}}%
  \def\holdover##1{%
    \RIfM@
      \penalty-\@M\null
      \hfil$\clubpenalty\z@\widowpenalty\z@\interlinepenalty\z@
      \offinterlineskip\endgraf
      \setbox\z@\lastbox\unskip \unpenalty
      \refconcat{##1}%
      \noindent
      $\hfil\penalty-\@M
    \else
      \endgraf\refconcat{##1}\noindent
    \fi}%
  \def\break{\holdover{\penalty-\@M}}%
  \let\vadjust@\vadjust
  \def\vadjust##1{\holdover{\vadjust@{##1}}}%
  \def\newpage{\vadjust{\vfill\break}}%
}
\def\refstyle#1{\uppercase{%
  \gdef\refstyle@{#1}%
  \if#1A\relax \def\keyformat##1{[##1]\enspace\hfil}%
  \else\if#1B\relax
    \refindentwd\parindent
    \def\keyformat##1{\aftergroup\kern
              \aftergroup-\aftergroup\refindentwd}%
  \else\if#1C\relax
    \def\keyformat##1{\hfil##1.\enspace}%
  \fi\fi\fi}
}
\refstyle{A}
\def\finalpunct{\ifnum\lastkern=\m@ne\unkern\else.\spacefactor2000 \fi
       \refquotes@\refbreak@}%
\def\continuepunct#1#2#3#4{}%
\def\endref@{%
  \keyhook@
  \def\nofrillscheck##1{%
    \def\do####1{\ifx##1####1\let\frills@\eat@\fi}%
    \let\frills@\identity@ \nofrills@list}%
  \ifvoid\bybox@
    \ifvoid\edbox@
    \else\setbox\bybox@\hbox{\unhbox\edbox@\breakcheck
      \nofrillscheck\edbox@\frills@{\space(\edtext)}\refbreak@}\fi
  \fi
  \ifvoid\keybox@\else\hbox to\refindentwd{%
       \keyformat{\unhbox\keybox@}}\fi
  \ifmoreref@
    \commaunbox@\morerefbox@
  \else
    \kern-\tw@ sp\kern\m@ne sp
  \fi
  \ppunbox@\empty@\empty@\bybox@\empty@
  \ifbook@ 
    \commaunbox@\bookbox@ \commaunbox@\bookinfobox@
    \ppunbox@\empty@{ (}\procinfobox@)%
    \ppunbox@,{ vol.~}\volbox@\empty@
    \ppunbox@\empty@{ (}\edbox@{, \edtext)}%
    \commaunbox@\publbox@ \commaunbox@\publaddrbox@
    \commaunbox@\yrbox@
    \ppunbox@,{ \pagestext~}\pagesbox@\empty@
  \else
    \commaunbox@\paperbox@ \commaunbox@\paperinfobox@
    \ifprocpaper@ 
      \commaunbox@\bookbox@
      \ppunbox@\empty@{ (}\procinfobox@)%
      \ppunbox@\empty@{ (}\edbox@{, \edtext)}%
      \commaunbox@\bookinfobox@
      \ppunbox@,{ \voltext~}\volbox@\empty@
      \commaunbox@\publbox@ \commaunbox@\publaddrbox@
      \commaunbox@\yrbox@
      \ppunbox@,{ \pagestext~}\pagesbox@\empty@
    \else 
      \commaunbox@\jourbox@
      \ppunbox@\empty@{ }\volbox@\empty@
      \ppunbox@\empty@{ (}\yrbox@)%
      \ppunbox@,{ \issuetext~}\issuebox@\empty@
      \commaunbox@\publbox@ \commaunbox@\publaddrbox@
      \commaunbox@\pagesbox@
    \fi
  \fi
  \commaunbox@\finalinfobox@
  \ppunbox@\empty@{ (}\miscnotebox@)%
  \finalpunct
  \ppunbox@\empty@{ (}\langbox@{)\spacefactor1001 }%
  \ifnum\spacefactor>\@m \ppunbox@{}{ }\MRbox@\empty@
  \else \commaunbox@\MRbox@
  \fi
  \defaultreftexts
}
\def\punct@#1{#1}
\def\ppunbox@#1#2#3#4{\ifvoid#3\else
  \let\prespace@\relax
  \ifnum\lastkern=\m@ne \unkern\let\punct@\eat@
    \ifnum\lastkern=-\tw@ \unkern\let\prespace@\ignorespaces \fi
  \fi
  \nofrillscheck#3%
  \punct@{#1}\refquotes@\refbreak@
  \let\punct@\identity@
  \prespace@
  \frills@{#2\eat@}\space
  \unhbox#3\breakcheck
  \frills@{#4\eat@}{\kern\m@ne sp}\fi}
\def\commaunbox@#1{\ppunbox@,\space{#1}\empty@}
\def\breakcheck{\edef\refbreak@{\ifnum\lastpenalty=\z@\else
  \penalty\the\lastpenalty\relax\fi}\unpenalty}
\def\endquotes{\def\refquotes@{''\let\refquotes@\empty@}}
\let\refquotes@\empty@
\let\refbreak@\empty@
\newif\ifmoreref@
\def\moreref{%
  \setbox\tw@\box\thr@@
  \makerefbox?\thr@@{\endgraf\egroup}%
  \let\savedef@\finalpunct  \let\finalpunct\empty@
  \endref@
  \def\punct@##1##2{##2;}%
  \global\let\nofrills@list\empty@ 
  \let\finalpunct\savedef@
  \moreref@true
  \def\curbox{\morerefbox@}%
  \setbox\morerefbox@\vbox\bgroup \hsize\maxdimen \noindent
}

\message{end of document,}
\ifx\plainend\undefined \let\plainend\end \fi
\outer\def\enddocument{\par
  \add@missing\endRefs
  \check@missing@document
  \add@missing\endroster \add@missing\endproclaim
  \add@missing\enddefinition
  \add@missing\enddemo \add@missing\endremark \add@missing\endexample
  \enddocument@text
  \vfill\supereject\plainend}

\def\check@missing@document{%
  \ifdocument@
  \else
    \Err@{You seem to have a missing or misspelled \string\document}%
  \fi}

\def\enddocument@text{%
  \ifmonograph@ 
  \else
    \nobreak
    \thetranslator@
    \count@\z@
    \loop\ifnum\count@<\addresscount@\advance\count@\@ne
      \csname address\number\count@\endcsname
      \csname email\number\count@\endcsname
      \csname urladdr\number\count@\endcsname
    \repeat
  \fi
}

\message{output routine,}
\def\folio{\ifnum\pageno<\z@ \romannumeral-\pageno\else\number\pageno \fi}
\def\foliofont@{\eightrm}
\def\headlinefont@{\eightpoint}
\def\leftheadline{\rlap{\foliofont@\folio}\hfill \iftrue\topmark\fi \hfill}
\def\rightheadline{\hfill \expandafter
  \hfill \llap{\foliofont@\folio}}
\newtoks\leftheadtoks
\newtoks\rightheadtoks
\def\leftheadtext{\let\savedef@\leftheadtext
  \def\leftheadtext##1{\let\leftheadtext\savedef@
    \leftheadtoks\expandafter{\frills@\uppercasetext@{##1}}%
    \mark{\the\leftheadtoks\noexpand\else\the\rightheadtoks}
    \ifsyntax@\setboxz@h{\def\\{\unskip\space\ignorespaces}%
        \headlinefont@##1}\fi}%
  \nofrillscheck\leftheadtext}
\def\rightheadtext{\let\savedef@\rightheadtext
  \def\rightheadtext##1{\let\rightheadtext\savedef@
    \rightheadtoks\expandafter{\frills@\uppercasetext@{##1}}%
    \mark{\the\leftheadtoks\noexpand\else\the\rightheadtoks}%
    \ifsyntax@\setboxz@h{\def\\{\unskip\space\ignorespaces}%
        \headlinefont@##1}\fi}%
  \nofrillscheck\rightheadtext}
\headline={\def\\{\unskip\space\ignorespaces}\headlinefont@
  \def\chapter{%
    \def\chapter##1{%
      \frills@{\afterassignment\chapterno@ \chaptercount@=}##1.\space}%
    \nofrillscheck\chapter}%
  \ifodd\pageno \rightheadline \else \leftheadline\fi}
\def\NoRunningHeads{\global\runheads@false\global\let\headmark\eat@}

\newif\iffirstpage@     \firstpage@true
\newif\ifrunheads@      \runheads@true
\output={\output@}
\newdimen\headlineheight \newdimen\headlinespace
\newdimen\dropfoliodepth

\dropfoliodepth=1pc
\headlineheight=5pt
\headlinespace=24pt

\def\pagewidth#1{\hsize#1%
   \captionwidth@\hsize \advance\captionwidth@-2\indenti}

\def\pageheight#1{%
  \vsize=#1 
  \advance\vsize -\headlineheight 
  \advance\vsize -\headlinespace 
  \advance\vsize \topskip 
}

\pagewidth{30pc}\pageheight{50.5pc}

\newinsert\copyins
\skip\copyins=12\p@
\dimen\copyins=40pc
\count\copyins=1000
\def\inslogo@{\insert\copyins{\logo@}}
\def\logo@{\rightline{\eightpoint Typeset by \AmSTeX}}
\def\nologo{\let\logo@\empty@ \let\inslogo@\empty@}
\let\flheadline\hfil \let\frheadline\hfil
\newif\ifplain@  \plain@false
\def\output@{%
  \def\break{\penalty-\@M}\let\par\endgraf
  \shipout\vbox{%
    \ifplain@
      \let\makeheadline\relax \let\makefootline\relax
    \else
      \iffirstpage@ \global\firstpage@false
        \let\rightheadline\frheadline
        \let\leftheadline\flheadline
      \else
        \ifrunheads@ \let\makefootline\relax
        \else \let\makeheadline\relax \fi
      \fi
    \fi
    \makeheadline \pagebody \makefootline
  }%
  \advancepageno \ifnum\outputpenalty>-\@MM\else\dosupereject\fi
}
\def\pagecontents{%
  \ifvoid\topins\else\unvbox\topins\fi
  \dimen@=\dp\@cclv \unvbox\@cclv 
  \ifvoid\footins
  \else 
    \vskip\skip\footins
    \footnoterule
    \unvbox\footins
  \fi
  \ifr@ggedbottom \kern-\dimen@ \vfil \fi
  \ifvoid\copyins \else \vskip\skip\copyins \unvbox\copyins \fi
}
\def\makeheadline{%
  \leftskip=\z@
  \vbox{%
    \vbox to\headlineheight{\vss
      \hbox to\hsize{\hskip\z@ plus\hsize\the\headline}%
      \kern-\prevdepth
    }%
    \vskip\headlinespace
    \vskip-\topskip
  }%
  \nointerlineskip
}
\def\makefootline{%
  \relax\ifdim\prevdepth>\z@ \ifdim\prevdepth>\maxdepth \else
    \vskip-\prevdepth \fi\fi
  \nointerlineskip
  \vbox to\z@{\hbox{}%
    \baselineskip\dropfoliodepth
    \hbox to\hsize{\hskip\z@ plus\hsize\the\footline}%
    \vss}}

\message{hyphenation exceptions (U.S. English)}
\hyphenation{acad-e-my acad-e-mies af-ter-thought anom-aly anom-alies
an-ti-deriv-a-tive an-tin-o-my an-tin-o-mies apoth-e-o-ses
apoth-e-o-sis ap-pen-dix ar-che-typ-al as-sign-a-ble as-sist-ant-ship
as-ymp-tot-ic asyn-chro-nous at-trib-uted at-trib-ut-able bank-rupt
bank-rupt-cy bi-dif-fer-en-tial blue-print busier busiest
cat-a-stroph-ic cat-a-stroph-i-cally con-gress cross-hatched data-base
de-fin-i-tive de-riv-a-tive dis-trib-ute dri-ver dri-vers eco-nom-ics
econ-o-mist elit-ist equi-vari-ant ex-quis-ite ex-tra-or-di-nary
flow-chart for-mi-da-ble forth-right friv-o-lous ge-o-des-ic
ge-o-det-ic geo-met-ric griev-ance griev-ous griev-ous-ly
hexa-dec-i-mal ho-lo-no-my ho-mo-thetic ideals idio-syn-crasy
in-fin-ite-ly in-fin-i-tes-i-mal ir-rev-o-ca-ble key-stroke
lam-en-ta-ble light-weight mal-a-prop-ism man-u-script mar-gin-al
meta-bol-ic me-tab-o-lism meta-lan-guage me-trop-o-lis
met-ro-pol-i-tan mi-nut-est mol-e-cule mono-chrome mono-pole
mo-nop-oly mono-spline mo-not-o-nous mul-ti-fac-eted mul-ti-plic-able
non-euclid-ean non-iso-mor-phic non-smooth par-a-digm par-a-bol-ic
pa-rab-o-loid pa-ram-e-trize para-mount pen-ta-gon phe-nom-e-non
post-script pre-am-ble pro-ce-dur-al pro-hib-i-tive pro-hib-i-tive-ly
pseu-do-dif-fer-en-tial pseu-do-fi-nite pseu-do-nym qua-drat-ic
quad-ra-ture qua-si-smooth qua-si-sta-tion-ary qua-si-tri-an-gu-lar
quin-tes-sence quin-tes-sen-tial re-arrange-ment rec-tan-gle
ret-ri-bu-tion retro-fit retro-fit-ted right-eous right-eous-ness
ro-bot ro-bot-ics sched-ul-ing se-mes-ter semi-def-i-nite
semi-ho-mo-thet-ic set-up se-vere-ly side-step sov-er-eign spe-cious
spher-oid spher-oid-al star-tling star-tling-ly sta-tis-tics
sto-chas-tic straight-est strange-ness strat-a-gem strong-hold
sum-ma-ble symp-to-matic syn-chro-nous topo-graph-i-cal tra-vers-a-ble
tra-ver-sal tra-ver-sals treach-ery turn-around un-at-tached
un-err-ing-ly white-space wide-spread wing-spread wretch-ed
wretch-ed-ly Brown-ian Eng-lish Euler-ian Feb-ru-ary Gauss-ian
Grothen-dieck Hamil-ton-ian Her-mit-ian Jan-u-ary Japan-ese Kor-te-weg
Le-gendre Lip-schitz Lip-schitz-ian Mar-kov-ian Noe-ther-ian
No-vem-ber Rie-mann-ian Schwarz-schild Sep-tem-ber}
\def\filename{amsppt.sti}
\def\fileversion{2.2}
\def\filedate{2001/08/07}
\loadeufm \loadmsam \loadmsbm
\message{symbol names}\UseAMSsymbols\message{,}




\tenpoint
\W@{}
\csname amsppt.sty\endcsname

\magnification 1200
\vsize = 9.5 true in
\hsize=6.2 true in
\NoRunningHeads        
\parskip=\medskipamount
        \lineskip=2pt\baselineskip=18pt\lineskiplimit=0pt
       
        \TagsOnRight
        \NoBlackBoxes

       \topmatter
        \title
        Eigenfunction localization for the $2D$ periodic Schr\"odinger operator      
        \endtitle
        \author
        W.-M.~Wang
        \endauthor
\address
D\'epartement de Math\'ematique, Universit\'e Paris-Sud, 91405 Orsay Cedex, FRANCE
\endaddress
        \email
{wei-min.wang\@math.u-psud.fr}
\endemail
\abstract
We prove that for any {\it fixed} trigonometric polynomial potential satisfying a genericity condition, the spectrum of the two
dimension periodic Schr\"odinger operator has finite multiplicity and the Fourier series of the eigenfunctions are uniformly exponentially localized about a finite number of frequencies. As a corollary, the $L^p$ norms of the eigenfunctions are bounded for all $p>0$, which answers a question of Toth and Zelditch \cite{TZ}.
\endabstract

        \bigskip\bigskip
        \bigskip
        \toc
        \bigskip
        \bigskip 
        \widestnumber\head {Table of Contents}
        \head 1. Introduction and statement of the theorem\endhead
        \head 2. Partition of the annuli and singular set 
        \endhead
        \head 3. Effective Hamiltonian and reduction to scalar
        \endhead
        \head 4. Polynomial approximation and generic $V$
        \endhead
        \head 5. Proof of the Theorem
        \endhead
        \endtoc
      \endtopmatter
        \vfill\eject
        \bigskip
\document
\head{\bf 1. Introduction and statement of the theorem}\endhead
We consider the Schr\"odinger operator on the square 2-torus $\Bbb T^2$:
$$H=-\frac{\partial^2}{\partial x^2}-\frac{\partial^2}{\partial y^2}+V\tag 1.1$$
on $L^2([-\pi, \pi]^2)$ with periodic boundary condition, where $V$ is real and as a 
function on $\Bbb R^2$ is $2\pi\times 2\pi$ periodic.

Let $V$ be a trigonometric polynomial of degree $k$. Assume $V$ is {\it generic} satisfying the genericity conditions (i, ii) at the end of this section. 
Here it suffices to remark that the genericity condition is explicit, for example $\cos x\cos y$ is generic. Moreover the non generic set is of codimension at least $1$. We postpone the discussion of
generic potentials until then, where we will also show that for $V$ of the form 
$V(x, y)=V_1(x)+V_2(y)$, $H$ always has uniformly bounded multiplicity.

Our main result is 
\proclaim{Theorem}
Let $V$ be a generic trigonometric polynomial of degree $k$. The spectrum of $H$ is of multiplicity at most $Ck^4$ and the Fourier series $\hat\phi$ of the eigenfunctions
$\phi$ with eigenvalues $E$ satisfy
$$|\hat\phi(j)|\leq C \sum_{|\ell|\leq Ck^4}e^{-|j-j_\ell|},\tag 1.2$$
where $C$ is uniform in $E$, while $\{j_\ell\}$ depends on $E$:
$$\{j_\ell\}\subset\{(m,n)\in\Bbb Z^2|\,|m^2+n^2-E|\leq\Vert V\Vert_{\infty}+1\}.\tag 1.3$$
\endproclaim
The above Theorem has the following consequences. Using  (1.2), 
$$\Vert\hat\phi\Vert_{\ell^1}\leq C',\tag 1.4$$
and we have $$\Vert\phi\Vert_{L^\infty}\leq C'.\tag 1.5$$
So we obtain 

\proclaim{Corollary} The eigenfunctions $\phi$ have bounded $L^p$ norms for all $p>0$:
$$\Vert \phi\Vert_{L^p(\Bbb T^2)}<C_p,\qquad \forall p>0.\tag 1.6$$
\endproclaim
\smallskip
\noindent {\it Motivation for the Theorem.}

Our motivation is threefold. The first comes from spectral theory. Consider the Laplacian on the $d$-torus. When $d=1$, the periodic Schr\"odinger operator 
is also called the Hill operator. Its spectral properties are well known. There is an extensive literature on the subject
starting from the 1946 paper of Borg on Sturm-Liouville problems \cite{Bor}.
The main point here is that the equation $n^2=E$ ($E\neq 0$) has only two solutions. The spectrum 
is therefore of multiplicity at most two. When $d>1$, the number of solutions to 
$$n_1^2+n_2^2+\cdots+n_d^2=E\tag 1.7$$
grows with $E$. The spectrum of the Laplacian has unbounded multiplicity. The problem here is therefore basic, namely how to do perturbation theory when there is {\it unbounded degeneracy}. 

Using separation properties of integer solutions to (1.7) and more generally to the inequality:
$$|n_1^2+n_2^2+\cdots+n_d^2-E|\leq A,$$
we prove that when $d=2$ for generic trigonometric polynomial potentials, the spectrum 
of the periodic Schr\"odinger operator in (1.1) has finite multiplicity and the Fourier series of the eigenfunctions are uniformly exponentially localized about a finite number of frequencies, hence
solving a basic problem in spectral theory.

There are previous results on some related problems. For the integrated density of states of the 
corresponding Schr\"odinger operator on $L^2(\Bbb R^2)$, see the recent paper \cite{PS}, cf. also \cite{So}. There are
results on the Schr\"odinger operators when $\Bbb T^2$ is replaced by $\Bbb R^d/\Gamma$, 
where $\Gamma$ is a generic lattice. Hence the multiplicity of the spectrum of the Laplacian
is typically finite \cite{FKT1, 2}.

 A related motivation is the $L^p$ bounds of eigenfunctions on compact manifolds $X$.
 Let $\lambda$ be an eigenvalue of a self-adjoint operator $H$ on $X$. Define 
 $$M_p{\overset\text{def }\to =}\sup_{\Sb\phi\\H\phi=\lambda\phi\endSb}\frac{\Vert\phi\Vert_{L^p}}{\Vert\phi\Vert_{L^2}}.\tag 1.8$$
 
 Assume $\lambda$ has multiplicity $\mu(\lambda)$. Taking $p=\infty$, it is easy to see that 
 $$M_{\infty}\geq\sqrt{\frac{\mu}{\text {vol } X}}.\tag 1.9$$
 by taking the eigenfunction $\psi(x)=\sum_{j=1}^\mu\bar\phi_j(x_0)\phi_j(x)$, where $\{\phi_j\}_{j=1}^\mu$ 
 is an orthonormal basis for the eigenspace corresponding to $\lambda$ and $x_0$ is the point where 
 $\sum_{j=1}^\mu|\phi_j(x_0)|^2\geq\mu/\text{vol } X$. Such an $x_0$ always exists, since 
 $\int \sum_{j=1}^\mu|\phi_j(x)|^2=\mu$. 
 On the other hand, there is the general upper bound from \cite{H, SS}:
 $$M_{\infty}\leq \lambda^{\frac{d-1}{4}}.\tag 1.10$$
 On the sphere (1.9, 1.10) are of the same order, where there is maximal eigenfunction growth. 
 
 On the flat 
 torus $\Bbb T^d=\Bbb R^d/\Bbb Z^d$, (1.10) is far from optimal. For example, when $d=2$ 
 simple number theory consideration gives 
 $$M_\infty\ll \lambda^\epsilon,\qquad\forall \epsilon>0,\tag 1.11$$
 and there are $\lambda$, where  $M_\infty(\lambda)$ are at least logarithmic in $\lambda$.  When $p=4$, 
 Zygmund \cite{Zy} proved however that $M_4(\lambda)\leq 5^{1/4}$.
 
With the addition of a generic polynomial potential $V$ to the Laplacian, the theorem says that on $\Bbb T^2$, 
the Schr\"odinger operator $H$ has finite multiplicity and 
$$M_\infty\leq C.\tag 1.12$$

The corollary answers a question in \cite{TZ, conj. 4.4}, where it is further stipulated that minimal growth criterion
similar to (1.12) characterize flat manifolds under classical integrability conditions. In this context,
see \cite{Bou} for an example where (1.11) is violated by a change of metric. For a general survey on the
subject with connections to number theory and quantum chaos, see \cite{Sa}.

Our motivation for the present problem also comes from parameter dependent situations, e.g., time
dependent or nonlinear perturbations of linear Schr\"odinger equations, where typically the frequency 
(of the perturbation in the linear case and of the quasi-periodic solutions in the nonlinear case) is an 
essential parameter in order to exclude resonances. 

The situation in (1.1) roughly corresponds to the resonant case, where the Theorem shows that there 
is {\it uniform} Fourier restriction and  (1.2) hold. (cf. \cite{W} for a related result in the time dependent case.)  
The small divisors are overcome {\it deterministically} using the separation property of integer solutions to 
$|m^2+n^2-E|\leq \Vert V\Vert_\infty+1$.
\bigskip
\noindent{\it Method of the proof and genericity}

Using the Fourier basis, $H$ is unitarily transformed to a matrix operator $\hat H$:
$$\hat H=\text {diag }(m^2+n^2)+\hat V*$$
on $\ell^2(\Bbb Z^2)$, where $\hat V$ is the Fourier series of $V$. To prove the Theorem, it suffices to
control local eigenvalue spacing. For a given $E$ in the spectrum of $H$:
$\sigma (H)$, we only need to consider the level set $L=\{(m,n)|\,|m^2+n^2-E|\leq \Vert V\Vert_\infty+1\}$, which
is the resonant set.
Using the separation property of $L$ over $\Bbb Z^2$, the local Hamiltonians can be reduced to effective matrices  $\Cal M$ of rank at most $\kappa$, where $\kappa$ is uniform in $E$. 

To investigate $\Cal M$, we first exclude a geometric singular set:
$$\{(m,n)|\,|m\alpha+n\beta|<K, \, (0, 0)\neq (\alpha, \beta)\in\Bbb Z^2,\, |\alpha|,\, |\beta|\leq ck,\, \Bbb Z\ni c,\, K>1\},$$ 
which includes rays determined by the Fourier support of $V$: $\text{supp }\hat V$ . For $\Cal M$ which do not involve resonant sites in the geometric singular set, the sites are at least at a distance $ck$ apart. So we can approximate
$\Cal M$ by a direct sum of scalar ($1\times 1$ matrices) functions. These scalars $M$ correspond to the same function, but 
at different angles $\theta$, which in turn enable us to make $1$ variable polynomial approximations of   
$\det \Cal M$ leading to the genericity conditions on $V$.

We define the geometric support of $\hat V$ to be 
$$\text{gsupp }\hat V{\overset\text{def }\to =}\{(a,b)\in\Bbb Z^2\backslash\{0\}|\exists s\geq 1\text{ such that } (sa, sb)\in\text{supp }\hat V\}, $$
and
$$\align v_{a,b}{\overset\text{def }\to =}&\sum_{1\leq s\leq k}|\hat V(sa,sb)|^2(a^2+b^2),\quad (a,b)\in\text{gsupp }\hat V,
\\
g_{a,b;c,d}{\overset\text{def }\to =}&\sum_{\Sb-k\leq s, s'\leq k\\s,s'\neq 0\endSb}\frac{\hat V(sa,sb)\hat V(s'c-sa,s'd-sb)\hat V(s'c,s'd)}{4s^2s'}, 
\quad (a,b), (c,d) \in\text{gsupp }\hat V.\endalign$$

Let 
$$ \align &f=(1+x^2)\big(\sum_{(a,b)\in\text{gsupp }\hat V} v_{a,b}\frac{ax-b}{(a+bx)^3}
+\sum_{(a,b), (c,d) \in\text{gsupp }\hat V}g_{a,b;c,d}\frac{ax-b}{(a+bx)^2(c+dx)}\big),\\
&\qquad\qquad\qquad\qquad\text{where }a+bx\neq 0, \, c+dx\neq 0 ;\tag 1.13\endalign$$
or $$\align &f=(1+x^2)\big(\sum_{(a,b)\in\text{gsupp }\hat V} v_{a,b}\frac{a-bx}{(ax+b)^3}
+\sum_{(a,b), (c,d) \in\text{gsupp }\hat V}g_{a,b;c,d}\frac{a-bx}{(ax+b)^2(cx+d)}\big),\\
&\qquad\qquad\qquad\qquad\text{where } ax+b\neq 0, \, cx+d\neq 0.\tag 1.14\endalign$$
Both $f$ are rational functions and can be written as 
$$f=\frac{P_1}{P_2}\tag 1.15$$ with $P_1$, $P_2$ polynomials in $x$ of degrees at most $\Cal O(k^4)$ and whose coefficients only depend on $\hat V$ and $\text{supp }\hat V$. 
\smallskip
\noindent {\it Definition.} $V$ is generic, if 
 
\item{(i)} $$\sum_{1\leq s\leq k}|\hat V(sa,sb)|^2(a^2+b^2)-\sum_{1\leq |s|, |s'|\leq k}\frac
{\hat V(sa, sb)\hat V(s'a,s'b)\hat V((s-s')a,(s-s')b))}{4s^2s'}\neq 0.$$
\item{(ii)} $\text {Resultant } (P_1, P_1')\neq 0$
for both $P_1$ defined from (1.13-1.15).

\noindent
(i, ii) show that $\cos x\cos y$ are indeed generic as claimed earlier.

The analysis of the derivative of the scalar function $M$ uses the resolvent expansion. The generic
condition (i) ensures that when $\theta$ is close to the angle of a ray in the geometric support of $\hat V$, the
first two terms dominate and the derivative is away from zero, cf. (4.4, 4.5).

When $\theta$ is otherwise, we make polynomial approximations. The genericity condition (ii)
comes from requiring both $P_1$ to have only simple zeroes so that the excised set contains  at most $\Cal O(k^4)$ 
sites. In studying these polynomials,
we also used a second separation property, namely, if $v_1$ and $v_2$ are two non colinear vectors in 
the Fourier support of $V$, then the angle between them is of order $1$. (For more details, see sect. 4.)

Consequently we show that for generic $V$, there are at most $\Cal O(k^4)$ local eigenvalues which are ``close" to any given $E$. Using again the separation property of the resonant sites in the level sets $L$ over $\Bbb Z^2$ mentioned above, we prove the Theorem. 

We note that when $V$ has separation of variables, i.e., $V(x, y)=V_1(x)+V_2(y)$, the geometric support of 
$\hat V$ is of dimension $2$ and the relevant polynomials can be written in terms of $x=m^2$ only, 
where $m$ is the horizintal coordinate, and are of uniformly bounded degree independent of $k$, cf., (4.13, 4.19- 4.21). It is easy to show that the multiplicity of the eigenvalues are uniformly bounded in $R$ and hence
agree with the known results. 

The scheme presented here to localize an individual eigenfunction is essentially the general one. 
Moreover it is intrinsically independent of self-adjointness. Instead it relies on the geometry of the
Fourier support of $V$, which is more intrinsic. 
In higher dimensions, there are counterparts to the techniques used here, which might be worth pursuing.

I thank P. Sarnak for several stimulating conversations and helpful comments on specific points.
\bigskip
\head{\bf 2. Partition of the annuli and singular set}\endhead
Let 
$$H=-\frac{\partial^2}{\partial x^2}-\frac{\partial^2}{\partial y^2}+V\tag 2.1$$
on $L^2([-\pi, \pi]^2)$ with periodic boundary condition as in section 1. Since $V$ is a real trigonometric polynomial of degree $k$, the Fourier series satisfies
$$\text{supp }\hat V\subset \{(a, b)\in\Bbb Z^2||a|,\, |b|\leq k\}.$$
Without loss, we may assume $\hat V(0, 0)=0$. Otherwise it attributes an overall constant. So
$$\text{supp }\hat V\subset\{(a, b)\in\Bbb Z^2\backslash\{(0,0)\}| |a|,\, |b|\leq k\}.\tag 2.2$$

$H$ is unitarily equivalent to 
$$\hat H=\text {diag }(m^2+n^2)+\hat V*\tag 2.3$$
on $\ell^2(\Bbb Z^2)$, which is the operator that we will work with in the rest of the paper. From now on we write $H$ for $\hat H$ and $\hat V$ for $\hat V*$. 

Assume $E$ is an eigenvalue of $H$, 
$$E\in R+[-1/2, 1/2]\tag 2.4$$
for some $\Bbb Z\ni R\geq -\Vert V\Vert_\infty$.
To deduce localization properties of the eigenfunctions, we only need to be concerned with the 
annulus:
$$S{\overset\text{def }\to =}\{(m,n)\in\Bbb Z^2|\,|m^2+n^2-R|\leq \Vert V\Vert_\infty+1\},\tag 2.5$$
as
$$\Vert (H_{\Bbb Z^2\backslash S}-E)^{-1}\Vert\leq 2,\tag 2.6$$
where for any set $A$, $A\subset\Bbb Z^2$, $H_A$ denotes the restricted operator:
$$\align H_{A}(i, j)&=H(i,j)\quad
(i, j)\in A\times A,\tag 2.7\\
&=0,\quad\qquad\text{otherwise}.\tag 2.8\endalign$$

The following separation property plays a crucial role:
\proclaim{Lemma 2.1} Let $S'$ be the annulus over $\Bbb R^2$: $|x^2+y^2-R|\leq \Vert V\Vert_\infty+1$, $R\in\Bbb N$. There exist $\Bbb N\ni\kappa>0$ (uniform in $R$) and
$\Pi$ a partition of $S'$ such that if $\Bbb R^2\supset p\in\Pi$, then 
$$\align &\bullet\quad |p\cap S'|=\Cal O (R^{1/6}),\tag 2.9\\
&\bullet\quad\#\{p\cap S\}\leq \kappa,\tag 2.10\\
&\bullet\quad\text{dist }( \{p\cap S\},\partial p)=\Cal O (R^{1/6}),\tag 2.11\endalign$$
\noindent where $|\,|$ in (2.9) denotes the length.
\endproclaim
\noindent{\it Remark.} It follows from (2.9, 2.10) that $\# S\leq \Cal O(R^{1/3})$. Estimates on the divisor function give a better bound $\# S\ll R^{\epsilon}$ for all $\epsilon>0$ (leading in particular to (1.11)), but with no geometric information on the integers 
$(m, n)$.

\demo{Proof} We use the argument of Janick \cite {J}, which extends to all stritly convex annuli, cf. \cite{CW}. For completeness we reproduce the proof for the circular annuli. 

We first let  $A_1$, $A_2$ and $A_3$ be $3$ integers (in this order) on the circle $\tilde S$ over $\Bbb R^2$ centered at $O=(0,0)$ of radius $R^{1/2}$,
$A_1\neq A_2\neq A_3$. In view of (2.9), it suffices to assume that $\text {max }(|A_1A_2|, |A_2A_3|)\leq \Cal O(R^{1/6})$.
Since they are not colinear, the area $S_1$ of the triangle formed by $A_1$, $A_2$ and $A_3$ satisfy
$$\Bbb N/2\ni S_1\geq 1/2.\tag 2.12$$
From convexity the area $S_2$ formed by the arc $A_1A_2A_3$ and the straight segment $A_1A_3$ satisfies
$$S_2\geq S_1\geq 1/2.$$
But 
$$S_2=\frac{\theta}{2\pi}\cdot \pi R-\frac{1} {2}R\sin\theta\asymp R\theta^3\geq \frac{1}{2},\tag 2.13$$
where $\theta$ is the angle formed by $OA_1$ and $OA_3$. So
$\theta\geq\Cal O(1/R^{1/3})$ and 
$$|A_1A_3|\geq\Cal O(R^{1/6})\tag 2.14$$ using that the radius is $R^{1/2}$.

We now let  $A_1$, $A_2$ and $A_3$ be any $3$ non colinear integers in $S'$, 
$A_1\neq A_2\neq A_3$ and $\text {max }(|A_1A_2|, |A_2A_3|)\leq \Cal O(R^{1/6})$. (2.12) holds. 
Let $A_j'=OA_j\cap\tilde S$, $j=1,\,2,\,3$. The area $S_1'$ formed by $A_j'$ satisfies:
$S_1'\geq 1/2-\Cal O(1)R^{1/6}\cdot R^{-1/2}>1/4$ using bilinearity. So (2.14) holds. Since the number 
of colinear integers in $S'$ is bounded (uniformly in $R$), (2.9-2.11) follow by choosing the 
$\Cal O(R^{1/6})$ smaller than that in (2.14).
\hfill []
\enddemo

Assume $p\in\Pi$ is such that $p\cap S\neq\emptyset$. Let $H_p$ be defined as in (2.8), where for simplicity 
we also used $p$ to denote $p\cap\Bbb Z^2$. In section 3, we reduce the study of $\sigma(H_p)\cap R+[-1/2, 1/2]$
to that of an effective matrix $\Cal M$, where $\Cal M$ is at most a $\kappa\times \kappa$ matrix. 

To further the analysis, we need to examine the sets $p\in\Pi$.
\proclaim{Lemma 2.2} Let $(x,y)\in\Bbb R^2$ satisfy 
$$|x\alpha+y\beta|\leq K,\qquad (K>0, \text{ independent of } R)\tag 2.15$$
for some $(\alpha,\beta)\in\Bbb Z^2\backslash\{0\}$, $|\alpha|,\,|\beta|\leq ck$, ($\Bbb N\ni c>1$)
and $k$ is the degree of the polynomial $V$. Let 
$$\Pi'=\{p\in\Pi|(2.15) \text{ is violated on } p\cap S'\},\tag 2.16$$
where $S'$ is as defined in Lemma 2.1. Then 
$$|\Pi\backslash\Pi'|\leq 17 c^2k^2.\tag 2.17$$
Assume $K>c^2k^2+\Vert V\Vert_\infty+1$, we have more over that for $p\in\Pi'$, if
$(m,n)$, $(m', n')\in p\cap S'$ satisfying  $(m-m', n-n')\in\Bbb Z^2\backslash\{0\}$, then 
$$\text {sup }  (|m-m'|, |n-n'|)>ck.\tag 2.18$$
\endproclaim
\demo{Proof}
Since $\alpha$ and $\beta$ are integers, (2.15) contains at most $(2ck+1)^2$ tubes $T$ bounded by 
the straight lines 
$$x\alpha+y\beta=\pm K.\tag 2.19$$
Since for each $T$, $T\cap S'$ contains $2$  ``arcs" of length $\Cal O(2K+1)\ll \Cal O(R^{1/6})$, it can intersect
at most $4$ $p\in\Pi$ in view of (2.9), which leads to (2.17). 

To prove (2.18), write $m'-m=\alpha$, $n'-n=\beta$, $(\alpha,\beta)\neq (0, 0)$. We have 
$$
\cases m^2+n^2=R'',\\
(m+\alpha)^2+(n+\beta)^2=R',
\endcases
\tag 2.20
$$
with $|R'-R''|\leq 2(\Vert V\Vert_\infty+1)$.
So $$|m\alpha+n\beta|=\frac{1}{2}|(R'-R'')-(\alpha^2+\beta^2)|.\tag 2.21$$
On the other hand, since $p\in\Pi'$, for all $\Bbb Z^2\ni (\alpha,\beta)\neq (0, 0)$, $|\alpha|$, $|\beta|\leq ck$,
$$|m\alpha+n\beta|>K>c^2k^2+\Vert V\Vert_\infty+1.\tag 2.22$$
(2.21, 2.22) imply that 
$$\text {sup }  (|\alpha|, |\beta|)=\text {sup }  (|m-m'|, |n-n'|)>ck.$$
\hfill []
\enddemo
From now on, $\Pi'$ is to denote the set satisfying (2.16) with $K>c^2k^2+\Vert V\Vert_\infty+1$.

\noindent{\it Remark.} It is important to note that (2.15, 2.17) are {\it independent} of $R$. 
They only depend on the degree $k$ of the trigonometric polynomial $V$. $\Pi\backslash\Pi'$ 
contains the ``singular" set. The effective Hamiltonian reduction will only be used in $\Pi'$,
where the resonant sites are at least at a distance $ck$ apart. 
\bigskip
\head{\bf 3. Effective Hamiltonian and reduction to scalar}\endhead
We now assume $p\in\Pi'$ and use the Schur complement reduction \cite{Sc1, 2} to investigate 
$\sigma (H_p)\cap [R-1/2, R+1/2]$, $R\in\Bbb N$, where $H_p$ is as defined in (2.8). 
Let $S=\{(m,n)\in\Bbb Z^2|\,|m^2+n^2-R|\leq \Vert V\Vert_\infty+1\}$. Assume $p\cap S\neq\emptyset$. Let
$P$ be the projection onto $p\cap S$ and $P_c$ onto $p\backslash S$. 

We have the following equivalence relation:
$$E\in\sigma(H_p)\cap [R-1/2, R+1/2]\iff 0\in\sigma(\Cal M),\tag 3.1$$
where $$\Cal M=E-PH_pP+PH_pP_c(E-P_cH_pP_c)^{-1}P_cH_pP,\tag 3.2$$
cf. \cite{Sect. 2.3, SZ}.
Since $\text{Rank }P\leq \kappa$, $\Cal M$ is at most rank $\kappa$, i.e., a $\kappa\times \kappa$ matrix. Moreover $\Cal M$
is analytic in $E$ for $E\in (R-1/2, R+1/2)$. Since $p\in\Pi'$, in view of (2.18), the
first two terms in (3.2) are diagonal. In the following, we view $E$ as a parameter.

Assume $c>8$ in (2.18). For all $i\in p\cap S$, define:
$$\align \Lambda_0&=i+[-4k-1, 4k+1]^2, \qquad(\text {So }\Lambda_0\cap S=\{i\}.)\tag 3.3\\
M_0&=H_{\Lambda_0\backslash\{i\}},\tag 3.4\\
M'_{ii}&=E-|i|^2+[\hat V(E-M_0)^{-1}\hat V](i, i),\tag 3.5\\
M'&=M'_{ii},\qquad\qquad\qquad\quad \text{ if } |p\cap S|=1,\\
&=\oplus_iM'_{ii}, \quad\qquad\qquad\quad\text{ if } |p\cap S|\geq 2.\tag 3.6\endalign$$
$M'$ is analytic in $E$ for $E\in(R-1/2, R+1/2)$.

\proclaim{Proposition 3.1} For $p\in\Pi'$, 
$$\align \Vert (\Cal M-M')\pi_i\Vert_{\ell^2\to\ell^2}&\leq\Cal O(1)\sum_{(a_\ell, b_\ell)\in\text {supp }\hat V}
\frac{1}{\prod_{\ell=1}^8|ma_\ell+nb_\ell+\frac{a_\ell^2+b_\ell^2}{2}-\frac{\lambda}{2}|}\tag 3.7\\
&<\Cal O(\frac{k^{16}}{K^8}),\tag 3.8\endalign$$
for all $E\in[R-1/2, R+1/2]$, where $i\in p\cap S$, $i=(m,n)$, $\lambda=E-|i|^2$ and $\pi_i$ is the projection onto $\delta_i$, provided $K>c^2k^2+c\Vert V\Vert_\infty+1$ ($\Bbb N\ni c>8$).
\endproclaim
\noindent{\it Remark.} It is important to note that the right side of (3.7) only depends on $i=(m,n)$ and $\text{supp }\hat V$. 

The following lemma is crucial to prove the proposition, in fact to all subsequent analysis. We first define a
few notions. For $(m,n)\in S'=\{(x,y)|\,|x^2+y^2-R|\leq \Vert V\Vert_\infty+1\}$, write $\bar m=(m,n)$. Let 
$$\Bbb Z^2_{\bar m}{\overset\text{def }\to =}\bar m+\Bbb Z^2.\tag 3.9$$
For $j\in \Bbb Z^2_{\bar m}\backslash\{\bar m\}$, define 
$$\align D_{jj}&=E-|j|^2 \tag 3.10\\
D&=\text{diag }D_{jj}.\tag 3.11\endalign$$
Assume $D^{-1}$ exists and define
$$\align F(\bar a)=&\hat V D^{-1}(\bar m, \bar m+\bar a)\\
{\overset\text{def }\to =}&\hat V D^{-1}(\cdot, \cdot+\bar a),\qquad \bar a\in\text{supp }\hat V;\tag 3.12\endalign$$
$$\align F(\bar a_1, \bar a_2)=&\quad\hat V D^{-1}(\cdot, \cdot+\bar a_1)\hat V D^{-1}(\cdot+\bar a_1,  \cdot+\bar a_1+\bar a_2)\\
&+\hat V D^{-1}(\cdot, \cdot+\bar a_2)\hat V D^{-1}(\cdot+\bar a_1, \cdot+\bar a_1+\bar a_2)\\
{\overset\text{def }\to =}&\sum_{\text {perm }(\bar a_1,\bar a_2)}\hat V D^{-1}(\cdot, \cdot+\bar a_1)\hat V D^{-1}(\cdot+\bar a_1,  \cdot+\bar a_1+\bar a_2),\qquad \bar a_1, \bar a_2\in\text{supp }\hat V.\tag 3.13\\
&\vdots\endalign$$
$$\align F(\bar a_1, \bar a_2,\cdots\bar a_s)=\sum_{\text {perm }(\bar a_1,\bar a_2\cdots\bar a_s)}&\hat V D^{-1}(\cdot, \cdot+\bar a_1)\hat V D^{-1}(\cdot+\bar a_1,  \cdot+\bar a_1+\bar a_2)\\
&\cdots \hat V D^{-1}(\cdot+\sum_{\ell=1}^{s-1}\bar a_\ell,  \cdot+\sum_{\ell=1}^{s}\bar a_\ell),\quad \bar a_1, \bar a_2,\cdots
\bar a_s\in\text{supp }\hat V.\tag 3.14\endalign$$

\proclaim{Lemma 3.2} Assume $\bar m=(m,n)\in\Pi'\cap S'$ and increase $K$ to $K>c^2k^2+c\Vert V\Vert_\infty+1$ ($\Bbb N\ni c>8$), so
$|m\alpha+n\beta|>K>c^2k^2+c\Vert V\Vert_\infty+1$, for all $(\alpha,\beta)\in\Bbb Z^2\backslash\{0\}$, $|\alpha|$, $|\beta|\leq ck$.
Then $$|F(\bar a_1, \bar a_2,\cdots\bar a_s)|\leq \Cal O(1)\prod_{\ell=1}^s\frac{|\hat V(a_\ell, b_\ell)|}
{|ma_\ell+nb_\ell+\frac{a_\ell^2+b_\ell^2}{2}-\frac{\lambda}{2}|},\tag 3.15$$
where $\bar a_1,\bar a_2,\cdots\bar a_s\in\text{supp }\hat V,\, 1\leq s\leq c$ and $\lambda=E-m^2-n^2$ as in (3.7).
\endproclaim
\demo{Proof} When $s=1$, (3.15) follows from the definition (3.12). When $s=2$
$$\aligned 4F(\bar a_1, \bar a_2)=&\quad\frac{\hat V(a_1, b_1)\hat V(a_2, b_2)}
{(ma_1+nb_1+\frac{a_1^2+b_1^2}{2}-\frac{\lambda}{2})(m(a_1+a_2)+n(b_1+b_2)+\frac{(a_1+a_2)^2+(b_1+b_2)^2}{2}-\frac{\lambda}{2})}\\
&+\frac{\hat V(a_2, b_2)\hat V(a_1, b_1)}
{(ma_2+nb_2+\frac{a_2^2+b_2^2}{2}-\frac{\lambda}{2})(m(a_1+a_2)+n(b_1+b_2)+\frac{(a_1+a_2)^2+(b_1+b_2)^2}{2}-\frac{\lambda}{2})}.\endaligned\tag 3.16$$
To simplify notation, let 
$$A_\ell=ma_\ell+nb_\ell+\frac{a_\ell^2+b_\ell^2}{2}-\frac{\lambda}{2},\tag 3.17$$
and more generally
$$A_{\ell_1\cdots\ell_s}=m\sum_{\ell=1}^s a_\ell+n\sum_{\ell=1}^s b_\ell+\frac{(\sum_{\ell=1}^s a_\ell)^2+(\sum_{\ell=1}^s b_\ell)^2}{2}-\frac{\lambda}{2}.\tag 3.18$$
So $$\aligned (3.16)&=\hat V(a_1, b_1)\hat V(a_2, b_2)[\frac{1}{A_1A_{12}}+\frac{1}{A_2A_{12}}]\\
&=\hat V(a_1, b_1)\hat V(a_2, b_2)[\frac{1}{A_1A_{2}}+\frac{\Cal O(1)}{A_1A_2A_{12}}].\endaligned\tag 3.19$$
Taking the absolute value, we obtain (3.15) for $s=2$. 

We now make an induction on $s$. Assume
$$F(\bar a_1, \bar a_2,\cdots, \bar a_s)=\prod_{\ell=1}^s \hat V(a_\ell, b_\ell)\cdot [\frac{1}{A_1A_2\cdots A_s}
+\frac{\Cal O(1)}{A_1A_2\cdots A_s A_{1\cdots s}}]\tag 3.20$$
holds at $s$. To arrive at $s+1$, we write
$$F(\bar a_1, \bar a_2,\cdots\bar a_{s+1})=\sum_\sigma F(\sigma)\frac{\hat V(\sigma^c)}{A_{1\cdots s+1}},\tag 3.21$$
where $\sigma\subset\{\bar a_1,\bar a_2,\cdots,\bar a_{s+1}\}$, $|\sigma|=s$, $\sigma^c$ is the complement, which 
only has one element. Using (3.20) for $F(\sigma)$, we obtain 
$$\aligned F(\bar a_1,\bar a_2,&\cdots,\bar a_{s+1})=\prod_{\ell=1}^{s+1}\hat V(a_\ell, b_\ell)\cdot\\
&[\frac{A_1+A_2+\cdots+A_{s+1}}{A_1A_2\cdots A_s A_{s+1} A_{1\cdots s+1}}+
\sum_\sigma\frac{\Cal O(1)}{(\prod_{\ell_s\in\sigma}A_{\ell_s})A_\sigma A_{1\cdots s+1}}].\endaligned\tag 3.22$$
$$A_1+A_2+\cdots+A_{s+1}=A_{1\cdots s+1}+\Cal O(1)\tag 3.23$$ 
and since $$\frac{1}{A_\sigma A_{1\cdots s+1}}=[\frac{1}{A_\sigma}-\frac{1}{A_{1\cdots s+1}}]\cdot
\frac{\Cal O(1)}{A_{\sigma^c}},\tag 3.24$$
(3.22) gives 
$$F(\bar a_1,\bar a_2,\cdots,\bar a_{s+1})=\prod_{\ell=1}^{s+1} \hat V(a_\ell, b_\ell)\cdot[\frac{1}{A_1A_2\cdots A_{s+1}}
+\frac{\Cal O(1)}{A_1A_2\cdots A_{s +1}A_{1\cdots s+1}}].\tag 3.25$$
Increasing $K$ to $K>c^2k^2+c\Vert V\Vert_\infty+1$ in view of the $\Cal O(1)$ in (3.23) and taking the absolute value, we obtain (3.15).\hfill []
\enddemo

\demo{Proof of Proposition 3.1}
Assume $|p\cap S|\geq 2$, otherwise set $\Cal M_{ij}=0$ $(i\neq j, i,j\in\{p\cap S\})$ in the argument below. Let 
$$\aligned &i\in p\cap S,\, \Lambda_0^c=p\backslash\{\Lambda_0\cup\{p\cap S\}\},\\
&M^c=P_cH_pP_c,\, M_0^c=H_{\Lambda_0^c},\\
&\Gamma=M^c-(M_0\oplus M_0^c),\endaligned\tag 3.26$$ 
where $\Lambda_0$ and $M_0$ as defined in (3.3, 3.4).

Using (3.2, 3.5, 3.26) and the resolvent equation, we have 
$$\aligned &\Cal M_{ii}-M'_{ii}=\hat V(E-M_0)^{-1}\Gamma (E-M^c)^{-1}\hat V\\
=&\sum_{i',i''}[\hat V D^{-1}]^4(i,i')[\hat V(E-M_0)^{-1}\Gamma (E-M^c)^{-1}\hat V](i',i'')[D^{-1}\hat V]^4(i'',i),\endaligned\tag 3.27$$
where we used the fact that (3.3) implies 
$$\text{dist }(i,\text{supp }\Gamma)>3k.\tag 3.28$$
Using (3.15) for $s=4$ and 
$$\aligned &\Vert (E-M_0)^{-1}\Vert =\Cal O(1/K),\\
 &\Vert (E-M^c)^{-1}\Vert=\Cal O(1),\\
 &\Vert \Gamma\Vert =\Cal O(1),\endaligned$$
 we obtain 
 $$|\Cal M_{ii}-M'_{ii}|\leq\Cal O(1/K)\cdot \big[\frac{(2k+1)^8}{4!}\big]^2\cdot\text{sup }_{(a_\ell, b_\ell)\in\text{supp }\hat V}
\prod_{\ell=1}^8\frac{1}
{|ma_\ell+nb_\ell+\frac{a_\ell^2+b_\ell^2}{2}-\frac{\lambda}{2}|},\tag 3.29$$
where $(m,n)=i$.  Similarly,
$$\Cal M_{ij}=\sum_{i'}[\hat VD^{-1}]^8(i,i')[\hat V(E-M^c)^{-1}\hat V](i',j),\quad i\neq j,\, i,j\in\{p\cap S\},\tag 3.30$$
where we used $|i-j|_\infty>ck>8k$.
So $$|\Cal M_{ij}|\leq\frac{(2k+1)^{16}}{8!}\text{min }_{(m,n)=i,j}\cdot\text{sup }_{(a_\ell, b_\ell)\in\text{supp }\hat V}\prod_{\ell=1}^8\frac{1}
{|ma_\ell+nb_\ell+\frac{a_\ell^2+b_\ell^2}{2}-\frac{\lambda}{2}|}\tag 3.31$$
using (3.15) and $\Cal M_{ij}=\Cal M_{ji}$. (3.29, 3.31) imply (3.7, 3.8).\hfill []
\enddemo
\smallskip
\noindent{\it The scalar Hamiltonian.}

Let $$\aligned M_{ii}&=M'_{ii}-E+|i|^2\\
&=M'_{ii}-\lambda\\
&=[\hat V(E-M_0)^{-1}\hat V](i, i), \quad i\in p\cap S\subset\Bbb Z^2,\endaligned\tag 3.32$$
from (3.3-3.5). $M_{ii}$ is the scalar Hamiltonian that we will study in detail in section 4. Here it 
suffices to note that for fixed $E$ and $|i|$, $M_{ii}$ is only a function of the angle $\theta$: $M_{ii}=M_{ii}(\theta)$. 
Moreover for $i\in p\cap S'\subset\Bbb R^2$, defining $\Lambda_0$ as in (3.3), $\Lambda_0\subset\Bbb Z_i^2=i+\Bbb Z^2$,
$M_0$ defined in (3.4) extends to a matrix on $\ell^2(\Lambda_0)$. 

For $i$ such that $\Lambda_0\subset p\in\Pi'$, $(E-M_0)^{-1}$ is well defined for $E\in[R-1/2, R+1/2]$ with
$$\Vert (E-M_0)^{-1}\Vert_{\ell^2(\Lambda_0\backslash\{i\})}\leq\Cal O(1/K).\tag 3.33$$
This is because for any two points $(m',n')$, $(m'', n'')\in\Lambda_0$, assuming $(m',n')\neq (m'', n'')$, one 
writes $m''=m'+\alpha$, $n''=n'+\beta$ with $\Bbb Z^2\ni (\alpha,\beta)\neq 0$.  Assuming
$$\aligned &{m'}^2+{n'}^2=R',\\ 
&{m''}^2+{n''}^2=R'',\endaligned\tag 3.34$$
one obtains 
$$|m'\alpha+n'\beta|=\frac{1}{2}|R'-R''-(\alpha^2+\beta^2)|>K>c^2k^2+c\Vert V\Vert_\infty+1,\tag 3.35$$
since $p\in\Pi'$and we used the larger $K$ from Proposition 3.1. 

From (3.3), $\text{sup }(|\alpha|, |\beta|)\leq 8k+2$, so $(\alpha^2+\beta^2)/2\leq (8k+2)^2$. Using this in (3.35) we
have $|R'-R''|>2(K-(8k+2)^2)>2((c^2-65)k^2+c\Vert V\Vert_\infty+1)$. Since $i\in\Lambda_0$ satisfies $||i|^2-R|\leq\Vert V\Vert_\infty+1$ and $\Bbb N\ni c>8$, this proves (3.33).
We now view $M_{ii}$ as a function of $\theta$, defined on appropriate arcs of the circle 
$S_i'=\{(x,y)|\,|x^2+y^2=|i|^2\}\subset\Bbb R^2$.
\bigskip
\head{\bf 4. Polynomial approximation and generic $V$}\endhead
In this section, we investigate the scalar Hamiltonian
$$M_{ii}=[\hat V(E-M_0)^{-1}\hat V](i, i)\tag 4.1$$
as defined in (3.32, 3.3-3.5). From (3.33), $|M_{ii}|\leq\Cal O(1/K)$, where $K>c^2k^2+c\Vert V\Vert_\infty+1$ and 
$\Bbb N\ni c>8$.
Since $E\in[R-1/2, R+1/2]$, for the purpose
of this section, we only need to consider $i$ such that 
$i=(m,n)=(\sqrt R\cos\theta, \sqrt R\sin\theta)$.  Let $\tilde S=\{(m,n)\in\Bbb R^2|m^2+n^2=R\}$
and $D$ the diagonal part of $E-M_0$ as before:
$$D_{\ell\ell}=E-|\ell|^2=R+\lambda-|\ell|^2,\quad \lambda\in[-1/2,1/2].\tag 4.2$$

Writing $\partial$ for $\partial_\theta D$ and using the resolvent equation, we have 
$$\align \frac{\partial M_{ii}}{\partial\theta}=&[\hat V(E-M_0)^{-1}(\partial)(E-M_0)^{-1}\hat V](i, i)\tag 4.3\\
=&(\hat VD^{-1})\partial(D^{-1}\hat V)\tag 4.4\\
+&[(\hat VD^{-1})\partial(D^{-1}\hat V)^2+(\hat VD^{-1})^2\partial(D^{-1}\hat V)]\tag4.5\\
+&[(\hat VD^{-1})^2\partial(D^{-1}\hat V)^2+(\hat VD^{-1})\partial(D^{-1}\hat V)^3+(\hat VD^{-1})^3\partial(D^{-1}\hat V)]
\tag 4.6\\
+&[(\hat VD^{-1})\partial((E-M_0)^{-1}\hat V)(D^{-1}\hat V)^3+(\hat VD^{-1})^3(\hat V(E-M_0)^{-1})\partial(D^{-1}\hat V)\\
&\quad+(\hat VD^{-1})^2\partial(D^{-1}\hat V)^3+(\hat VD^{-1})^3\partial(D^{-1}\hat V)^2]\tag 4.7\\
+&[(\hat VD^{-1})^2\partial((E-M_0)^{-1}\hat V)(D^{-1}\hat V)^3+
(\hat VD^{-1})^3(\hat V(E-M_0)^{-1})\partial(D^{-1}\hat V)^2\\
&\quad+(\hat VD^{-1})^3\partial(D^{-1}\hat V)^3]\tag 4.8\\
+&[(\hat VD^{-1})^3\partial((E-M_0)^{-1}\hat V)(D^{-1}\hat V)^3+
(\hat VD^{-1})^3(\hat V(E-M_0)^{-1})\partial(D^{-1}\hat V)^3]\tag 4.9\\
+&[(\hat VD^{-1})^3(\hat V(E-M_0)^{-1})\partial((E-M_0)^{-1}\hat V)(D^{-1}\hat V)^3],\tag 4.10
\endalign$$
where it is understood that (4.4-4.10) pertain to the $(i,i)$ entry. It follows immediately from Lemma 3.2
and $\Vert\partial_\theta D\Vert=\Cal O(\sqrt R)$:
\proclaim{Lemma 4.1} $$\frac{1}{\sqrt R}|[(4.6)+\cdots+(4.10)]|\leq \Cal O(1)\sum_{(a_\ell, b_\ell)\in\text {supp }\hat V}
\frac{1}{\prod_{\ell=1}^4|ma_\ell+nb_\ell+\frac{a_\ell^2+b_\ell^2}{2}-\frac{\lambda}{2}|},\tag 4.11$$
where $(m,n)=i$. 
\endproclaim
\smallskip
The rest of this section is devoted to estimate the main terms (4.4, 4.5). Before that we first estimate (4.1),
which gives an appoximation to $\lambda$.
\proclaim{Lemma 4.2} $$|M_{ii}|\leq\Cal O(1)\big[\sum_{(a_\ell, b_\ell)\in\text {supp }\hat V}
\frac{1}{|ma_\ell+nb_\ell|}\big]^2,\tag 4.12$$
where $(m,n)=i$ and $|ma_\ell+nb_\ell|>K>c^2k^2+c\Vert V\Vert_\infty+1$ ($\Bbb N\ni c>8$).
\endproclaim
\demo{Proof}
Using the resolvent equation, we have 
$$M_{ii}=\hat VD^{-1}\hat V+\hat VD^{-1}\hat VD^{-1}\hat V+\hat VD^{-1}\hat V(E-M_0)^{-1}\hat VD^{-1}\hat V,\tag 4.13$$
where the right side only refers to the $(i,i)$ entry.

For any $(a,b)$, $(a',b')\in\Bbb Z^2\backslash\{0\}$, we say $(a,b)\sim (a', b')$ if $(a,b)=s(a',b')$ or $(a',b')=s(a,b)$,  
$s\in\Bbb Z \backslash\{0\}$. We call this equivalent class $\Cal C_{a,b}$ if $(a, b)$ is such that $a\geq 0$, $a+b\geq 0$
and $|(a,b)|\leq |(a',b')|$ for all $(a',b')$ such that $(a,b)\sim (a', b')$. We define the geometric support of $\hat V$ to be 
$$\text{gsupp }\hat V{\overset\text{def }\to =}\{(a,b)\in\Bbb Z^2\backslash\{0\}|\exists s\geq 1\text{ such that } (sa, sb)\in\text{ supp }\hat V\}.\tag 4.14$$
Assume $(a,b)\in\text{gsupp }\hat V$, we define
$$\align \lambda_{a, b}=&\sum_{(a'b')\in\Cal C_{a, b}}\hat V(a',b')D^{-1}(i+(a',b'), i+(a',b'))\hat V(-a',-b')\tag 4.15\\
=&\sum_{\Sb s\geq 1\\(sa, sb)\in\text{supp }\hat V\endSb }\hat V(sa,sb)D^{-1}(i+(sa,sb), i+(sa,sb))\hat V(-sa,-sb)\tag 4.16\\
&\qquad\qquad\quad\quad +\hat V(-sa,-sb)D^{-1}(i-(sa,sb), i-(sa,sb))\hat V(sa,sb)\tag 4.17\endalign$$
Using the above, we have 
$$\hat VD^{-1}\hat V=\sum_{(a,b)\in\text{gsupp }\hat V}\lambda_{a,b},\tag 4.18$$
where 
$$\aligned \lambda_{a,b}&=\sum_{s\geq 1}|\hat V(sa,sb)|^2[-\frac{1}{2sam+2sbn+s^2a^2+s^2b^2-\lambda}+\frac{1}{2sam+2sbn-(s^2a^2+s^2b^2-\lambda)}]\\
&=\sum_{s\geq 1}\frac{|\hat V(sa,sb)|^2\cdot (s^2a^2+s^2b^2-\lambda)}{2s^2}\cdot\frac{1}{(ma+nb)^2-(\frac{s^2a^2+s^2b^2-\lambda}{2s})^2},\endaligned\tag 4.19$$
$|a|$, $|b|\leq k$, $(a,b)\neq (0, 0)$ and $\lambda\in[-1/2, 1/2]$.
So $$\hat VD^{-1}\hat V=\Cal O(1)\sum_{(a_\ell, b_\ell)\in\text {gsupp }\hat V}
\frac{1}{(ma_\ell+nb_\ell)^2},\tag 4.20$$
if $|ma_\ell+nb_\ell|>K>c^2k^2+c\Vert V\Vert_\infty+1$ ($\Bbb N\ni c>8$).  

The second and third terms in the right side of (4.13) are bounded above by
$$\Cal O(1)\big[\sum_{(a_\ell, b_\ell)\in\text {supp }\hat V}\frac{1}{|ma_\ell+nb_\ell|}\big]^2\text{  and  }\Cal O(1/K)\big[\sum
_{(a_\ell, b_\ell)\in\text {supp }\hat V}\frac{1}{|ma_\ell+nb_\ell|}\big]^2.\tag 4.21$$
(4.20, 4.21) imply (4.12).\hfill []
\enddemo

Since $(4.4)=\frac{\partial}{\partial\theta} (4.18)$,  we take the derivative of (4.19) and have 
$$\aligned \frac{1}{\sqrt R}\frac{\partial}{\partial\theta}\lambda_{a, b}=-\sum_{1\leq s\leq k}&\frac{|\hat V(sa,sb)|^2(s^2a^2+s^2b^2-\lambda)}{s^2}\\
&(a\sin\theta-b\cos\theta)\cdot \frac{(ma+nb)}{[(ma+nb)^2-(\frac{s^2a^2+s^2b^2-\lambda}{2s})^2]^2},\endaligned\tag 4.22$$
$|a|$, $|b|\leq k$, $(a,b)\neq (0, 0)$. 

For a fixed $(a, b)\in\text{gsupp }\hat V$, (4.22) have a sign. More precisely, the vectors $(a,b)$ and $(-b, a)$ divide 
$\Bbb R^2$ into four quadrants. If $(m,n)$ is in the first and third, $\frac{1}{\sqrt R}\frac{\partial}{\partial\theta}\lambda_{a, b}>~0$,
otherwise $\frac{1}{\sqrt R}\frac{\partial}{\partial\theta}\lambda_{a, b}<0$. But the quadrants vary according to $(a,b)$
leading to cancellations in the sum:
$$\frac{1}{\sqrt R}\sum_{(a,b)\in\text{gsupp }\hat V}\frac{\partial}{\partial\theta}\lambda_{a, b}=\frac{(4.4)}{\sqrt R}.\tag 4.23$$
The following separation property of $\text{gsupp }\hat V$ plays an essential role in determining zeroes of (4.23).
\proclaim{Lemma 4.3} Let $(m,n)\in \tilde S$, the circle centered at $(0,0)$ of radius $\sqrt R$ in $\Bbb R^2$. If there exists 
$(a,b)\in\text{gsupp }\hat V$, $\Bbb N\ni |a|$, $|b|\leq k$, such that 
$$|ma+nb|<\epsilon\sqrt R,\quad\epsilon>0,\tag 4.24$$
then for all $$(a',b')\in\text{gsupp }\hat V\backslash\{(a,b)\}\tag 4.25$$
$$|ma'+nb'|>[\Cal O(1/k)-\epsilon]\sqrt R.\tag 4.26$$
\endproclaim
\demo{Proof}
$$ma+nb=(m,n)\cdot (a,b)=\sqrt R\cdot\sqrt{a^2+b^2}\cdot\cos\theta=\sqrt R\cdot\sqrt{a^2+b^2}\cdot\sin\phi,\tag 4.27$$
where $\theta$ is the angle between $(m,n)$ and $(a,b)$, $\phi=\pi/2-\theta$. 
$$\aligned ma'+nb'=&\sqrt R\cdot\sqrt {{a'}^2+{b'}^2}\cdot\cos\theta'=\sqrt R\cdot \sqrt {{a'}^2+{b'}^2}\cdot\sin\phi'\\
=&\sqrt R\cdot\sqrt {{a'}^2+{b'}^2}\cdot\sin(\phi'-\phi+\phi),\endaligned\tag 4.28$$
where $\theta'$ is the angle between $(m,n)$ and $(a',b')$, $\phi'=\pi/2-\theta'$. Since $\min |\phi'\pm\phi|=\Cal O(1/k)$
for $(a'b')$ satisfying (4.25), using (4.24) in (4.26, 4.27), we obtain (4.26).\hfill []
\enddemo
\noindent Using Lemma 4.3 in (4.12), we obtain that $\lambda\in[-\Cal O(1/K^2),\Cal O(1/K^2) ]$ in (4.22).

\proclaim{Lemma 4.4}
Let $(m,n)\in \tilde S$, the circle centered at $(0,0)$ of radius $\sqrt R$ in $\Bbb R^2$. If there exists $(a,b)\in\text{gsupp }\hat V$
such that  $$|ma+nb|<\epsilon\sqrt R,\quad(0<\epsilon<1/k^3),\tag 4.29$$
then $$\aligned \frac{|(4.4)+(4.5)|}{\sqrt R}&>\Cal O(1)\frac{1}{|ma+nb|^3}\\
&>\Cal O(\frac{1}{\epsilon^3})\frac{1}{R^{3/2}},\endaligned\tag 4.30$$
if 
$$\sum_{1\leq s\leq k}|\hat V(sa,sb)|^2(a^2+b^2)-\sum_{1\leq |s|, |s'|\leq k}\frac
{\hat V(sa, sb)\hat V(s'a,s'b)\hat V((s-s')a,(s-s')b))}{4s^2s'}\neq 0.\tag 4.31$$
\endproclaim
\demo{Proof} We first assume $a\geq 0$ and $b\geq 0$. Write $m=\sqrt R\cos\theta$ and $n=\sqrt R\sin\theta$.  
We distinguish in (4.5) the terms only involve $\hat V(sa, sb)$ with $(a,b)\in\text{gsupp }\hat V$
satisfying (4.29),
$1\leq|s|\leq k$ and call the sum $\mu_{a,b}$. We have 
$$\aligned \mu_{a,b}=\sum_{1\leq |s|, |s'|\leq k}&\frac
{\hat V(sa, sb)\hat V(s'a,s'b)\hat V((s-s')a,(s-s')b))}{4s^2s'}\\
&(a\sin\theta-b\cos\theta)\cdot\frac{1}{(ma+nb)^3}+\Cal O(\frac{1}{(ma+nb)^5}).\endaligned\tag 4.32$$
There are three cases:
$a>0$, $b>0$; $a=0$, $b>0$ and $a>0$, $b=0$.
\item{(i)}  $a>0$, $b>0$, 
(4.29) implies $\cos\theta\sin\theta\leq 0$. Otherwise $|ma+nb|>|m|+|n|>\sqrt R$. So 
$|a\sin\theta-b\cos\theta|\geq 1$ and
$$\frac{1}{\sqrt R}|\frac{\partial}{\partial\theta}\lambda_{a, b}+\mu_{a,b}|>\Cal O(1)\cdot \frac{1}{|ma+nb|^3}$$
from (4.22, 4.32), where we used (4.31, 4.12). 

\item{(ii)} $a=0$, $b>0$, (4.29) reduces to 
$$|nb|<\epsilon\sqrt R.$$
So $$|\sin\theta|<\Cal O(\epsilon),\quad |\cos\theta|>1-\Cal O(\epsilon).\tag 4.33$$
Using (4.29, 4.33) in (4.22, 4.32), we obtain 
$$\frac{1}{\sqrt R}|\frac{\partial}{\partial\theta}\lambda_{0, b}+\mu_{0,b}|>\Cal O(1)\cdot \frac{1}{|ma+nb|^3}\tag 4.34$$

\item{(ii)} $a=0$, $b>0$, similarly,
$$\frac{1}{\sqrt R}|\frac{\partial}{\partial\theta}\lambda_{a, 0}+\mu_{a, 0}|>\Cal O(1)\cdot \frac{1}{|ma+nb|^3}\tag 4.35$$

Clearly, same estimates hold for $a\geq 0$ and $b\leq 0$. So we have 
$$\aligned \frac{1}{\sqrt R}|\frac{\partial}{\partial\theta}\lambda_{a, b}+\mu_{a, b}|&>\Cal O(1)\frac{1}{|ma+nb|^3}\\
&>\Cal O(\frac{1}{\epsilon^3})\cdot\frac{1}{R^{3/2}}>\Cal O(\frac{1}{\epsilon^3})\cdot\sum_{\Sb(a'b')\neq(a,b)\\(a'b')\in\text{gsupp }\hat V\endSb}\frac{1}{|ma'+nb'|^3}\\
&>\Cal O(\frac{1}{\epsilon^3})\frac{1}{\sqrt R}\sum_{\Sb(a'b')\neq(a,b)\\(a'b')\in\text{gsupp }\hat V\endSb}
|\frac{\partial}{\partial\theta}\lambda_{a', b'}|,
\endaligned\tag 4.36$$
where we used (4.29, 4.32-4.35, 4.22) and Lemma 4.3.

So 
$$\aligned \frac{|(4.4)+\mu_{a, b}|}{\delta^2\sqrt R}&>\frac{1}{\sqrt R}\big[|\frac{\partial}{\partial\theta}\lambda_{a, b}+\mu_{a, b}|-\sum_{\Sb(a'b')\neq(a,b)\\(a'b')\in\text{gsupp }\hat V\endSb}|\frac{\partial}{\partial\theta}\lambda_{a', b'}|\big]\\
&>(1-\Cal O(\epsilon^3k^2))\frac{1}{\sqrt R}|\frac{\partial}{\partial\theta}\lambda_{a, b}+\mu_{a, b}|\\
&>\Cal O(1) \frac{1}{|ma+nb|^3}\\
&>\Cal O(\frac{1}{\epsilon^3})\frac{1}{R^{3/2}}.\endaligned\tag 4.37$$

Since each term in (4.5) aside from $\mu_{a, b}$ is third order in $D^{-1}$ involving at least one $(a',b')\in\text{gsupp }\hat V$, $(a',b')\neq (a,b)$, (4.24, 4.26) imply
$$|\frac{(4.5)-\mu_{a, b}}{\sqrt R}|<\Cal O(k^2)\frac{1}{|ma+nb|^2|ma'+nb'|}<\Cal O(\epsilon k^3) \frac{1}{|ma+nb|^3}.$$
Combining with (4.37), this proves the lemma for $0<\epsilon<1/k^3$. \hfill []
\enddemo

Combining Lemme 4.1 and 4.4, we have 
\proclaim{Proposition 4.5} Let $(m,n)\in \tilde S\cap\Pi'$. If there exists $(a,b)\neq (0,0)$, $(a,b)\in\text{gsupp }\hat V$
such that 
$$|ma+nb|<\epsilon\sqrt R,\quad(0<\epsilon<1/k^3),\tag 4.38$$
then $$\frac{1}{\sqrt R}| \frac{\partial M_{ii}}{\partial\theta}|>\frac{\Cal O(1)}{|ma+nb|^3}
>\frac{\Cal O(1)}{R^{3/2}},\quad (|ma+nb|>K>c^2k^2+c\Vert V\Vert_\infty+1, \Bbb N\ni c>8),\tag 4.39$$
provided (4.31) holds.
\endproclaim
\demo{Proof} This follows immediately from (4.30, 4.11, 4.3).\hfill []
\enddemo
\noindent{\it Polynomial approximation.}

It follows from Proposition 4.5 that in order to control the zeroes of $\frac{\partial M_{ii}}{\partial\theta}$, we only 
need to restrict to $(m,n)$ such that 
$$|ma+nb|\geq\epsilon\sqrt R\tag 4.40$$
for all $(a,b)\in\text{gsupp }\hat V$. More precisely, in view of Lemma 4.1, we want to exclude $\theta$ satisfying (4.40)
such that 
$$\frac{1}{\sqrt R}\big|\sum_{(a,b)\in\text{gsupp }\hat V} \frac{\partial}{\partial\theta}\lambda_{a, b}+(4.5)\big|
=|\frac{(4.4)+(4.5)}{\sqrt R}|\leq\frac{\Cal O(1)}{R^2}.\tag 4.41$$

Let 
$$\align
 \Lambda{\overset\text{def }\to =}&\frac{1}{\sqrt R}\big(\sum_{(a,b)\in\text{gsupp }\hat V} \frac{\partial}{\partial\theta}\lambda_{a, b}+(4.5)\big),\tag 4.42\\
v_{a,b}{\overset\text{def }\to =}&\sum_{1\leq s\leq k}|\hat V(sa,sb)|^2(a^2+b^2),\quad (a,b)\in\text{gsupp }\hat V,
\tag 4.43\\
g_{a,b;c,d}{\overset\text{def }\to =}&\sum_{\Sb-k\leq s, s'\leq k\\s,s'\neq 0\endSb}\frac{\hat V(sa,sb)\hat V(s'c-sa,s'd-sb)\hat V(s'c,s'd)}{4s^2s'}, 
\quad (a,b), (c,d) \in\text{gsupp }\hat V.\tag 4.44\endalign$$
Assume $(m,n)$ such that (4.40) hold for all $(a,b)\in\text{gsupp }\hat V$ and $\lambda=\Cal O(1/R)$, then 
$$\align\Lambda=&\frac{1}{R^{3/2}}\big(\sum_{(a,b)\in\text{gsupp }\hat V} v_{a,b}\frac{a\sin\theta-b\cos\theta}{(a\cos\theta+b\sin\theta)^3}\\
&\qquad\quad+\sum_{(a,b), (c,d) \in\text{gsupp }\hat V}g_{a,b;c,d}\frac{a\sin\theta-b\cos\theta}{(a\cos\theta+b\sin\theta)^2(c\cos\theta+d\sin\theta)}\big)\tag 4.45\\
&\qquad\qquad+\Cal O(\frac{1}{R^2})+\Cal O(\frac{1}{R^{5/2}})\\
{\overset\text{def }\to =}&\frac{1}{R^{3/2}}\Lambda_1+\Cal O(\frac{1}{R^2})+\Cal O(\frac{1}{R^{5/2}}),\tag 4.46
\endalign$$
where $v_{a,b}$ and $g_{a,b;c,d}$ as in (4.43).

Let $\nu=|\text{gsupp }\hat V|\leq\Cal O(k^2)$, (4.40) define $2\nu$ arcs $\Gamma'$ of the circle $\tilde S$. Let $(a,b)_{\perp}$ be the
ray perpendicular to $(a,b)\in\text{gsupp }\hat V: (a,b)_{\perp}\cdot(a,b)=0$. Then for all $\Gamma'$, 
$\Gamma'\cap (a,b)_{\perp}=\emptyset$, for all $(a,b)\in\text{gsupp }\hat V$, and $\Lambda_1$ is well defined on $\Gamma'$.

Let $$\aligned &x=\tan\theta,\quad \text { when } |\tan\theta|\leq 1,\\
&x=\coth \theta,\quad \text {when } |\coth\theta|<1.\endaligned\tag 4.47$$
Rewrite $\Lambda_1$ in terms of $x$ and call the resulting function $f$. We have 
$$\align f=(1+x^2)&\big(\sum_{(a,b)\in\text{gsupp }\hat V} v_{a,b}\frac{ax-b}{(a+bx)^3}
+\sum_{(a,b), (c,d) \in\text{gsupp }\hat V}g_{a,b;c,d}\frac{ax-b}{(a+bx)^2(c+dx)}\big),\\
&\qquad\qquad\qquad\qquad\qquad  |\tan\theta|\leq 1, x=\tan\theta, |x|\leq 1,\tag 4.48\\
f=(1+x^2)&\big(\sum_{(a,b)\in\text{gsupp }\hat V} v_{a,b}\frac{a-bx}{(ax+b)^3}
+\sum_{(a,b), (c,d) \in\text{gsupp }\hat V}g_{a,b;c,d}\frac{a-bx}{(ax+b)^2(cx+d)}\big),\\
&\qquad\qquad\qquad\qquad\qquad  |\coth\theta|\leq 1, x=\coth\theta, |x|<1.\tag 4.49
\endalign$$
Both $f$ are rational functions and can be written as 
$$f=\frac{P_1}{P_2},\tag 4.50$$
where $P_1$ and $P_2$ are polynomials in $x$ of degrees at most $3(\nu^2+\nu)<4\nu^2$ and 
$$0<|P_2|<\Cal O(1)\tag 4.51$$
on arcs $\Gamma'$, defined above (4.47). Moreover $P_1$ is a polynomial whose coefficients {\it only}
depends on $\hat V$ and $\text{supp }\hat V$ in view of (4.43, 4.48, 4.49). It is of the form 
$$P_1=A_px^p+A_{p-1}x^{p-1}+\cdots+A_0,\quad A_p\neq 0, 0<p< 4\nu^2,\tag 4.52$$
and $$A_j=A_j(\hat V, \text{supp }\hat V).\tag 4.53$$

From (4.51), the set 
$$I{\overset\text{def }\to =}\{x||f(x)|<\frac{1}{\sqrt R}\}\subseteq I_1{\overset\text{def }\to =}\{x||P_1(x)|<\frac{\Cal O(1)}{\sqrt R}\}.
\tag 4.54$$
To bound the measure of $I_1$, we use the resultant. From (4.52),
$$P'_1=pA_px^{p-1}+(p-1)A_{p-1}x^{p-2}+\cdots+A_1.\tag 4.55$$
By definition, 
$$\text{Resultant }(P_1, P'_1)=\det\pmatrix A_p&A_{p-1}&A_{p-2}&\cdots&A_1&A_0&\cdots&\cdots\\
0&A_p&A_{p-1}&\cdots&\cdots&A_1&A_0&\cdots\\
\vdots&{}&\vdots&{}&A_p&A_{p-1}&\cdots&A_0\\
pA_p&(p-1)A_{p-1}&(p-2)A_{p-2}&\cdots&A_1&0&\cdots&\cdots\\
0&\cdots&0&\cdots&\cdots&pA_p&\cdots&A_1
\endpmatrix, \tag 4.56$$
and let  $D(\hat V)$ denote the above resultant. $\quad\qquad\qquad\qquad\quad\quad\qquad\qquad\quad\quad\quad\,$(4.57)

If  $D(\hat V)\neq 0$, $P_1$ and $P_1'$ have no 
common roots. 
Let $\Gamma\subset \tilde S$ be the largest set such that on $\Gamma$, (4.40) hold for all $(a,b)\in\text{gsupp }\hat V$.
Since $R$ is fixed, we also use $\Gamma$ to denote the corresponding set of angles $\theta\in[0,2\pi)$. 
\proclaim{Lemma 4.6} Assume $V$ is such that $D(\hat V)\neq 0$ for both $P_1$ defined from (4.48-4.50), 
and $\lambda=\Cal O(1/R)$, then 
$$\text{mes }\{\theta\in\Gamma||\Lambda(\theta)|\leq \frac{\gamma}{R^2}\}\leq \frac {\gamma C_V}{\sqrt R},\tag 4.58$$
where $\Lambda$ is as defined in (4.42). Moreover the set in (4.58) has at most $\Cal O(k^4)$ connected components.
\endproclaim
\demo{Proof}
$P_1$ is of degree at most $4\nu^2$, with $\nu=|\text{gsupp }\hat V|=\Cal O(k^2)$. So $P_1$ has at most
$4\nu^2$ zeroes. Since $D(\hat V)\neq 0$, 
$$\text{min }\{|P'_1(x)|P_1(x)=0\}>\frac{1}{C_{V}}>0.\tag 4.59$$
So $$\text{mes }\{x|P_1(x)|\leq \frac{\gamma}{\sqrt R}\}\leq \frac {\gamma C_{V}}{\sqrt R},\tag 4.60$$
(4.59, 4.48-4.50, 4.45) and the fact that 
$$d\theta=\pm\frac{1}{1+x^2}dx$$
imply (4.58).\hfill []
\enddemo
\bigskip
\head{\bf 5. Proof of the Theorem}\endhead
Assume $V$ is a generic trigonometric polynomial of degree $k$ satisfying the genericity conditions (i, ii) in sect. 1,
so that Lemmas 4.4 and 4.6 are available. Let $\tilde S$ be the
circle over $\Bbb R^2$ of radius $\sqrt R$, $R\in\Bbb N$ as before. Take $c=9$ in Lemma 2.2
and define the geometric singular set 
$$\Theta_g{\overset\text{def }\to =}\{\theta\in[0,2\pi)||\alpha\cos\theta+\beta\sin\theta|\leq\frac{K}{\sqrt R}\,\text{ for some }
(\alpha,\beta)\in[-9k, 9k]^2\backslash\{0\},\tag 5.1$$
where $$K>c^2k^2+c\Vert V\Vert_\infty+1=81 k^2+9\Vert V\Vert_\infty+1.\tag 5.2$$
$\Theta_g$ has at most $\Cal O(k^2)$ connected components and 
$$\text{mes }\Theta_g=\frac{\Cal O(1)}{\sqrt R}\,\text{ on }[0,2\pi).\tag 5.3$$
We also use $\Theta_g$ to denote the corresponding arcs of $\tilde S$, As before, let 
$$\Gamma=\{(m,n)\in \tilde S||ma+nb|\geq \epsilon\sqrt R,\quad \forall (a,b)\in\text{gsupp }\hat V\},\tag 5.4$$
where $0<\epsilon<1/k^3$.

Assume $\lambda=\Cal O(1/R)$ and let $\Theta_a\subset\Gamma$ be the
algebraic singular set defined in (4.58) with $\gamma>k^{10}/\epsilon^4>k^{22}$ in view of (4.11) and (4.40), 
$$\text{mes }\Theta_a\leq \frac{\Cal O(1)}{\sqrt R}\,\text{ on }[0,2\pi),\tag 5.5$$
where $\Cal O(1)=\gamma C_V$.
Define $$\Theta=\Theta_g\cup\Theta_a\tag 5.6$$
and let $\Theta$ also denote the corresponding set on $\tilde S$. $\Theta$ has at most $\Cal O(k^4)$ connected components,
$$\text{mes }\Theta=\Cal O(1)\,\text{ on } \tilde S.\tag 5.7$$
Let $p\in\Pi$, $\Pi$ as defined in Lemma 2.1 and $\bar S=\tilde S\cap\Bbb Z^2$. Assume
$$p\cap\Theta=\emptyset,\quad p\cap \bar S\neq\emptyset.\tag 5.8$$
Define 
$$\Cal M{\overset\text{def }\to =}\Cal M_p\tag 5.9$$
as in (3.2), $M'$ as in (3.3-3.6) and $M_{ii}$ as in (4.1) first for $i\in p\cap \bar S$, then for $i\in \tilde S\backslash\Theta_g$.
\proclaim{Lemma 5.1} Assume $$\lambda'+M_{ii}(\lambda')=0,\tag 5.10$$ 
where $\lambda'=E-|i|^2=E-R$.
Then on $\Gamma$ defined in (5.4),
$$\lambda'=\Cal O(\frac{1}{R});\tag 5.11$$
and on each connected component of 
$$S''{\overset\text{def }\to =}\tilde S\backslash\Theta\tag 5.12$$
either 
$$\align&\frac{1}{\sqrt R}\frac{d\lambda'}{d\theta}\geq \Cal O(1)\sum_{(a_\ell, b_\ell)\in\text {gsupp }\hat V}
\frac{1}{\prod_{\ell=1}^4|ma_\ell+nb_\ell|},\tag 5.13\\
\text{or }&\frac{1}{\sqrt R}\frac{d\lambda'}{d\theta}\leq -\Cal O(1)\sum_{(a_\ell, b_\ell)\in\text {gsupp }\hat V}
\frac{1}{\prod_{\ell=1}^4|ma_\ell+nb_\ell|},\tag 5.14\endalign$$
where $i=(m,n)=\sqrt R (\cos\theta,\sin\theta)$.
\endproclaim
\demo{Proof}
(5.11) follows from Lemma 4.2.  Using (4.30) or (4.58) in (4.4, 4.5), Lemma 4.1 in (4.6-4.10), we obtain
$$\frac{1}{\sqrt R}\big|\frac{\partial M_{ii}}{\partial\theta}\|\geq\Cal O(1)\sum_{(a_\ell, b_\ell)\in\text {gsupp }\hat V}
\frac{1}{\prod_{\ell=1}^4|ma_\ell+nb_\ell|},\tag 5.15$$
$\Bbb R^2\ni i=(m,n)$ on $\tilde S\backslash\Theta$. Here we also used (5.11), when $\theta\in\Gamma$. Moreover it
is sign definite on each connected component of $\tilde S\backslash\Theta=S''$. From (5.10)
$$-\frac{d\lambda'}{d\theta}=\frac{\partial M_{ii}}{\partial\theta}+\frac{\partial M_{ii}}{\partial\lambda'}\cdot \frac{d\lambda'}{d\theta}.\tag 5.16$$
So $$\frac{d\lambda'}{d\theta}=(-1+\Cal O(1/K^2))\cdot \frac{\partial M_{ii}}{\partial\theta},\tag 5.17$$
where we used (3.33). Using (5.15), we obtain the lemma. $\hfill$ []
\enddemo
\proclaim{Proposition 5.2}
The set $S''=\tilde S\backslash\Theta$ has at most $\Cal O(k^4)$ connected components. Let $\Gamma''\subset S''$ be 
a connected component. Assume $p$, $p'\in\Pi$ and $p$, $p'\cap \bar S\neq \emptyset$ be 
such that $p$, $p'\cap \tilde S\subset\Gamma''$.  Let $i\in \{p, p'\cap \bar S\}$, $\lambda_i=E-|i|^2=E-R\in[-1/2,1/2]$ 
be such that 
$$0\in\sigma(\Cal M(\lambda_i)),\tag 5.18$$
where $\Cal M=\Cal M_p$ or $\Cal M_{p'}$. Let $\lambda_i'=E'-|i|^2=E'-R\in[-1/2,1/2]$ be such that 
$$0=\lambda_i'+M_{ii}(\lambda'_i).\tag 5.19$$
Then $$|\lambda_i-\lambda'_i|=\Cal O(1)\sum_{(a_\ell, b_\ell)\in\text {gsupp }\hat V}
\frac{1}{\prod_{\ell=1}^8|ma_\ell+nb_\ell|},\quad (m,n)=i,\tag 5.20$$
and $$|\lambda_i-\lambda_{j'}|>\sup_{(m,n)=i,j'}\sum_{(a_\ell, b_\ell)\in\text {gsupp }\hat V}
\frac{\Cal O(1)}{\prod_{\ell=1}^4|ma_\ell+nb_\ell|},\quad i,j'\in \{p, p'\cap \bar S\}, i\neq j'.\tag 5.21$$
\endproclaim
\demo{Proof} The number of connected components follow from the definition of $\Theta$ in (5.1, 4.58, 5.6).
Assume $\Cal M=\Cal M_p$ has rank $\geq 2$. Write the right side of (3.29) as $\Cal O_i$. We have 
$$\aligned F(E)&=\det (\Cal M(E))\\
&=\prod_j(E-|j|^2+M_{jj}(E)+\Cal O_j)+\Cal O(\sum
\prod\Cal M_{ij}),\endaligned\tag 5.22$$
where the second product contains at least two off diagonal elements.
So $$\aligned F'(E)=&\sum_j\prod_{j'\neq j}(E-|j'|^2+M_{jj}(E)+\Cal O_{j'}) (1+\frac{\partial M_{jj}}{\partial E}(E)+\Cal O_j)\\
&+\Cal O_{ij},
\endaligned\tag 5.23$$
where we used $\Cal O_{ij}$ to denote the right side of (3.31) and analyticity in $E$ to reach (5.23).

Let $E=|i|^2+\lambda'_i$ and write $F(\lambda_i')$, $F'(\lambda_i')$ for $F(|i|^2+\lambda_i')$, $F'(|i|^2+\lambda_i')$
respectively. Since $\frac{\partial M_{jj}}{\partial E}=\Cal O(1/K^2)$, (5.23) gives 
$$|F'(\lambda_i')|\geq(1-\Cal O(1/K^2))\cases |\lambda_i'-\lambda_{j'}'|-\max_{ij}\Cal O_{ij},\quad \text{if } \exists j'\in p,\, 
j'\neq i,\,|j'|^2=R,\\
\Cal O(1),\qquad\qquad\qquad\qquad \text{otherwise},\endcases
\tag 5.24$$
where we also used (4.12).
Using (5.13) or (5.14), we have for any $\tilde j\in\Gamma''$, $|\tilde j-i|\geq 1$ and
$$\aligned |\lambda_i'-\lambda_{\tilde j}'|&=|\int_{\tilde j}^i (\frac{d\lambda'}{d\theta})d\theta|\geq
\max_{(m,n)=i,\tilde j}\sum_{(a_\ell, b_\ell)\in\text {gsupp }\hat V}
\frac{1}{\prod_{\ell=1}^4|ma_\ell+nb_\ell|}\\
&\quad\gg\max_{i\tilde j}\Cal O_{i\tilde j}.\endaligned \tag 5.25$$
So $$|F'(\lambda_i')|\geq(1-\Cal O(1/K^2))\cases |\lambda_i'-\lambda_{j'}'|,\quad \text{if } \exists j'\in p,\, 
j'\neq i,\,|j'|^2=R,\\
\Cal O(1),\qquad\qquad\qquad\qquad \text{otherwise}.\endcases
\tag 5.26$$

From (5.22), 
$$|F(\lambda_i')|\leq \Cal O_i\cases |\lambda_i'-\lambda_{j'}'|,\quad \text{if } \exists j'\in p,\, 
j'\neq i,\,|j'|^2=R,\\
\Cal O(1),\qquad\qquad\qquad\qquad \text{otherwise}.\endcases
\tag 5.27$$
Let $0<a\ll 1$.
$$F(\lambda_i'\pm a)=F(\lambda_i')\pm a F'(\lambda_i')+\Cal O(a^2)=F'(\lambda_i')(\frac{F(\lambda_i')}{F'(\lambda_i')}\pm a)
+\Cal O(a^2).\tag 5.28$$
Since $$|\frac{F(\lambda_i')}{F'(\lambda_i')}|\leq\Cal O_i{\overset\text{def }\to =}\Cal O(1)\sum_{(a_\ell, b_\ell)\in\text {gsupp }\hat V}\frac{1}{\prod_{\ell=1}^8|ma_\ell+nb_\ell|}\tag 5.29$$
from (3.29), for $a>10\Cal O_i$,
$F(\lambda'+a)$ and $F(\lambda'-a)$ have opposite signs. Since $F$ is analytic, this implies 
$F(\lambda_i)=0$ for some 
$$\lambda_i\in \lambda'_i+(-11\Cal O_i, 11\Cal O_i),\tag 5.30$$
which proves (5.20). Since similar statements hold for $\lambda_{j'}$, we obtain
$$ |\lambda_i-\lambda_{j'}|>\frac{1}{2} |\lambda_i'-\lambda_{j'}'|,\tag 5.31$$
implying (5.21) by using (5.25).

Clearly simpler arguments apply when $\Cal M$ is a scalar as $F'(\lambda_i')=(1-\Cal O(1/K^2))>1/2$ and $F'(\lambda_i')=
\Cal O_i$. Combining the two cases, we obtain the proposition.\hfill []
\enddemo
\smallskip
\noindent{\it Proof of the Theorem.}

$$\sigma(H)=\sigma(\hat H)\subseteq\cup_{R\in\Bbb Z}[R-1/2, R+1/2].\tag 5.32$$
In the following lines, we go back to the convention of writing $H$ for $\hat H$. 
Since $R\geq -\Vert V\Vert_\infty$ and for $-\Vert V\Vert_\infty<R\leq 0$, (1,2, 1.3) are obvious, we only need to be concerned with $R\in\Bbb N$.
From Proposition 5.2, given 
$E\in[R-1/2, R+1/2]$, $R\in\Bbb N$, there exist at most $\Cal O(k^4)$ $p\in\Pi$, such that 
$$\text{dist }(E,\sigma(H_p))\leq o(\frac{1}{R^2}).\tag 5.33$$ 
First recall  $$\aligned &S'=\{(x,y)\in\Bbb R^2|\,|x^2+y^2-R|\leq\Vert V\Vert_\infty+1\},\quad S=S'\cap\Bbb Z^2,\\
&\tilde S=\{(x,y)\in\Bbb R^2|\,x^2+y^2=R\},\quad \bar S=\tilde S\cap\Bbb Z^2\endaligned$$
and $\Theta$ as defined in (5.1, 4.58, 5.6).
This is because
\item{$\bullet$} if $p\cap\Theta=\emptyset$ and $p\cap\bar S=\emptyset$, then
$$\text{dist } (\sigma(H_p), [R-1/2, R+1/2])\geq\Cal O(1),$$
from Proposition 3.1, Lemma 4.2 and analyticity of $\Cal M_p$ in $E$,
\item{$\bullet$} $S''=\tilde S\backslash\Theta$ has at most $\Cal O(k^4)$ connected components $\Gamma''$.
On each $\Gamma''$, (5.21) hold,
\item{$\bullet$} $\Theta$ has at most $\Cal O(k^4)$ connected components and 
$\text{mes }\Theta=\Cal O(1)$ on $\tilde S$,
\item{$\bullet$} for all $p\in\Pi$, $|p\cap \tilde S|=\Cal O(R^{1/6})|$, $|p\cap S|\leq \kappa$.

Assume $p\cap S\neq\emptyset$. Since $$\text{dist }(\{p\cap S\},\partial p)=\Cal O(R^{1/6})\tag 5.34$$
by construction (Lemma 2.1), this implies 
$$\text{dist }(\sigma(H_p), \sigma(H))=\Cal O(e^{-R^{1/6}}).\tag 5.35$$
This is because if $\hat\phi$ is an eigenfunction of $H_p$ with eigenvalue $E$, then $(H-E)\hat\phi=\Cal O(e^{-R^{1/6}})$,
which implies (5.35).

In fact more generally, for all $\Lambda$ such that either $\Lambda\subset p$ or $\Lambda\supseteq p$: 
$$\text{dist }(\sigma(H_p), \sigma(H_\Lambda))=\Cal O(e^{-\text{min }(d_1,d_2)}).\tag 5.36$$
where $$\align &d_1=\text{dist }(\{p\cap S\},\partial\Lambda\},\tag 5.37\\
&d_2=\text{dist }(\{p\cap S\},\partial p\}.\tag 5.38\endalign$$

Let $E\in\sigma(H)$, since each $H_p$ has at most $\kappa$ eigenvalues in $[R-1/2, R+1/2]$, (5.33, 5.35) 
give that $\sigma(H)$ is of multiplicity at most $\Cal O(k^4)$. 

To prove localization of the Fourier series $\hat\phi$ of the eigenfunction $\phi$, we proceed as follows. Let $p\in\Pi$
be such that $\text{dist }(E,\sigma(H_p))\leq o(\frac{1}{R^2})$.
Let $\Cal S$ be this set of singular $p$. From the argument above, there are only $\Cal O(k^4)$ such $p$. Let 
$$\Cal R=\{(m,n)\in S\cap\Cal S\}.\tag 5.39$$
Then $|\Cal R|=\Cal O(k^4)$, since $|p\cap S|\leq \kappa$.
(Note that $|\Cal R|\geq 1$ from (5.35). So the following construction is not empty.)

Since $\hat\phi\in\ell^2$, we may assume $\Vert\hat\phi\Vert_{\infty}\leq 1$ by normalization: $\Vert\hat\phi\Vert_2=1$.
So $$|\hat\phi(j)|\leq 1\,\text{ for }j\in\Cal R.\tag 5.40$$
To prove decay of $\hat\phi(j)$ for $j\notin\Cal R$, we let
$i_1\in\Cal R$ be such that 
$$|i_1-j|=\text{min }_{i\in\Cal R}|i-j|.\tag 5.41$$
(If there are two sites which are minimal, choose one and name it $i_1$.) Let $\Lambda$ be a square of size $\Cal O(|j-i_1|)$
such that $i_1\in\Lambda$, $j\in\Lambda$ and  
$$\text{dist }(j,\partial\Lambda)=2|j-i_1|.\tag 5.42$$
Let $$\tilde\Lambda=\Lambda\backslash\Cal R.\tag 5.43$$
\item{(i)} If $|j-i_1|\leq R^{1/7}$, then
$$\Vert (H_{\tilde\Lambda}-E)^{-1}\Vert\leq \Cal O(1),\tag 5.44$$
since $\tilde\Lambda\cap S=\emptyset$.
\item{(ii)} Otherwise $$\Vert (H_{\tilde\Lambda}-E)^{-1}\Vert\leq \Cal O(R^2)\tag 5.45$$
from (5.21, 5.36).

Define $$\Cal V=H-(H_{\tilde\Lambda}\oplus H_{\Bbb Z^2\backslash\tilde\Lambda}).\tag 5.46$$
Since $$(H-E)\hat\phi=0,\tag 5.47$$
we have $$\Pi_{\tilde\Lambda}\hat\phi=\Pi_{\tilde\Lambda}(H_{\tilde\Lambda}-E)^{-1}\Cal V\hat\phi.\tag 5.48$$

\item{(i)} $$|\hat\phi(j)|\leq C \sum_{j_\ell\in\Cal R\cap\Lambda}e^{-|j-j_\ell|}\tag 5.49$$
follows from Neumann series expansion about the diagonal.
\item{(ii)} Let $\Cal R'=\tilde\Lambda\cap S$. For $i'\in\Cal R'$, let $\Lambda'$ be the square centered at $i'$ of size 
$L'=(\log R)^2$. There are two possibilities: 
$\text{dist }(\{\Lambda'\cap S\},\,\partial\Lambda')=\Cal O((\log R)^2)$ or $\text{dist }(\{\Lambda'\cap S\},\,\partial\Lambda')<\Cal O((\log R)^2)$. In the
latter case, let $L''=100 L'$ and $\Lambda''$ be the square 
centered at $i'$ of size $L''$. By construction
$$\text{dist }(\{\Lambda''\cap S\},\partial\Lambda'')=\Cal O((\log R)^2),\tag 5.50$$
this is because from Lemma 2.1, for a given integer in $S$ there is at most $1$ other integer in $S$ which is at distance $\asymp\Cal O((\log R)^2)$ apart. Rename $\Lambda''$ as $\Lambda'$.

We have from (5.36, 5.21)
$$\text{dist }(E,\sigma(H_{\Lambda'}))\geq\Cal O(\frac{1}{R^2})\tag 5.51$$
and moreover 
$$|(H_{\Lambda'}-E)^{-1}(x,y)|\leq e^{-|x-y|}\tag 5.52$$
for $|x-y|\geq L'/10$ by using Neumann series, (5.51) and the fact that  
$|x_1-x_2|\leq L'/100$ for all $x_1$, $x_2\in\{\Lambda'\cap S\}$.
Clearly (5.51, 5.52) hold for all $\Lambda'$ of size $\Cal O((\log R)^2)$, $\Lambda'\subset\tilde\Lambda$, $\Lambda'\cap S=\emptyset$. 

Expanding $(H_{\tilde\Lambda}-E)^{-1}$ repeatedly in $(H_{\Lambda'}-E)^{-1}$ using the resolvent equation:
$$(H_{\tilde\Lambda}-E)^{-1}=(H_{\Lambda'}-E)^{-1}\tilde\Gamma(H_{\tilde\Lambda}-E)^{-1},$$
where $\tilde\Gamma{\overset\text{def }\to =}H_{\tilde\Lambda}-(H_{\Lambda'}\oplus H_{\tilde\Lambda\backslash\Lambda'})$,
(5.51, 5.52, 5.45) give 
$$|\hat\phi(j)|\leq C \sum_{j_\ell\in\Cal R\cap\Lambda}e^{-|j-j_\ell|}.$$
Combining cases (i,ii), we obtain (1.2, 1.3) and hence the Theorem.\hfill []

\enddocument
\end